\documentclass{amsart}
\usepackage{amssymb} 
\usepackage{color}

\usepackage{graphicx}

\usepackage[all]{xy}
 \newtheorem{theorem}{Theorem}[section]

 \newtheorem{proposition}[theorem]{Proposition}
 \newtheorem{lemma}[theorem]{Lemma}
 \newtheorem{corollary}[theorem]{Corollary}
 \newtheorem{remark}[theorem]{Remark}
 \newtheorem{definition}[theorem]{Definition}
 
 \newtheorem{example}[theorem]{Example}

 \newtheorem{defprop}[theorem]{Definition and Proposition}

 \numberwithin{equation}{section}
\allowdisplaybreaks[1]

\begin{document}

\title
[Applications of Hochschild cohomology to the moduli of subalgebras
]
{
Applications of Hochschild cohomology to the moduli of subalgebras of the full matrix ring 
}
\author{Kazunori Nakamoto and Takeshi Torii}
\address{%Center for Life Science Research,
%Interdisciplinary Graduate School of Medicine and Engineering,
Center for Medical Education and Sciences,
Faculty of Medicine, 
University of Yamanashi,
Yamanashi 409--3898, Japan}
\email{nakamoto@yamanashi.ac.jp}

\address{Department of Mathematics, 
%Faculty of Science, 
Okayama University,
Okayama 700--8530, Japan}
\email{torii@math.okayama-u.ac.jp}
\thanks{The first author was partially supported by 
JSPS KAKENHI Grant Numbers JP15K04814, JP20K03509. 
The second author was partially supported by
JSPS KAKENHI Grant Numbers JP22540087, JP25400092, JP17K05253.}

\subjclass[2010]{Primary 16E40; Secondary 14D22, 16S50, 16S80.}

\keywords{Hochschild cohomology, Subalgebra, Matrix ring, Moduli of molds} 

\date{June 14, 2020\ ({\tt version~1.0.0})}

\begin{abstract}
Let ${\rm Mold}_{n, d}$ be the moduli of rank $d$ subalgebras of ${\rm M}_n$ over ${\Bbb Z}$. For $x \in {\rm Mold}_{n, d}$, let ${\mathcal A}(x) \subseteq {\rm M}_n(k(x))$ be the subalgebra of ${\rm M}_n$ corresponding to $x$, where $k(x)$ is the residue field of $x$.  In this article, we apply Hochschild cohomology to ${\rm Mold}_{n, d}$.  
The dimension of the tangent space $T_{{\rm Mold}_{n, d}/{\Bbb Z}, x}$ of ${\rm Mold}_{n, d}$ over ${\Bbb Z}$ at $x$ can be calculated by the Hochschild cohomology $H^{1}({\mathcal A}(x), {\rm M}_n(k(x))/{\mathcal A}(x))$.  We show that $H^{2}({\mathcal A}(x), {\rm M}_n(k(x))/{\mathcal A}(x)) = 0$ is a sufficient condition  for the canonical morphism ${\rm Mold}_{n, d} \to {\Bbb Z}$ being smooth at $x$.  We also calculate $H^{i}(A, {\rm M}_n(k)/A)$ for several $R$-subalgebras $A$ of ${\rm M}_n(R)$ over a commutative ring $R$. 
In particular, we summarize the results on $H^{i}(A, {\rm M}_n(k)/A)$ for all $k$-subalgebras $A$ of ${\rm M}_n(k)$ over an algebraically closed field $k$ in the case $n=2, 3$.     
\end{abstract}

\maketitle

\section{Introduction}
By a rank $d$ {\it mold} ${\mathcal A}$ of degree $n$ on a scheme $X$, we mean a subsheaf of ${\mathcal O}_X$-algebras of ${\rm M}_n({\mathcal O}_X)$ such that ${\mathcal A}$ is a rank $d$ subbundle of 
${\rm M}_n({\mathcal O}_X)$ (Definition \ref{def:molds}).   
Let ${\rm Mold}_{n, d}$ be the moduli of rank $d$ molds of degree $n$ over ${\Bbb Z}$ (Definition and Proposition \ref{defandprop:moduliofmolds}). 
Roughly speaking, ${\rm Mold}_{n, d}$ is the moduli of $d$-dimensional subalgebras of the full matrix ring ${\rm M}_n$ over ${\Bbb Z}$. The moduli ${\rm Mold}_{n, d}$ is a closed subscheme of the Grassmann scheme ${\rm Grass}(d, {\rm M}_n)$ and has rich information on subalgebras of the full matrix ring ${\rm M}_n$. 

Let ${\mathcal A}$ be the universal mold on ${\rm Mold}_{n, d}$. 
For $x \in {\rm Mold}_{n, d}$, denote by ${\mathcal A}(x) = {\mathcal A}\otimes_{{\mathcal O}_{{\rm Mold}_{n, d}}} k(x) \subset {\rm M}_n(k(x))$ the mold corresponding to $x$, where $k(x)$ is the residue field of $x$.  
For investigating ${\rm Mold}_{n, d}$, it is useful to calculate Hochschild cohomology $H^{i}({\mathcal A}(x), {\rm M}_n(k(x))/{\mathcal A}(x))$ for each point $x \in {\rm Mold}_{n, d}$. The dimension of the tangent space $T_{{\rm Mold}_{n, d}/{\Bbb Z}, x}$ of ${\rm Mold}_{n, d}$ over ${\Bbb Z}$ at $x$ can be calculated by the following theorem: 

\begin{theorem}[{\it cf.} Corollary~\ref{cor:dimtangent}]\label{th:main1} 
For each point $x \in {\rm Mold}_{n, d}$, 
\[ 
\dim_{k(x)} T_{{\rm Mold}_{n, d}/{\Bbb Z}, x} = 
\dim_{k(x)} H^{1}({\mathcal A}(x), {\rm M}_n(k(x))/{\mathcal A}(x)) + n^2 - \dim_{k(x)} 
N({\mathcal A}(x)), 
\] 
where $N({\mathcal A}(x)) = \{ b \in {\rm M}_n(k(x)) \mid [b, a]= ba-ab \in {\mathcal A}(x) \mbox{ for any } 
a \in {\mathcal A}(x) \}$. 
\end{theorem} 

\noindent By using $H^{2}({\mathcal A}(x), {\rm M}_n(k(x))/{\mathcal A}(x))$, we obtain a sufficient condition for the canonical morphism ${\rm Mold}_{n, d} \to {\Bbb Z}$ being smooth at $x$: 

\begin{theorem}[Theorem~\ref{th:smooth}]\label{th:main2}   
Let $x \in {\rm Mold}_{n, d}$.  
If $H^2({\mathcal A}(x), {\rm M}_n(k(x))/{\mathcal A}(x)) = 0$, then the canonical morphism 
${\rm Mold}_{n, d} \to {\Bbb Z}$ is smooth at $x$. 
\end{theorem} 

\noindent For a rank $d$ mold ${\mathcal A}$ of degree $n$ on a locally noetherian scheme $S$, we can consider a ${\rm PGL}_{n, S}$-orbit $\{ P^{-1} {\mathcal A}P \mid P \in {\rm PGL}_{n, S} \}$ in ${\rm Mold}_{n, d}\otimes_{{\Bbb Z}} S$, where 
${\rm PGL}_{n, S} = {\rm PGL}_{n}\otimes_{{\Bbb Z}} S$. By using $H^{1}({\mathcal A}(x), {\rm M}_n(k(x))/{\mathcal A}(x))$, we also have: 

\begin{theorem}[Corollary~\ref{cor:openorbit}]\label{th:main3} 
Assume that $H^{1}({\mathcal A}(x), {\rm M}_n(k(x))/{\mathcal A}(x)) = 0$ for each $x \in S$. 
Then the ${\rm PGL}_{n, S}$-orbit $\{ P^{-1}{\mathcal A}P \mid P \in {\rm PGL}_{n, S} \}$ is open in ${\rm Mold}_{n, d}\otimes_{{\Bbb Z}} S$. 
\end{theorem} 

\noindent These theorems are useful to investigate the moduli ${\rm Mold}_{n, d}$. 
We will describe the moduli ${\rm Mold}_{n, d}$ in the case $n=3$ in \cite{Nakamoto-Torii:classification}. 

For $k$-subalgebras $A, B \subseteq {\rm M}_n(k)$ over a field $k$, we say that $A$ and $B$ are {\it equivalent} if there exists $P \in {\rm GL}_n(k)$ such that $P^{-1}AP = B$ (Definition~\ref{def:equivalentsubalgebra}).  
In the case $n=2$, there exist $5$ equivalence classes of $k$-subalgebras of ${\rm M}_2(k)$ over an algebraically closed field $k$ (Proposition~\ref{prop:deg2}).  
In the case $n=3$, there exist $26$ equivalence classes of $k$-subalgebras of ${\rm M}_3(k)$ over an algebraically closed field $k$ (Theorem~\ref{th:deg3}).  For each $k$-subalgebra $A \subseteq {\rm M}_n(k)$ ($n=2, 3$), we calculate the Hochschild cohomology $H^{i}(A, {\rm M}_n(k)/A)$ in Section 5.  We will use the results on $H^{i}(A, {\rm M}_n(k)/A)$ for describing ${\rm Mold}_{n, d}$ in the case $n=3$ in \cite{Nakamoto-Torii:classification}. 

This article is the detailed version of \cite{Nakamoto-Torii:51st} and \cite{Nakamoto-Torii:52nd}.  
In the proof of Theorem~\ref{th:S11} of this paper,  
the Fibonacci numbers appear as the ranks of free modules in the cochain complex for calculating $H^{n}({\rm S}_{11}(R), {\rm M}_3(R)/{\rm S}_{11}(R))$, which seems strange to us, while we have shown another proof of Theorem~\ref{th:S11} using spectral sequence in \cite{Nakamoto-Torii:52nd}. 
We need to point out that our results are closely related with the variety ${\rm Alg}_n$ of $n$-dimensional algebras in the sense of Gabriel in \cite{Gabriel}. Our results can be regarded as a reformulation of Gabriel's theory in the $k$-subalgebra case.  
We will explain  the relation between ${\rm Alg}_{n}$ and ${\rm Mold}_{n, n}$ in another paper. 

\bigskip 

The organization of this paper is as follows: in Section 2, we review Hochschild cohomology. 
For calculating $H^{i}(A, {\rm M}_n/A)$, we introduce several results on Hochschild cohomology.  
In Section 3, we review the moduli of molds. For describing the moduli of molds, we introduce several applications of Hochschild cohomology to the moduli of molds such as Theorems \ref{th:main1}--\ref{th:main3} (Corollary~\ref{cor:dimtangent}, Theorem~\ref{th:smooth}, Corollary~\ref{cor:openorbit}, etc.).  
In Section 4, we explain how to calculate Hochschild cohomology. By using Ciblis's result (Proposition~\ref{prop:cibils}), we can calculate Hochschild cohomology for several cases. We also explain several techniques and perform several calculations. In Section 5, we introduce the classification of $k$-subalgebras of ${\rm M}_n(k)$ over an algebraically closed field $k$ in the case $n=2, 3$. For each $k$-subalgebra $A$ of ${\rm M}_{n}(k)$ ($n=2, 3$), we calculate $H^{i}(A, {\rm M}_n(k)/A)$ for $i \ge 0$. 
In Section 6, we summarize the results on $H^{i}(A, {\rm M}_n(R)/A)$ for $R$-subalgebras $A$ of ${\rm M}_{n}(R)$ over a commutative ring $R$ in the case $n=2, 3$ as Tables \ref{table:deg2} and \ref{table:deg3}.  

\bigskip 

%Let $\mathbb{Z}$, $\mathbb{Q}$,
%$\mathbb{R}$, $\mathbb{C}$
%be the ring of integers,
%the field of rational numbers,
%the field of real numbers,
%the field of complex numbers,
%respectively.

For a commutative ring $R$, we denote by $I_n$ the identity matrix of ${\rm M}_n(R)$. 
We denote by $E_{ij} \in {\rm M}_n(R)$ the matrix with entry $1$ in the $(i, j)$-component and $0$ the other components. 
Set $[I_n] = (I_n \mod R^{\times} \cdot I_n) \in {\rm PGL}_n(R) = {\rm GL}_n(R)/(R^{\times}\cdot I_n)$ for a local ring $R$.  
We also denote by $(R, m, k)$ the triple of a local ring $R$, a maximal ideal $m$ of $R$, and $k=R/m$.  
By a {\it module} $M$ over an associative algebra $A$, we mean a left module $M$ over $A$, unless stated otherwise.

\section{Preliminaries on Hochschild cohomology}
In this section
we give a review of %product structure on
Hochschild cohomology groups
(cf.~\cite{Gerstenhaber} and \cite{Witherspoon}). 
Throughout this section, $R$ denotes a commutative ring, $A$ an associative algebra over $R$, and $M$ an $A$-bimodule over $R$.

\begin{definition}\label{def:hochschildcohomology}\rm 
%Let $A$ be an associative algebra over a commutative ring $R$, and let $M$  be an 
%$A$-bimodule over $R$.  
Assume that $A$ is a projective module over $R$. 
Let $A^e = A\otimes_{R} A^{op}$ be the enveloping algebra of $A$. 
For $A$-bimodules $A$ and $M$ over $R$, we can regard them as $A^e$-modules.  
We define the $i$-th Hochschild cohomology group $H^i (A, M)$ as 
${\rm Ext}^{i}_{A^{e}} (A, M)$.  
%We can calculate $H^i (A, M)$ by taking the cohomology groups of the bar complex 
%$(C^{i}(A, M), d^{i})_{i \in  {\Bbb Z}}$  which is given by 
%\[
%C^{i}(A, M) = \left\{ 
%\begin{array}{cc} 
%{\rm Hom}_{R}(A^{\otimes i}, M) & (i \ge 0) \\ 
%0 & (i<0) \\ 
%\end{array} 
%\right. 
%\] 
%and  $d^{i} : C^{i}(A, M) \to C^{i+1}(A, M)$ $(i \ge 0)$ defined by 
%\begin{multline*} 
%\displaystyle d^{i}(f)(a_1\otimes a_2\otimes \cdots \otimes a_{i+1})  \\ = 
%a_1 f(a_2 \otimes \cdots \otimes a_{i+1}) + 
%\sum_{j=1}^{i} (-1)^j f(a_1\otimes \cdots \otimes a_{j}a_{j+1}\otimes \cdots \otimes a_{i+1}) 
%+ (-1)^{i+1}f(a_1\otimes a_2\otimes \cdots \otimes a_{i})a_{i+1}
%\end{multline*}   for $f \in C^{i}(A, M)$.   
\end{definition}

%Let $R$ be a commutative ring and
%let $A$ be an associative algebra over $R$.
%We assume that $A$ is projective as an $R$-module. 
We denote by
$B_{\ast}(A,A,A)$
the bar resolution of $A$ as $A$-bimodules over $R$.
For $p\ge 0$,
we have
\[ B_p(A,A,A)=A\otimes_R \overbrace{A\otimes_R\cdots\otimes_RA}^p
              \otimes_R A.\]
For an $A$-bimodule $M$ over $R$,
we define a cochain complex $C^{\ast}(A, M)$ to be 
\[ {\rm Hom}_{A^e}(B_{\ast}(A,A,A),M). \]
We can identify
$C^p(A, M)$ with an $R$-module
\[ {\rm Hom}_R(\overbrace{A\otimes_R\cdots\otimes_RA}^{p},M). \]
Under this identification,
the coboundary map $d^p:C^p(A, M)\to C^{p+1}(A, M)$ is given by
\[ \begin{array}{rcl}
   d^p(f)(a_1\otimes\cdots\otimes a_{p+1})&=&
   a_1\cdot f(a_2\otimes \cdots\otimes a_{p+1})\\
  && \displaystyle +\sum_{i=1}^p(-1)^i f(a_1\otimes\cdots\otimes a_ia_{i+1}
                     \otimes\cdots\otimes a_{p+1})\\
  &&+(-1)^{p+1}f(a_1\otimes\cdots\otimes a_p)\cdot a_{p+1}\\
  \end{array}\]  
for $f\in C^p(A, M)$.
The Hochschild cohomology group $H^{\ast}(A, M)$ of $A$ with coefficients 
in $M$ can be calculated by taking the cohomology of the cochain 
complex $C^{\ast}(A, M)$:
\[ H^{\ast}(A, M)=H^{\ast}(C^{\ast}(A, M)).\]     

\begin{remark}\rm
In Definition \ref{def:hochschildcohomology}, the assumption that $A$ is a projective module over $R$ is needed for 
${\rm Ext}^{i}_{A^{e}} (A, M) \cong H^{i}(C^{\ast}(A, M))$ for $i \ge 0$. 
\end{remark} 

Let $N$ be another $A$-bimodule over $R$.
We define a map 
\[ \cup: C^{\ast}(A, M)\times C^{\ast}(A, N)\longrightarrow
   C^{\ast}(A, M\otimes_A N)\]
by
\[ (f\cup g)(a_1\otimes\cdots\otimes a_p\otimes 
    b_1\otimes\cdots\otimes b_q)=
   f(a_1\otimes\cdots \otimes a_p)\otimes
   g(b_1\otimes \cdots\otimes b_q)\]
for $f\in C^p(A, M)$ and $g\in C^q(A, N)$.
The map $\cup$ is $R$-bilinear and satisfies
\[ d^{p+q}(f\cup g)=
   d^p(f)\cup g+ (-1)^pf\cup d^q(g).\]
Hence the map $\cup$ induces a map
\[ H^p(A, M)\otimes_RH^q(A, N)\longrightarrow
   H^{p+q}(A, M\otimes_A N)\]
of $R$-modules.

By the above construction,
we see that the Hochschild cohomology
$H^{\ast}(A, -)$
defines a lax monoidal functor
from the monoidal category of $A$-bimodules over $R$
to the monoidal category of graded $R$-modules.
Hence, 
$H^{\ast}(A, M)$ is a 
graded associative algebra over $R$ 
if $M$ is a monoid object in
the category of $A$-bimodules over $R$.

Suppose that the unit map 
$R\to A$ is a split monomorphism.
We set $\overline{A}=A/RI$,
where $I\in A$ is the image of $1\in R$ under the unit map.
Let $\overline{B}_{\ast}(A,A,A)$ be 
the reduced bar resolution of $A$ as $A$-bimodules over $R$.
We have
\[ \overline{B}_p(A,A,A)\cong
   A\otimes_R \overbrace{\overline{A}\otimes_R\cdots\otimes_R
   \overline{A}}^p\otimes_R A\]
for $p\ge 0$.
For an $A$-bimodule $M$ over $R$, 
we denote  the cochain complex
${\rm Hom}_{A^e}(\overline{B}_{\ast}(A,A,A),M)$
by $\overline{C}^{\ast}(A, M)$.
The cochain complex
$\overline{C}^{\ast}(A,M)$ is a subcomplex of $C^{\ast}(A, M)$.
Recall that the reduced bar resolution $\overline{B}_{\ast}(A,A,A)$
is chain homotopy equivalent to
the bar resolution $B_{\ast}(A,A,A)$,
and hence that the inclusion
$\overline{C}^{\ast}(A, M)\to C^{\ast}(A, M)$ induces
an isomorphism 
\[ H^{\ast}(\overline{C}(A, M))\cong H^{\ast}(A, M). \] 
We observe that the map
$\cup: C^{\ast}(A, M)\times C^{\ast}(A, N)\to
C^{\ast}(A, M\otimes_A N)$ induces an $R$-bilinear map
\[ \cup: \overline{C}^{\ast}(A, M)\times \overline{C}^{\ast}(A, N)\longrightarrow
   \overline{C}^{\ast}(A, M\otimes_A N),\]
where $N$ is another $A$-bimodule over $R$.
Hence the map $\cup: \overline{C}^{\ast}(A, M)\times \overline{C}^{\ast}(A, N)\longrightarrow
   \overline{C}^{\ast}(A, M\otimes_A N)$
induces the same map
$H^p(A, M)\otimes_RH^q(A, N)\longrightarrow
   H^{p+q}(A, M\otimes_A N)$
of $R$-modules as before.

\bigskip 

The following proposition is a basic result. 

\begin{proposition}\label{prop:hochextfield}
Let $A$ be a finite-dimensional associative algebra over a field $k$. 
Let $M$ be an $A$-bimodule over $k$.  
For an extension field $K$ of $k$, we have 
$H^{i}(A\otimes_{k}K, M\otimes_{k}K) \cong 
H^{i}(A, M)\otimes_{k}K$ for $i \ge 0$. 
\end{proposition}

\proof 
The bar complex $C^{i}(A\otimes_{k}K, M\otimes_{k}K) = {\rm Hom}_{K}((A\otimes_{k}K)^{\otimes i}, 
M\otimes_{k}K)$ is isomorphic to $C^{i}(A, M)\otimes_{k}K$. 
Hence $H^{i}(A\otimes_{k}K, M\otimes_{k}K) \cong 
H^{i}(A, M)\otimes_{k}K$ for $i \ge 0$. 
\qed

\begin{corollary}
Let $A, M, k, K$ be as in Proposition \ref{prop:hochextfield}. 
Then $H^{i}(A, M)=0$ if and only if $H^{i}(A\otimes_{k}K, M\otimes_{k}K)=0$. 
\end{corollary} 

\bigskip 

For an $A$-bimodule $M$ over $R$, suppose that there exists a filtration of $A$-bimodules over $R$: 
\[
0 = F^m \subset F^{m-1} \subset \cdots \subset F^{1} \subset F^{0} = M.  
\]
We denote by ${\rm Gr}^{p}(M)$ the $p$-th associated graded module $F^{p}/F^{p+1}$. 
The filtration induces a long exact sequence
\[ \cdots \to H^{p+q}(A, F^{p+1})\to
   H^{p+q}(A, F^p)\to H^{p+q}(A, {\rm Gr}^p(M))
   \to H^{p+q+1}(A, F^{p+1})\to\cdots. \]
We set 
\[ \begin{array}{rcl}
     D^{p,q}&=&H^{p+q}(A, F^p)\\
     E^{p,q}&=&H^{p+q}(A, {\rm Gr}^p(M))\\
   \end{array} \]
We obtain an exact couple
\[
\xymatrix{
D\ar[rr]^{i}&\ar@{}[d]|{}&D\ar[dl]^{j}\\
&E\ar[ul]^{k}&
}
\]
where $D=\oplus_{p,q}D^{p,q}$ and $E=\oplus_{p,q}E^{p,q}$.
By standard construction,
we obtain a spectral sequence: 

\begin{proposition}\label{prop:spectralsequence} 
For a filtration of $A$-bimodules over $R$:   
\[
0 = F^m \subset F^{m-1} \subset \cdots \subset F^{1} \subset F^{0} = M,   
\]
there exists a spectral sequence  
\[ E_1^{p,q}=H^{p+q}(A, {\rm Gr}^p(M))\Longrightarrow
   H^{p+q}(A, M)\]
of $R$-modules with
\[ d_r: E_r^{p,q}\longrightarrow E_r^{p+r,q-r+1} \]
for $r\ge 1$, where ${\rm Gr}^p(M) = F^{p}/F^{p+1}$. 
Here $d_{1} : E_1^{p,q}\longrightarrow E_1^{p+1,q}$ is identified with the connecting homomorphism 
$H^{p+q}(A, {\rm Gr}^{p}(M)) \to H^{p+q+1}(A, {\rm Gr}^{p+1}(M))$ 
of the long exact sequence 
\[ 
\cdots \to H^{\ast}(A, {\rm Gr}^{p+1}(M)) \to H^{\ast}(A, F^{p}/F^{p+2}) \to H^{\ast}(A, {\rm Gr}^{p}(M)) 
\to H^{\ast+1}(A, {\rm Gr}^{p+1}(M)) \to \cdots
\]  
induced by the short exact sequence $0 \to {\rm Gr}^{p+1}(M) \to F^{p}/F^{p+2} \to {\rm Gr}^{p}(M) \to 0$. 
Moreover, the spectral sequence collapses at the $E_m$-page. 
\end{proposition} 

\proof 
See, for example, \cite[\S2.2]{McCleary} for construction
of spectral sequences.
\qed

\section{Applications of Hochschild cohomology groups to the moduli of molds}\label{section:moduli-of-molds}

In this section, we apply Hochschild cohomology to the moduli ${\rm Mold}_{n, d}$ of molds, that is, 
the moduli of subalgebras of the full matrix ring. In Section~\ref{subsection:3.1}, we give a review of the moduli of molds. In Section~\ref{subsection:3.2}, we show that the tangent space of the moduli of molds over ${\Bbb Z}$ at each point $x$ can be calculated by $H^{1}({\mathcal A}(x), {\rm M}_n(k(x))/{\mathcal A}(x))$ (Corollary~\ref{cor:dimtangent}).  In Section~\ref{subsection:3.3}, we construct an obstruction for the canonical morphism ${\rm Mold}_{n, d} \to {\Bbb Z}$ to be smooth at $x$ as a cohomology class of $H^{2}({\mathcal A}(x), {\rm M}_n(k(x))/{\mathcal A}(x))$. Hence $H^{2}({\mathcal A}(x), {\rm M}_n(k(x))/{\mathcal A}(x)) = 0$ is a sufficient condition for ${\rm Mold}_{n, d} \to {\Bbb Z}$ being smooth at $x$ (Theorem~\ref{th:smooth}). In Section~\ref{subsection:3.4}, we discuss the morphism $\phi_{{\mathcal A}} : {\rm PGL}_{n, S} \to {\rm Mold}_{n, d}\otimes_{{\Bbb Z}} S$ defined by $P \mapsto P^{-1}{\mathcal A}P$ for a rank $d$ mold ${\mathcal A} \subseteq {\rm M}_n({\mathcal O}_S)$ on a locally noetherian scheme $S$, where ${\rm PGL}_{n, S} = {\rm PGL}_n\otimes_{{\Bbb Z}} S$.  We show that $\phi_{{\mathcal A}}$ is smooth if and only if $H^{1}({\mathcal A}(x), {\rm M}_n(k(x))/{\mathcal A}(x))=0$ for each $x \in S$ (Theorem~\ref{th:orbitsmooth}).     

\subsection{The moduli of molds}\label{subsection:3.1} 
In this subsection, 
we introduce the 
notion of mold. 
We use \cite{Nakamoto2} as our main reference.

\begin{definition}[{\cite[Definition~1.1]{Nakamoto2}}]\label{def:molds}\rm   
  Let $X$ be a scheme. A subsheaf of ${\mathcal O}_X$-algebras 
${\mathcal A} \subseteq {\rm M}_n({\mathcal O}_X)$ is said to be
a {\it mold} of degree $n$ on $X$ if 
${\mathcal A}$ and ${\rm M}_n({\mathcal O}_X)/{\mathcal A}$
are locally free sheaves on $X$. 
We denote by ${\rm rank} {\mathcal A}$ the rank of ${\mathcal A}$ as a locally free sheaf on X.
For a commutative ring $R$, we say that an $R$-subalgebra 
$A \subseteq {\rm M}_n(R)$ is a {\it mold} of degree $n$ 
over $R$ if $A$ is a mold of degree $n$ on ${\rm Spec}R$.
\end{definition}

\begin{definition}[{\cite[Definition~1.2]{Nakamoto2}}]\label{def:locallyequivalent}\rm
  Let ${\mathcal A}$ and ${\mathcal B}$ be molds of degree 
$n$ on a scheme $X$. 
We say that ${\mathcal A}$ and ${\mathcal B}$ are 
{\it locally equivalent} if for each $x \in X$ 
there exist a neighborhood $U$ of $x$ and 
$P_x \in {\rm GL}_n({\mathcal O}_X(U))$ such that 
$P_x^{-1}({\mathcal A}\mid_U)P_x = {\mathcal B}\mid_U \subseteq 
{\rm M}_n({\mathcal O}_U)$.
\end{definition}

When $X$ is  ${\rm Spec}k$ with a field $k$, we define:  

\begin{definition}\label{def:equivalentsubalgebra}\rm 
Let $k$ be a field. Let $A$ and $B$ be $k$-subalgebras of 
${\rm M}_n(k)$. We say that $A$ and $B$ are {\it equivalent} (or $A \sim B$) if 
there exists $P \in {\rm GL}_n(k)$ such that $P^{-1}AP = B$. 
\end{definition}

We can construct the moduli of molds: 

\begin{defprop}[{\cite[Definition and Proposition~1.1]{Nakamoto2}}]\label{defandprop:moduliofmolds}
  The following contravariant functor is representable by 
a closed subscheme of the Grassmann scheme ${\rm Grass}(d, {\rm M}_n)$:  
\[
\begin{array}{ccccl}
  {\mathcal Mold}_{n, d} & : & ({\bf Sch})^{op} & \to & ({\bf Sets}) \\
  & & X & \mapsto & \{ {\mathcal A} \mid 
\mbox{ a mold of degree } n \mbox{ on } 
X \mbox{ with } {\rm rank}{\mathcal A} = d \}.
\end{array}
\]
We denote by ${\rm Mold}_{n, d}$ the scheme representing the functor 
${\mathcal Mold}_{n, d}$. 
\end{defprop}

\bigskip 

Here we review ${\rm Mold}_{2, d}$ for $d=1, 2, 3, 4$. 
\begin{example}[{\cite[Example~1.1]{Nakamoto2}}]\rm
 In the case $n = 2$, we have
\begin{eqnarray*}
  {\rm Mold}_{2, 1}&  = & {\rm Spec}{{\Bbb Z}}, \label{ex1-1} \\
  {\rm Mold}_{2, 2}&  = & {\Bbb P}^2_{\Bbb Z}, \label{ex1-2} \\
  {\rm Mold}_{2, 3}&  = & {\Bbb P}^1_{\Bbb Z}, \label{ex1-3} \\
  {\rm Mold}_{2, 4}&  = & {\rm Spec}{\Bbb Z}.  \label{ex1-4}
\end{eqnarray*}
\end{example}

\subsection{Tangent spaces of the moduli of molds}\label{subsection:3.2}  

Let $k$ be a field. 
Let $A_0 \in {\rm Mold}_{n, d}(k)$. 
In other words, $A_0$ is a $d$-dimensional 
$k$-subalgebra of ${\rm M}_n(k)$.  
Let ${\mathcal R}$ be the category of Artin local rings with residue field $k$ 
and local homomorphisms. 

\begin{definition}\rm 
We define the covariant functor ${\rm Def}_{A_0} : {\mathcal R} \to ({\bf Sets})$ by 
\[
{\rm Def}_{A_0}(R) = \left\{ 
A \in {\rm Mold}_{n, d}(R) \;  
\begin{array}{|cccc}  
& A\otimes_{R}k & \stackrel{\cong}{\to} & A_0  \\ 
A \mbox{ satisfies that } & \rotatebox{90}{$\hookleftarrow$}  & 
\circlearrowright & 
\rotatebox{90}{$\hookleftarrow$}   \\
& {\rm M}_n(R)\otimes_{R} k & \stackrel{\cong}{\to} & {\rm M}_n(k)   
\end{array}  
\right\} 
\] 
for $R \in {\mathcal R}$. We also define the covariant functor 
$G : {\mathcal R} \to ({\bf Groups})$ by 
\[
G(R) = \left\{ 
P \in {\rm PGL}_n(R) \;  \mid 
P \equiv [I_n] \mod m 
\right\},  
\] 
where $m$ is the maximal ideal of $R \in {\mathcal R}$ and 
$({\bf Groups})$ is the category of groups. 
Then $G(R)$ acts on ${\rm Def}_{A_0}(R)$ from the right by 
\[
\begin{array}{ccc}
{\rm Def}_{A_0}(R) \times G(R) & \to & {\rm Def}_{A_0}(R) \\
(A, P) & \mapsto & P^{-1}AP.  
\end{array}
\]
\end{definition} 

\begin{definition}\rm 
Denote by $({\bf Groupoids})$ the category of groupoids. 
We can regard $({\rm Def}_{A_0}(R), G(R)) \in ({\bf Groupoids})$ for 
$R \in {\mathcal R}$. 
We define the covariant functor 
$F : {\mathcal R} \to ({\bf Groupoids})$ by 
\[
\begin{array}{ccccl}
F & : & {\mathcal R} & \to & ({\bf Groupoids}) \\
& & R & \mapsto & ({\rm Def}_{A_0}(R), G(R)).   
\end{array}
\] 
We also define the covariant functor  
$\pi_0 : ({\bf Groupoids}) \to ({\bf Sets})$ by 
\[
\begin{array}{ccccl}
\pi_0  & : & ({\bf Groupoids}) & \to & ({\bf Sets}) \\
& & G & \mapsto & \{ \mbox{ isomorphism classes of objects of } G \}.   
\end{array}
\] 
Then we have the following composition 
\[
\begin{array}{ccccl}
\pi_0 \circ F & : & {\mathcal R} & \to & ({\bf Sets}) \\
& & R & \mapsto & {\rm Def}_{A_0}(R)/G(R).    
\end{array}
\] 
\end{definition} 

\bigskip

We define the $k$-vector space of derivations ${\rm Der}_{k}(A_0, {\rm M}_n(k)/A_0)$ 
by ${\rm Der}_{k}(A_0, {\rm M}_n(k)/A_0) = 
\{ f \in {\rm Hom}_{k}(A_0, {\rm M}_n(k)/A_0) \mid f(ab) 
= af(b)+f(a)b \mbox{ for } a, b \in A_0 \}$.  

\begin{proposition}\label{prop:defisom}
There exists an isomorphism 
\[ 
{\rm Def}_{A_0}(k[\epsilon]/(\epsilon^2)) \cong 
{\rm Der}_{k}(A_0, {\rm M}_n(k)/A_0). 
\]
\end{proposition} 

\proof
For $\theta \in {\rm Der}_{k}(A_0, {\rm M}_n(k)/A_0)$, take a $k$-linear map 
$\theta' : A_0 \to {\rm M}_n(k)$ as a lift of $\theta$.  
We define $A(\theta) = (k[\epsilon]/(\epsilon^2)) \{ \ a+\theta'(a)\epsilon 
\mid a \in A_0 \} \subset {\rm M}_n(k[\epsilon]/(\epsilon^2))$. 
It is easy to check that the definition of $A(\theta)$ does not 
depend on the choice of $\theta'$. We define a map 
${\rm Der}_{k}(A_0, {\rm M}_n(k)/A_0) \to {\rm Def}_{A_0}(k[\epsilon]/(\epsilon^2))$ by 
$\theta \mapsto A(\theta)$. We can easily prove that this map is bijective. 
\qed

\begin{definition}[{\cite[16.5.13]{EGA4}, \cite[Definition~0B2C]{Stacks}}]\rm
Let $f : X \to S$ be a morphism of schemes. 
Let $x \in X$ and $s = f(x) \in S$. We denote by $k(x)$ and $k(s)$ the residue fields of 
$x$ and $s$, respectively. The field extension $k(s) \subseteq k(x)$ induces 
$k(s)$-algebra homomorphisms $k(s) \stackrel{\varphi_1}{\to} k(x)[\epsilon]/(\epsilon^2) \stackrel{\varphi_2}{\to} 
k(x)$ such that $\varphi_1(a) = a+0\cdot \epsilon$ for $a \in k(s)$, $\varphi_2(b+c\epsilon)=b$ for $b, c \in k(x)$,  and $\varphi_2\circ \varphi_1$ is the inclusion $k(s) \hookrightarrow k(x)$. By $\varphi_1$ and $\varphi_2$, we obtain morphisms ${\rm Spec} \: k(x) \to {\rm Spec} \: k(x)[\epsilon]/(\epsilon^2) \to {\rm Spec} \: k(s)$. 
By a {\it tangent vector} of $X/S$ at $x$, we mean an $S$-morphism 
$\psi : {\rm Spec} \: k(x)[\epsilon]/(\epsilon^2) \to X$ 
such that the following diagram is commutative: 
\[
\begin{array}{rcc}
{\rm Spec} \: k(x) & & \\ \vspace*{1ex} 
\downarrow {\varphi_2^{\ast}}  & \searrow & \\ 
{\rm Spec} \: k(x)[\epsilon]/(\epsilon^2) & \stackrel{\psi}{\to} & X \\ \vspace*{1ex} 
\downarrow {\varphi_1^{\ast}} &  &  \downarrow \\
{\rm Spec} \: k(s) & \to & S.  
\end{array} 
\]
We call the set of tangent vectors of $X/S$ at $x$ the {\it tangent space} $T_{X/S, x}$ of $X$ over $S$ at $x$, which has a canonical $k(x)$-vector space structure. 
\end{definition} 

\begin{remark}[{\cite[16.5.13.1]{EGA4}, \cite[(0BEA) and Lemma~0B2D]{Stacks}}]\rm 
Let $X_{s}$ be the scheme-theoretic fiber $f : X \to S$ over $s=f(x)$. Then there exists a canonical isomorphism 
$T_{X/S, x} \cong T_{X_{s}/{\rm Spec}\: k(s), x}$ as $k(x)$-vector spaces. 
Let ${\Omega}_{X/S}$ be the sheaf of relative differentials of $X$ over $S$.  
We also have a canonical isomorphism 
\[
T_{X/S, x} \cong {\rm Hom}_{{\mathcal O}_{X, x}}(\Omega_{X/S, x}, k(x)) 
\]
as $k(x)$-vector spaces, where $\Omega_{X/S, x}$ is the stalk of $\Omega_{X/S}$ at $x$. 
\end{remark} 

Let ${\mathcal A}$ be the universal mold on ${\rm Mold}_{n, d}$.  
For a point $x$ of ${\rm Mold}_{n, d}$, we denote by ${\mathcal A}(x) = {\mathcal A}\otimes_{{\mathcal O}_{{\rm Mold}_{n, d}}} k(x) \subseteq {\rm M}_n(k(x))$ the mold corresponding to $x$. 

\begin{corollary}\label{cor:defisom}  
Let $x$ be a point of ${\rm Mold}_{n, d}$.  
The tangent space $T_{{\rm Mold}_{n, d}/{\Bbb Z}, x}$ of ${\rm Mold}_{n, d}$ over ${\Bbb Z}$ at $x$ is isomorphic to 
${\rm Der}_{k(x)}({\mathcal A}(x), {\rm M}_n(k(x))/{\mathcal A}(x))$ as $k(x)$-vector spaces. 
\end{corollary} 

\proof  
We see that $T_{{\rm Mold}_{n, d}/{\Bbb Z}, x} 
\cong {\rm Def}_{{\mathcal A}(x)}(k(x)[\epsilon]/(\epsilon^2)) \cong {\rm Der}_{k(x)}({\mathcal A}(x), {\rm M}_n(k(x))/{\mathcal A}(x))$ by the definition of tangent space and Proposition \ref{prop:defisom}. We also see that $T_{{\rm Mold}_{n, d}/{\Bbb Z}, x}$ is canonically isomorphic to ${\rm Der}_{k(x)}({\mathcal A}(x), {\rm M}_n(k(x))/{\mathcal A}(x))$ as $k(x)$-vector spaces. 
\qed 

\bigskip 

Let us define $d : {\rm M}_n(k) \to {\rm Der}_{k}(A_0, {\rm M}_n(k)/A_0)$ by 
$d(X)(a) = ([X, a] \mod A_0)= (Xa-aX \mod A_0)$ for $X \in {\rm M}_n(k)$ and $a \in A_0$. 
It is easy to check that $d(X) \in {\rm Der}_{k}(A_0, {\rm M}_n(k)/A_0)$.

\begin{proposition}\label{prop:h1} 
There exists an isomorphism 
\[ 
H^{1}({A_0}, {\rm M}_n(k)/A_0) \cong 
{\rm Der}_{k}(A_0, {\rm M}_n(k)/A_0)/{\rm Im}\: d. 
\]
\end{proposition} 

\proof
Let us consider the bar complex 
\[ 
0 \to C^{0}(A_0, {\rm M}_n(k)/A_0) 
\stackrel{d^0}{\to} C^{1}(A_0, {\rm M}_n(k)/A_0) \stackrel{d^1}{\to} 
C^{2}(A_0, {\rm M}_n(k)/A_0)\to\cdots. 
\] 
Note that ${\rm Ker} \:d^{1} = {\rm Der}_{k}(A_0, {\rm M}_n(k)/A_0) \supseteq {\rm Im} \:d^0 = 
{\rm Im}\: d$. Hence we have 
$H^{1}(A_0, {\rm M}_n(k)/A_0) \cong {\rm Der}_{k}(A_0, {\rm M}_n(k)/A_0)/{\rm Im}\: d$. 
\qed 

\bigskip 

Let $N(A_0) = \{ X \in {\rm M}_n(k) \mid [X, a] \in A_0 \mbox{ for  any } a \in A_0 \}$.  
The $k$-linear map $d : {\rm M}_n(k) \to {\rm Der}_{k}(A_0, {\rm M}_n(k)/A_0)$ induces 
a $k$-linear map $\overline{d} : {\rm M}_n(k)/A_{0} \to {\rm Der}_{k}(A_0, {\rm M}_n(k)/A_0)$. 
Then we have  

\begin{corollary}\label{cor:dimtangent} 
There exists an exact sequence 
\[
0 \to N(A_0)/A_0 \to {\rm M}_n(k)/A_0 \stackrel{\overline{d}}{\to} {\rm Der}_{k}(A_0, {\rm M}_n(k)/A_0)  
\to H^{1}({A_0}, {\rm M}_n(k)/A_0)  \to 0. 
\]
In particular, 
$\dim_{k(x)} T_{{\rm Mold}_{n, d}/{\Bbb Z}, x} = 
\dim_{k(x)} H^{1}({\mathcal A}(x), {\rm M}_n(k(x))/{\mathcal A}(x)) + n^2 - \dim_{k(x)} 
N({\mathcal A}(x))$ for any point $x \in {\rm Mold}_{n, d}$. 
\end{corollary}

\proof 
By Proposition \ref{prop:h1}, ${\rm M}_n(k)/A_0 \stackrel{\overline{d}}{\to} 
{\rm Der}_{k}(A_0, {\rm M}_n(k)/A_0) \to H^{1}({A_0}, {\rm M}_n(k)/A_0)  \to 0$ is exact. 
The kernel of $\overline{d}$ is equal to $N(A_0)/A_0$. Hence we have the exact sequence above. 
The last statement follows from the fact that $T_{{\rm Mold}_{n, d}/{\Bbb Z}, x}   
\cong {\rm Der}_{k(x)}({\mathcal A}(x), 
{\rm M}_n(k(x))/{\mathcal A}(x))$ by Corollary \ref{cor:defisom}. 
\qed 

\bigskip 

By the definition, $G(k[\epsilon]/(\epsilon^2)) = \{ [I_n+X\epsilon] \in 
{\rm PGL}_n(k[\epsilon]/(\epsilon^2)) \mid X \in {\rm M}_n(k) \}$. Note that 
$[I_n+X\epsilon] = [I_n+Y\epsilon]$ if and only if there exists $c \in k$ such that 
$X = Y+ cI_n$. Hence we have a group isomorphism 
$G(k[\epsilon]/(\epsilon^2)) \cong {\rm M}_n(k)/kI_n$ defined by $[I_n+X\epsilon] \mapsto (X \mod kI_n)$.  
Recall $A(\theta) \in {\rm Def}_{A_0}(k[\epsilon]/(\epsilon^2))$ defined in 
Proposition \ref{prop:defisom} for 
$\theta \in {\rm Der}_{k}(A_0, {\rm M}_n(k)/A_0)$. 
For $P = [I_n + X\epsilon]  \in G(k[\epsilon]/(\epsilon^2))$, 
\begin{eqnarray*}
P^{-1}A(\theta)P & = & [I_n-X\epsilon] A(\theta) [I_n+X\epsilon] \\ 
 & = &  (k[\epsilon]/(\epsilon^2))\{ a + (\theta'(a)+ aX-Xa)\epsilon \mid a \in A_0 \}  \\ 
 & = & A(\theta - d(X)),   
\end{eqnarray*} 
where $\theta' : A_0 \to {\rm M}_n(k)$ is a lift of $\theta$. Hence we have
 
\begin{proposition}
We have isomorphisms 
\[
H^{1}({A_0}, {\rm M}_n(k)/A_0) \cong {\rm Def}_{A_0}(k[\epsilon]/(\epsilon^2))/G(k[\epsilon]/(\epsilon^2)) 
\cong \pi_0 \circ F(k[\epsilon]/(\epsilon^2)). 
\]
\end{proposition}  

\proof
By the discussion above and Propositions \ref{prop:defisom} and  
\ref{prop:h1}, we can prove the statement. 
\qed 

\bigskip 

Let us consider $H^{0}({A_0}, {\rm M}_n(k)/A_0)$.  

\begin{definition}\rm 
We define the trivial deformation $A_{\epsilon} = A_0\otimes_{k} (k[\epsilon]/(\epsilon^2)) 
\in {\rm Def}_{A_0}(k[\epsilon]/(\epsilon^2))$ of $A_0$ to $k[\epsilon]/(\epsilon^2)$. 
Note that $A_{\epsilon} = A(0)$, where $0 \in {\rm Der}_{k}(A_0, {\rm M}_n(k)/A_0)$.  
We also define $G_{\epsilon} = \{ P \in G(k[\epsilon]/(\epsilon^2)) \mid P^{-1}A_{\epsilon}P = A_{\epsilon} \}$. Then $G_{\epsilon}$ is equal to the stabilizer group 
${\rm Aut}_{F(k[\epsilon]/(\epsilon^2))}(A_{\epsilon})$ of $A_{\epsilon}$ in the 
groupoid $F(k[\epsilon]/(\epsilon^2)) = ({\rm Def}_{A_0}(k[\epsilon]/(\epsilon^2)), 
G(k[\epsilon]/(\epsilon^2)))$. 
\end{definition} 

\begin{proposition}
There exists an exact sequence 
\[
0 \to I \to G_{\epsilon} \to H^{0}({A_0}, {\rm M}_n(k)/A_0) \to 0, 
\]
where $I = \{ P \in G(k[\epsilon]/(\epsilon^2)) \mid P=[I_n+X\epsilon], X \in A_0 \}$. 
\end{proposition}

\proof
Recall the isomorphism 
$G(k[\epsilon]/(\epsilon^2)) \cong {\rm M}_n(k)/kI_n$ given by $[I_n+X\epsilon] \mapsto 
(X\!\mod kI_n)$. For $P=[I_n+X\epsilon] \in G(k[\epsilon]/(\epsilon^2))$, 
$P^{-1}A_{\epsilon}P=P^{-1}A(0)P=A(-d(X))$. Hence we have $G_{\epsilon} = 
\{ P \in G(k[\epsilon]/(\epsilon^2)) \mid P=[I_n+X\epsilon], d(X)=0 \}
\cong \{  X \in {\rm M}_n(k) \mid [A_0, X]\subseteq A_0 \}/kI_n$. 
Let $\overline{d} : {\rm M}_n(k)/A_{0} \to {\rm Der}_{k}(A_0, {\rm M}_n(k)/A_0)$ be the 
$k$-linear map induced by $d : {\rm M}_n(k) \to {\rm Der}_{k}(A_0, {\rm M}_n(k)/A_0)$. 
By the bar complex, $H^{0}(A_0, {\rm M}_n(k)/A_0) = {\rm Ker} \: \overline{d} = 
\{ [X] \in {\rm M}_n(k)/A_{0} \mid [A_0, X] \subseteq A_0 \}$. 
The canonical projection ${\rm M}_n(k)/kI_n \to {\rm M}_n(k)/A_{0}$ induces 
a surjective homomorphism $p : G_{\epsilon} \to H^{0}(A_0, {\rm M}_n(k)/A_0)$. 
The kernel of $p$ is $I = 
\{ P \in G(k[\epsilon]/(\epsilon^2)) \mid P=[I_n+X\epsilon], X \in A_0 \} 
\cong A_{0}/kI_n$. This complete the proof. 
\qed

\subsection{Smoothness of ${\rm Mold}_{n, d}$}\label{subsection:3.3}  

Let $(\widetilde{R}, \widetilde{m}, k)$ be an Artin local ring. 
Let $I$ be an ideal of $\widetilde{R}$ such that $\widetilde{m} I = 0$. 
Set $R= \widetilde{R}/I$ and $m= \widetilde{m}/I$. 
Then $(R, m, k)$ is also an Artin local ring. 
%Denote by $\pi : {\rm M}_n(\widetilde{R}) \to {\rm M}_n(R)$ the canonical projection. 
%Let $s : R \to \widetilde{R}$ be a set-theoretical section of $\pi$. 
Let $A \in {\rm Mold}_{n, d}(R)$. In other words, $A \subset {\rm M}_n(R)$ is  
a rank $d$ mold.  Since $R$ is a local ring, $A$ and ${\rm M}_n(R)/A$ are free modules over $R$. 
%Take a basis $a_1, a_2, \ldots, a_{n^2}$ of ${\rm M}_n(R)$ over $R$ 
%such that 
Take a basis $a_1, a_2, \ldots, a_{d}$ of $A$ over $R$. 
For $1 \le i \le d$, choose a lift $S(a_i) \in {\rm M}_n(\widetilde{R})$ of $a_i$.    
Since $a_i a_j \in A = Ra_1 \oplus Ra_2 \oplus \cdots \oplus Ra_d$, 
we can choose a lift $S(a_ia_j) \in \widetilde{R}S(a_1) \oplus \widetilde{R}S(a_2) \oplus \cdots \oplus \widetilde{R}S(a_d)$ of $a_i a_j$ for $1 \le i, j \le d$. Note that $S(a_ia_j) - S(a_i)S(a_j) \in {\rm M}_n(I)$. 
Let us define an $R$-linear map 
$c' : A\otimes_{R} A \to {\rm M}_n(I) \cong {\rm M}_n(k)\otimes_{k} I$ 
by $c'(\sum_{1 \le i, j \le d} r_{ij} a_i\otimes a_j) = \sum_{1 \le i, j \le d} \widetilde{r}_{ij}(S(a_ia_j)-S(a_i)S(a_j))$, 
where $\widetilde{r}_{ij} \in \widetilde{R}$ is a lift of $r_{ij} \in R$.  
The $R$-module structure of ${\rm M}_n(k)\otimes_{k} I$ is given by 
$a \cdot (X\otimes x) = (p(a)X) \otimes x$ for 
$a \in R, X \in {\rm M}_n(k)$,  and $x \in I$, where 
$p : R \to R/m = \widetilde{R}/\widetilde{m} = k$ is  
the canonical projection. 
By using $I^2=0$, we easily see that the definition of $c'$ does not depend on the choice of  lifts $\widetilde{r}_{ij}$ of $r_{ij}$. 

%Then we define $S : {\rm M}_n(R) \to {\rm M}_n(\widetilde{R})$ by 
%$S(\sum_{i=1}^{n^2} r_i a_i) = \sum_{i=1}^{n^2} s(r_i)S(a_i)$ for 
%$\sum_{i=1}^{n^2} r_i a_i \in {\rm M}_n(R)$. 
%Note that $S : {\rm M}_n(R) \to {\rm M}_n(\widetilde{R})$ does not necessarily 
%coincide with the map given by applying $s : R \to \widetilde{R}$ 
%to each entries of matrices in ${\rm M}_n(R)$. 
%We simply call such $s : R \to \widetilde{R}$ and $S : {\rm M}_n(R) \to {\rm M}_n(\widetilde{R})$ 
%{\it lifts} of the projections $\pi : \widetilde{R} \to R$ and 
% $\pi : {\rm M}_n(\widetilde{R}) \to {\rm M}_n(R)$.   
%
%Let us define an $R$-linear map 
%$c' : A\otimes_{R} A \to {\rm M}_n(I) \cong {\rm M}_n(k)\otimes_{k} I$ 
%by $c'(\sum_{1 \le i, j \le d} r_{ij} a_i\otimes a_j) = \sum_{1 \le i, j \le d} 
%s({r}_{ij}) (S(a_ia_j)-S(a_i)S(a_j))$ for $r_{ij} \in R$.  Note that the 
%$R$-module structure of ${\rm M}_n(k)\otimes_{k} I$ is given by 
%$a \cdot (X\otimes x) = (p(a)X) \otimes x$ for 
%$a \in R, X \in {\rm M}_n(k)$,  and $x \in I$, where 
%$p : R \to R/m \cong  \widetilde{R}/\widetilde{m} \cong k$ is  
%the canonical projection. 
%Since $s(a_ia_j)-s(a_i)s(a_j) \in {\rm M}_n(I)$, $c'$ is well-defined.  

Set $A_0 = A\otimes_{R} k \subseteq {\rm M}_n(k)$.  
Since $A = \oplus_{i=1}^{d} Ra_i$, we can write  
$A_0 = \oplus_{i=1}^{d} k\overline{a}_i$, where $\overline{a}_i = (a_i \mod m)$.  
We denote by $c''$ the composition $A\otimes_{R} A \stackrel{c'}{\to} {\rm M}_n(k)\otimes_{k} I 
\to ({\rm M}_n(k)/A_{0}) \otimes_{k} I$. 
It is easy to see that $c'' : A\otimes_{R} A \to ({\rm M}_n(k)/A_{0}) \otimes_{k} I$ goes through 
$A_{0}\otimes_{k} A_0$.  
Hence we have a $k$-linear map $c : A_{0}\otimes_{k} A_0 \to ({\rm M}_n(k)/A_{0}) \otimes_{k} I$. 
By $a\cdot (\overline{X}\otimes x) \cdot b = (a\overline{X}b)\otimes x$ for 
$\overline{X} \in {\rm M}_n(k)/A_{0}$, $x \in I$ and $a, b \in A_0$, we can 
regard $({\rm M}_n(k)/A_{0}) \otimes_{k} I$ 
as an $A_0$-bimodule.  
For the $A_0$-bimodule $({\rm M}_n(k)/A_{0}) \otimes_{k} I$, let us consider 
the bar complex $(C^i, d^i)_{i \in {\mathbb Z}}$, where 
$C^i = C^i(A_0, ({\rm M}_n(k)/A_0) \otimes_{k} I) = 
{\rm Hom}_{k}(A_{0}^{\otimes i}, ({\rm M}_n(k)/A_{0}) \otimes_{k} I)$ and 
$d^i : C^{i} \to C^{i+1}$ is a differential.   

\begin{lemma}\label{lemma:cocycle}
The $k$-linear map $c : A_{0}\otimes_{k} A_0 \to ({\rm M}_n(k)/A_{0}) \otimes_{k} I$ is a $2$-cocyle in $C^2$.  
\end{lemma} 

\proof 
Let us show that $d^2(c)=0$, where $d^2 : C^2 \to C^3$. 
% = {\rm Hom}_{k}(A_0\otimes_{k} A_0\otimes_{k} A_0, ({\rm M}_n(k)/A_{0}) \otimes_{k} I)$. 
%Since $A = \oplus_{i=1}^{d} Ra_i$, we can write  
%$A_0 = \oplus_{i=1}^{d} k\overline{a}_i$, where $\overline{a}_i = (a_i \mod m)$.  
It suffices to show that $d^2(c)(\overline{a}_i\otimes \overline{a}_j \otimes 
\overline{a}_l)=0$ for $1 \le i, j, l \le d$. 
By the definition, 
\begin{eqnarray*} 
d^2(c)(\overline{a}_i \otimes \overline{a}_j \otimes \overline{a}_l) 
=  \overline{a}_ic(\overline{a}_j\otimes \overline{a}_l) - 
 c(\overline{a}_i\overline{a}_j\otimes \overline{a}_l) 
+ c(\overline{a}_i\otimes \overline{a}_j\overline{a}_l) - 
c(\overline{a}_i\otimes \overline{a}_j)\overline{a}_l.  
\end{eqnarray*}  
For $1 \le i, j \le d$, there exist $c_{ij}^{s} \in R$ such that 
$a_i a_j 
= \sum_{s=1}^{d} c_{ij}^{s} a_s \in A$. 
Putting $\overline{c}_{ij}^{s} = (c_{ij}^{s} \mod m)$, we have 
$\overline{a}_i \overline{a}_j 
= \sum_{s=1}^{d} \overline{c}_{ij}^{s} \overline{a}_s \in {\rm M}_n(k)$. 
For verifying $d^2(c)(\overline{a}_i \otimes \overline{a}_j \otimes \overline{a}_k) = 0$ in 
$({\rm M}_n(k)/A_0)\otimes_{k} I = ({\rm M}_n(k)\otimes_{k} I) /(A_{0}\otimes_{k} I)$, 
we calculate $S(a_i)c'(a_j\otimes a_l), c'(a_ia_j\otimes a_l), c'(a_i\otimes a_ja_l)$, and 
$c'(a_i\otimes a_j)S(a_l)$   
in ${\rm M}_n(k)\otimes_{k} I = {\rm M}_n(I) \subset {\rm M}_n(\widetilde{R})$. 
%Choose a lift $\widetilde{c}_{ij}^{s} \in \widetilde{R}$ of $c_{ij}^{s} \in R$ for $1 \le i, j, s \le d$.   
We can write $S(a_ia_j) = \sum_{s=1}^{d}  \widetilde{c}_{ij}^{s} S(a_s) \in \oplus_{s=1}^{d} \widetilde{R} S(a_s) \subseteq {\rm M}_n(\widetilde{R})$, where $\widetilde{c}_{ij}^{s} \in \widetilde{R}$ is a lift of $c_{ij}^{s}$.      
Since 
\begin{eqnarray*} 
S(a_i)c'(a_j\otimes a_l) & = & S(a_i)(S(a_ja_l)-S(a_j)S(a_l)) \\ 
 & = & S(a_i)S(a_ja_l)-S(a_i)S(a_j)S(a_l), 
\end{eqnarray*} 
\begin{eqnarray*} 
c'(a_ia_j\otimes a_l) & = 
& c'\left(\sum_{s=1}^{d} c_{ij}^{s} a_s \otimes  a_l \right) \\ 
& = & \sum_{s=1}^{d} \widetilde{c}_{ij}^{s}( S(a_sa_l)-S(a_s)S(a_l)) \\ 
& = & \sum_{s, t=1}^{d} \widetilde{c}_{ij}^{s} \widetilde{c}_{sl}^{t} S(a_t)  - \sum_{s=1}^{d} \widetilde{c}_{ij}^{s} S(a_s)S(a_l) \\
& = & \sum_{s, t=1}^{d} \widetilde{c}_{ij}^{s} \widetilde{c}_{sl}^{t} S(a_t)  - S(a_ia_j)S(a_l),
\end{eqnarray*}
\begin{eqnarray*} 
c'(a_i\otimes a_ja_l) & = &   c'\left(\sum_{s=1}^{d} c_{jl}^{s} a_i\otimes a_s \right) \\  
& = & \sum_{s=1}^{d} \widetilde{c}_{jl}^{s}(S(a_ia_s)-S(a_i)S(a_s)) \\ 
& = & \sum_{s, t=1}^{d} \widetilde{c}_{jl}^{s} \widetilde{c}_{is}^{t}S(a_t) - \sum_{s=1}^{d} \widetilde{c}_{jl}^{s}S(a_i)S(a_s) \\  
& = & \sum_{s, t=1}^{d} \widetilde{c}_{jl}^{s} \widetilde{c}_{is}^{t}S(a_t) - S(a_i)S(a_ja_l), \\  
\end{eqnarray*} 
and 
\begin{eqnarray*} 
c'(a_i\otimes a_j)S(a_l) & = & (S(a_ia_j)-S(a_i)S(a_j)) S(a_l) \\ 
& = & S(a_ia_j)S(a_l) - S(a_i)S(a_j)S(a_l),  
\end{eqnarray*} 
we have 
\begin{eqnarray}\label{eq:differential}  
& & S(a_i)c'(a_j\otimes a_l)-c'(a_ia_j\otimes a_l)+c'(a_i\otimes a_ja_l)
- c'(a_i\otimes a_j)S(a_l)   \\ 
& = &  \sum_{s, t=1}^{d} ( \widetilde{c}_{is}^{t} \widetilde{c}_{jl}^{s}   
-  \widetilde{c}_{ij}^{s} \widetilde{c}_{sl}^{t} ) S(a_t).  \nonumber 
\end{eqnarray} 
The associativity $a_i(a_ja_l)=(a_ia_j)a_l$ implies that 
$\sum_{s=1}^{d} (\widetilde{c}_{is}^{t} \widetilde{c}_{jl}^{s}   
-  \widetilde{c}_{ij}^{s} \widetilde{c}_{sl}^{t}) \in I$ for each $1 \le t \le d$. 
The right hand side of (\ref{eq:differential}) is contained in 
$\sum_{i=1}^{d} I S(a_i)  = A_0 \otimes_{k} I \subset {\rm M}_n(k)\otimes_{k} I = {\rm M}_n(I)$.  
Thus, we have $d^2(c)(\overline{a}_i\otimes \overline{a}_j \otimes 
\overline{a}_l)=0$. 
\qed

\bigskip 

For lifts $S(a_i) \in {\rm M}_n(\widetilde{R})$ of $a_i$ ($1 \le i \le d$) and lifts 
$S(a_ia_j) \in \widetilde{R}S(a_1)\oplus \cdots \oplus \widetilde{R}S(a_d)$ of $a_ia_j$ ($1 \le i, j \le d$),  
we can define a $2$-cocycle $c_{S} \in C^2$ by Lemma \ref{lemma:cocycle}.  
We denote by $[c_{S}]$ the cohomology class of $c_{S}$ %defined by $(s, S)$  
in $H^2(A_0, ({\rm M}_n(k)/A_0)\otimes_{k} I)$. 

\begin{lemma}\label{lemma:independentlift}
The cohomology class $[c_{S}]$ in $H^2(A_0, ({\rm M}_n(k)/A_0)\otimes_{k} I)$ is 
independent from the choice of the lifts $S(a_i) \in {\rm M}_n(\widetilde{R})$  of $a_i$ $(1 \le i \le d)$ and the lifts $S(a_ia_j) \in  \widetilde{R}S(a_1)\oplus \cdots \oplus \widetilde{R}S(a_d)$ of $a_ia_j$ $(1 \le i, j \le d)$.  
\end{lemma} 

\proof 
Let $T(a_i) \in {\rm M}_n(\widetilde{R})$ and $T(a_ia_j) \in  \widetilde{R}T(a_1)\oplus \cdots \oplus \widetilde{R}T(a_d)$ be other lifts of $a_i$ ($1 \le i \le d$) and $a_ia_j$ ($1 \le i, j \le d$), respectively.  
We denote by $c_{S}, c_{T} : A_{0}\otimes_{k} A_{0} \to ({\rm M}_n(k)/A_0)\otimes_{k} I$ the 
$2$-cocycles defined by the lifts $\{ S(a_i) \} \cup \{ S(a_ia_j) \}$ and $\{ T(a_i) \} \cup \{ T(a_ia_j) \}$, respectively. 
We define the $k$-linear map $\theta : A_0 \to ({\rm M}_n(k)/A_0)\otimes_{k} I = ({\rm M}_n(k)\otimes_{k} I)/(A_{0}\otimes_{k} I)$ by 
$\overline{a}_i \mapsto (T(a_i)-S(a_i) \mod A_0\otimes_{k} I)$ for $1 \le i \le d$. 
Note that $T(a_i)-S(a_i) \in {\rm M}_n(I)={\rm M}_n(k)\otimes_{k} I $.  
Let us calculate $d^{1}(\theta)(\overline{a}_i\otimes \overline{a}_j)  =  
\overline{a}_i\theta(\overline{a}_j) - \theta(\overline{a}_i\overline{a}_j) + \theta(\overline{a}_i)
\overline{a}_j $. Put $a_ia_j = \sum_{s=1} c_{ij}^{s} a_s$ for $c_{ij}^{s} \in R$ and $\overline{c}_{ij}^{s} = (c_{ij}^s \mod m) \in k$. 
We can write $S(a_ia_j) = \sum_{s=1}^{d} \widetilde{c}_{ij}^{s} S(a_s)$ and  
$T(a_ia_j) = \sum_{s=1}^{d} \widetilde{d}_{ij}^{s} T(a_s)$ for $1 \le i, j \le d$, where 
$\widetilde{c}_{ij}^{s}, \widetilde{d}_{ij}^{s} \in {\rm M}_n(\widetilde{R})$ are lifts of $c_{ij}^{s}$. 
Using $\widetilde{c}_{ij}^{s} - \widetilde{d}_{ij}^{s} \in I$ and $(\widetilde{c}_{ij}^{s}-\widetilde{d}_{ij}^{s})T(a_s) \in A_0\otimes_{k}I$, we have  
\begin{eqnarray*} 
\theta(\overline{a}_i\overline{a}_j) &=& \theta \left(\sum_{s=1}^{d} \overline{c}_{ij}^{s} \overline{a}_s \right) \\ 
& = &  \sum_{s=1}^{d} \widetilde{c}_{ij}^{s}(T(a_s) - S(a_s)) \mod A_0\otimes_{k}I \\ 
& = & \sum_{s=1}^{d} \widetilde{d}_{ij}^{s}T(a_s) - \sum_{s=1}^d \widetilde{c}_{ij}^{s} S(a_s) + \sum_{s=1}^d (\widetilde{c}_{ij}^{s}-\widetilde{d}_{ij}^{s})T(a_s) \mod A_0\otimes_{k}I \\ 
&=& \sum_{s=1}^{d} \widetilde{d}_{ij}^{s}T(a_s) - \sum_{s=1}^d \widetilde{c}_{ij}^{s} S(a_s) \mod A_0\otimes_{k}I \\ 
&=&T(a_ia_j)-S(a_ia_j) \mod A_0\otimes_{k}I.  
\end{eqnarray*} 
%\begin{multline*} 
%\theta(\overline{a}_i\overline{a}_j) = \theta(\sum_{s=1}^{d} \overline{c}_{ij}^{s} \overline{a}_s) = \sum_{s=1}^{d} %\widetilde{c}_{ij}^{s}(T(a_s) - S(a_s)) \mod A_0\otimes_{k}I \\ 
%= \sum_{s=1}^{c} \widetilde{d}_{ij}^{s}T(a_s) - \sum_{s=1}^d \widetilde{c}_{ij}^{s} S(a_s) + \sum_{s=1}^d (\widetilde{c}_{ij}%^{s}-\widetilde{d}_{ij}^{s})T(a_s) \mod A_0\otimes_{k}I \\ 
%= \sum_{s=1}^{d} \widetilde{d}_{ij}^{s}T(a_s) - \sum_{s=1}^d \widetilde{c}_{ij}^{s} S(a_s) \mod A_0\otimes_{k}I \\ 
%=T(a_ia_j)-S(a_ia_j) \mod A_0\otimes_{k}I.  
%\end{multline*} 
Since $\overline{a}_i X = S(a_i)X$ and $X \overline{a}_j = XT(a_j)$ 
for each $X \in {\rm M}_n(I)$, we have 
\begin{eqnarray*}
& & d^{1}(\theta)(\overline{a}_i\otimes \overline{a}_j) \\ 
& = & \overline{a}_i\theta(\overline{a}_j) - \theta(\overline{a}_i\overline{a}_j) + \theta(\overline{a}_i)
\overline{a}_j \\ 
& = & S(a_i)(T(a_j)-S(a_j)) - (T(a_ia_j)-S(a_ia_j)) +(T(a_i)-S(a_i))T(a_j) \mod A_0\otimes_{k} I \\
& = & c_{S}(\overline{a}_i\otimes \overline{a}_j)-c_{T}(\overline{a}_i\otimes \overline{a}_j). 
\end{eqnarray*} 
Hence we have $c_{S} - c_{T} = d^{1}(\theta)$, which implies that $[c_{S}] = [c_{T}]$. 
\qed 

%\begin{remark}\rm
%%For lifts $s : R \to \widetilde{R}$ and $S :  {\rm M}_n(R) \to {\rm M}_n(\widetilde{R})$, 
%%the cocycle $c_s : A_0\otimes_{k} A_0 \to ({\rm M}_n(k)/A_{0})\otimes I$ is 
%%determined by $\{  c_{s}(\overline{a}_i \otimes \overline{a}_j) \mid 1 \le i, j \le d \}$, since 
%%$c_s$ is a $k$-linear map. 
%We can also construct a $2$-cocycle $c$  by another method. 
%Let $S(a_1), \ldots, S(a_d) \in {\rm M}_n(\tilde{R})$ be lifts of $a_1,  \ldots, a_d$, respectively. 
%For $1 \le i, j \le d$, we can choose a lift $S(a_ia_j) \in \widetilde{R}S(a_1) \oplus \cdots 
%\oplus \widetilde{R}S(a_d)$ of $a_ia_j$.  
%Note that $S(a_ia_j) - S(a_i)S(a_j) \in {\rm M}_n(I) \cong {\rm M}_n(k)\otimes_{k} I$. 
%Set $c(\overline{a}_i \otimes \overline{a}_j) = (S(a_ia_j) - S(a_i)S(a_j) \mod A_{0}\otimes_{k} I) 
%\in ({\rm M}_n(I)/A_{0})\otimes_{k} I$. 
%We define $c : A_0\otimes_{k} A_0 \to ({\rm M}_n(k)/A_{0})\otimes_{k} I$ as the $k$-linear map 
%determined by $\overline{a}_i \otimes \overline{a}_j \mapsto c(\overline{a}_i \otimes \overline{a}_j)$ for $1 \le i, j %%%\le d$. 
%Then we can easily check that 
%$c$ is a $2$-cocycle and that $[c] = [c_{(s, S)}] \in H^2(A_0, ({\rm M}_n(k)/A_0)\otimes_{k} I)$ 
%for arbitrary lifts $s : R \to \widetilde{R}$ and $S :  {\rm M}_n(R) \to {\rm M}_n(\widetilde{R})$. 
%For defining the cohomology class $[c]$, 
%we only need to choose lifts $S(a_i)$ and $S(a_ia_j)$ of $a_i$ and $a_ia_j$ for $1 \le i, j \le d$, respectively.  
%\end{remark} 

\bigskip 

For lifts $S(a_i) \in {\rm M}_n(\widetilde{R})$ of $a_i$ ($1 \le i \le d$) and lifts 
$S(a_ia_j) \in \widetilde{R}S(a_1)\oplus \cdots \oplus \widetilde{R}S(a_d)$ of $a_ia_j$ ($1 \le i, j \le d$),  
it seems that the $2$-cocycle $c_{S}$ depends on the choice of  an $R$-basis 
$a_1, \ldots, a_{d}$ of $A$.     
In fact, we see that the cohomology class $[c_{S}]$ is independent from the choice of an $R$-basis of $A$ 
by the following lemma.  

\begin{lemma} 
The cohomology class $[c_{S}]$ is independent from the choice of an $R$-basis 
$\{ a_1,\ldots,a_{d} \}$ of $A$.   
%Let $c'' : A\otimes_{R}A \to ({\rm M}_n(k)/A_0)\otimes_{k} I$ be as above. 
%Then $c''$ does not depend on the choice of  a basis $a_1,\ldots,a_{n^2}$ of ${\rm M}_n(R)$ 
%such that $A = Ra_{1} \oplus\cdots\oplus Ra_{d}$. 
\end{lemma}

\proof 
Suppose that $a_1,\ldots,a_{d}$ and $b_1,\ldots,b_{d}$ are bases of 
$A$ over $R$.  Then there exists $P=(p_{ij})\in {\rm GL}_d(R)$ such that 
\[ (a_1,\ldots,a_d)=(b_1,\ldots,b_d)P. \] 
Let $\widetilde{P}=(\widetilde{p}_{ij}) \in {\rm GL}_d(\widetilde{R})$ be a lift of $P$.  
Let us choose lifts $S(a_i) \in {\rm M}_n(\widetilde{R})$ of $a_i$ ($1 \le i \le d$) and lifts 
$S(a_ia_j) \in \widetilde{R}S(a_1)\oplus \cdots \oplus \widetilde{R}S(a_d)$ of $a_ia_j$ ($1 \le i, j \le d$).   
We define lifts $T(b_1),\ldots,T(b_d)\in {\rm M}_n(\widetilde{R})$ of 
$b_1, \ldots, b_d \in {\rm M}_n(R)$ by 
\[ (T(b_1),\ldots,T(b_d))=(S(a_1),\ldots,S(a_d))\widetilde{P}^{-1}. \]
Then 
\[ S(a_i)=\sum_{j=1}^d T(b_j)\widetilde{p}_{ji} \qquad (1\le i\le d) \] holds.  
Take lifts $T(a_ia_j) \in \widetilde{R}T(a_1)\oplus \cdots \oplus \widetilde{R}T(a_d)$ of $a_ia_j$ for $1 \le i, j \le d$.  

We define 
$c_{S}' : A\otimes_RA\to {\rm M}_n(I)$ by
\[ c_{S}'\left(\sum_{i,j=1}^d r_{ij}a_i\otimes a_j \right)=
   \sum_{i,j=1}^d \widetilde{r}_{ij} (S(a_ia_j)-S(a_i)S(a_j)),   \]
where $\widetilde{r}_{ij} \in \widetilde{R}$ is a lift of $r_{ij}$.  
Similarly, we also define 
$c_{T}': A\otimes_RA\to {\rm M}_n(I)$ by 
\[ c_{T}' \left(\sum_{i,j=1}^d r_{ij}b_i\otimes b_j \right)=
   \sum_{i,j=1}^d \widetilde{r}_{ij} (T(b_ib_j)-T(b_i)T(b_j)). \] 
Let $c''_{S}, c''_{T} : A\otimes_RA\to ({\rm M}_n(k)/A_0)\otimes_{k} I$ be 
the compositions of ${\rm M}_n(I)\cong {\rm M}_n(k)\otimes_k I 
\to ({\rm M}_n(k)/A_0)\otimes_k I$ with $c_{S}'$ and $c_{T}'$, respectively. 
%Similarly, we also define $c''_t : A\otimes_RA\to ({\rm M}_n(k)/A_0)\otimes_k I$. 
By lemma \ref{lemma:independentlift}, we only need to show the claim that $c''_{S}=c''_{T}$. 

Assume that 
\[ a_ia_j=\sum_{k=1}^d\alpha^k_{ij}a_k,\qquad
   b_ib_j=\sum_{k=1}^d\beta^k_{ij}b_k\]
for $1\le i,j\le d$, where $\alpha^k_{ij}, \beta^k_{ij} \in R$.  
We can write 
\[ S(a_ia_j)=\sum_{k=1}^d \widetilde{\alpha}^k_{ij}S(a_k),\qquad
   T(b_ib_j)=\sum_{k=1}^d \widetilde{\beta}^k_{ij} T(b_k), \]
where $\widetilde{\alpha}^k_{ij}, \widetilde{\beta}^k_{ij} \in  \widetilde{R}$ are 
lifts of $\alpha^k_{ij}$ and $\beta^k_{ij}$, respectively. 
Since 
\[ \begin{array}{rcl}
    a_ia_j&=&{\displaystyle 
    \sum_{k=1}^d\alpha^k_{ij}a_k}\\[4mm]
    &=&{\displaystyle \sum_{k=1}^d\alpha^k_{ij}
        \sum_{l=1}^d b_lp_{lk}}\\[4mm]
    &=&{\displaystyle \sum_{k,l=1}^d \alpha^k_{ij}p_{lk}b_l}\\[4mm]
   \end{array} \]  and 
\[ \begin{array}{rcl}
     a_ia_j&=&{\displaystyle (\sum_{l_1=1}^db_{l_1}p_{l_1i})
     (\sum_{l_2=1}^db_{l_2}p_{l_2j})}\\[4mm]
     &=&{\displaystyle \sum_{l_1,l_2=1}^d p_{l_1i}p_{l_2j}
         b_{l_1}b_{l_2}}\\[4mm]
     &=&{\displaystyle \sum_{l_1,l_2=1}^d p_{l_1i}p_{l_2j}
         \sum_{l=1}^d \beta^l_{l_1l_2}b_l}\\[4mm]
     &=&{\displaystyle \sum_{l,l_1,l_2=1}^d
         p_{l_1i}p_{l_2j}\beta^l_{l_1l_2}b_l}, \\[4mm]
   \end{array}\] we have 
\begin{equation}\label{eq:ab-formulta} 
\sum_{k=1}^d \alpha^k_{ij}p_{lk}=
   \sum_{l_1,l_2=1}^d p_{l_1i}p_{l_2j}\beta^l_{l_1l_2} 
\end{equation} 
for $1\le i,j,l\le d$. 
Let us show that $c_{S}'(x)-c_{T}'(x)\in I\otimes_k A_0$ for any 
$x\in A\otimes_RA$. 
Let ${\displaystyle x=\sum_{i,j=1}^d r_{ij}a_i\otimes a_j} \in A\otimes_{R}A$.  Let $\widetilde{r}_{ij} \in \widetilde{R}$ be a lift  of $r_{ij} \in R$ ($1 \le i, j \le d$).  
Then we have  
%\[ \begin{array}{rcl}
\begin{eqnarray*} 
      {\displaystyle c_{S}'(x)} \;& = & \; {\displaystyle c_{S}'\left(\sum_{i,j=1}^d r_{ij}a_i\otimes a_j\right)} \\ 
      &=&   {\displaystyle \sum_{i,j=1}^d \widetilde{r}_{ij}(S(a_ia_j)-S(a_i)S(a_j))}\\%[4mm]
    &=&{\displaystyle \sum_{i,j=1}^d \widetilde{r}_{ij} \left( \sum_{k=1}^d\widetilde{\alpha}^k_{ij} S(a_k)
                     -S(a_i)S(a_j)\right)} \\%[4mm] 
    &=&{\displaystyle \left(\sum_{i,j,k=1}^d \widetilde{r}_{ij}\widetilde{\alpha}^k_{ij}S(a_k)\right)
                     -\left(\sum_{i,j=1}^d \widetilde{r}_{ij}S(a_i)S(a_j)\right)}\\%[4mm] 
    &=&{\displaystyle \left(\sum_{i,j,k=1}^d \widetilde{r}_{ij}\widetilde{\alpha}^k_{ij} 
                      \sum_{l=1}^dT(b_l)\widetilde{p}_{lk}\right)} %\\[4mm] 
    %&&
    {\displaystyle -\left(\sum_{i,j=1}^d \widetilde{r}_{ij} 
                       \left(\sum_{l_1=1}^d T(b_{l_1})\widetilde{p}_{l_1i}\right) 
                      \left(\sum_{l_2=1}^dT(b_{l_2})\widetilde{p}_{l_2j}\right)\right)}\\%[4mm] 
    &=&{\displaystyle \left(\sum_{i,j,k,l=1}^d \widetilde{r}_{ij}\widetilde{\alpha}^k_{ij}\widetilde{p}_{lk}  
                       T(b_l)\right)}%\\[4mm] 
    %&&
    {\displaystyle -\left(\sum_{i,j,l_1,l_2=1}^d
                       \widetilde{r}_{ij}\widetilde{p}_{l_1i}\widetilde{p}_{l_2j} 
                       T(b_{l_1})T(b_{l_2})\right)}. \\%[4mm]
\end{eqnarray*} 
%   \end{array}\]
On the other hand, 
\begin{eqnarray*} 
%\[ \begin{array}{rcl} 
   {\displaystyle c_{T}'(x)} \; & = & \; 
    {\displaystyle c_{T}'\left(\sum_{i,j=1}^dr_{ij}a_i\otimes a_j\right)} \\ 
    &=& 
    {\displaystyle c_{T}'\left(\sum_{i,j=1}^dr_{ij}\left(\sum_{l_1=1}^d b_{l_1}p_{l_1i}\right)\otimes
         \left(\sum_{l_2=1}^d b_{l_2}p_{l_2j}\right)\right)}\\%[4mm]
    &=&{\displaystyle \sum_{l_1,l_2=1}^d \sum_{i, j=1}^{d} \widetilde{r}_{ij}\widetilde{p}_{l_1i}\widetilde{p}_{l_2j} 
                      (T(b_{l_1}b_{l_2})-T(b_{l_1})T(b_{l_2}))}\\%[4mm]
    &=&{\displaystyle \sum_{l_1,l_2=1}^d \sum_{i, j=1}^{d} \widetilde{r}_{ij}\widetilde{p}_{l_1i}\widetilde{p}_{l_2j}
                      \left(\sum_{l=1}^d \widetilde{\beta}^l_{l_1l_2}T(b_l)
                      -T(b_{l_1})T(b_{l_2})\right)}\\%[4mm]
%    &=&{\displaystyle \sum_{l_1,l_2=1}^d s(\sum_{i, j=1}^{d} r_{ij}p_{l_1i}p_{l_2j})
%                     (\sum_{l=1}^ds(\beta^l_{l_1l_2})
 %                      T(b_l)-T(b_{l_1})T(b_{l_2}))}\\[4mm]
    &=&{\displaystyle \left(\sum_{i, j, l,l_1,l_2=1}^d
                      \widetilde{r}_{ij}\widetilde{p}_{l_1i}\widetilde{p}_{l_2j}
                      \widetilde{\beta}^l_{l_1l_2}T(b_l)\right)}%\\[4mm]
    %&&
    {\displaystyle  -\left(\sum_{i, j, l_1,l_2=1}^d \widetilde{r}_{ij}\widetilde{p}_{l_1i}\widetilde{p}_{l_2j}
                      T(b_{l_1})T(b_{l_2})\right)}. \\ 
%   \end{array}\]
\end{eqnarray*} 
By (\ref{eq:ab-formulta}), 
\[ \sum_{i,j,k=1}^d \widetilde{r}_{ij}\widetilde{\alpha}^k_{ij}\widetilde{p}_{lk} 
   - \sum_{i, j, l_1,l_2=1}^d \widetilde{r}_{ij}\widetilde{p}_{l_1i}\widetilde{p}_{l_2j}\widetilde{\beta}^l_{l_1l_2} \in I \]
for $1\le l\le d$. 
Hence 
\[ \begin{array}{rl}
c_{S}'(x) -c_{T}'(x)  =  &{\displaystyle \sum_{l=1}^d \left(\sum_{i,j,k=1}^d \widetilde{r}_{ij}\widetilde{\alpha}^k_{ij}\widetilde{p}_{lk} - \sum_{i, j, l_1,l_2=1}^d \widetilde{r}_{ij}\widetilde{p}_{l_1i}\widetilde{p}_{l_2j}\widetilde{\beta}^l_{l_1l_2} \right)T(b_{l})}\\[5mm] 
   =&{\displaystyle \sum_{l=1}^d \left(\sum_{i,j,k=1}^d \widetilde{r}_{ij}\widetilde{\alpha}^k_{ij}\widetilde{p}_{lk} - \sum_{i, j, l_1,l_2=1}^d \widetilde{r}_{ij}\widetilde{p}_{l_1i}\widetilde{p}_{l_2j}\widetilde{\beta}^l_{l_1l_2} \right)\overline{b}_{l}}  
   \in I\otimes_k A_0,  
   \end{array} \]
where $\overline{b}_l=(b_l \mod m)$.  
Therefore, $c''_{S}=c''_{T}$.  
\qed 

%%%%%%%%%%%%%%%

\bigskip 

By the lemmas above, we have a unique cohomology class $[c] \in H^{2}(A_{0}, ({\rm M}_n(k)/A_{0})\otimes_{k} I)$ for $A \in {\rm Mold}_{n, d}(R)$ and $(\widetilde{R}, \widetilde{m}, k)$. Here we introduce the following definition: 

\begin{definition}\rm 
We call $[c] \in H^{2}(A_{0}, ({\rm M}_n(k)/A_{0})\otimes_{k} I)$ 
the {\it cohomology class} defined by $A$ and $(\widetilde{R}, \widetilde{m}, k)$.  
\end{definition}

\begin{proposition}\label{prop:cohomogyzerolift}
Let $(R, m, k)$, $(\widetilde{R}, \widetilde{m}, k)$, and $I$ be as above. 
Let $A \in {\rm Mold}_{n, d}(R)$ and $A_0=A\otimes_{R} k$. 
There exists $\widetilde{A} \in {\rm Mold}_{n, d}(\widetilde{R})$ such that 
$\widetilde{A}\otimes_{\widetilde{R}} R = A$ if and only if 
the cohomology class $[c]$ defined by $A$ and $(\widetilde{R}, \widetilde{m}, k)$ 
is zero in $H^2(A_0, ({\rm M}_n(k)/A_0)\otimes_{k} I)$. 
\end{proposition}

\proof 
Assume that there exists $\widetilde{A} \in {\rm Mold}_{n, d}(\widetilde{R})$ such that 
$\widetilde{A}\otimes_{\widetilde{R}} R = A$. 
For a basis $a_1, a_2, \ldots, a_{d}$ of $A$ over $R$, there exists a basis 
$\widetilde{a}_1, \widetilde{a}_2, \ldots, \widetilde{a}_d$ of $\widetilde{A}$ over 
$\widetilde{R}$ such that  $\pi(\widetilde{a}_i) = a_i$ for $i=1, 2, \ldots, d$, where 
$\pi : {\rm M}_n(\widetilde{R}) \to {\rm M}_n(R)$ is the projection. 
Set $S(a_i) = \widetilde{a}_i$ for $1 \le i \le d$ and  
$S(a_ia_j) = \widetilde{a}_i\widetilde{a}_j \in \widetilde{R} S(a_1) \oplus \cdots \oplus \widetilde{R} S(a_d) = \widetilde{A}$ for $1 \le i, j \le d$.  Put $\overline{a}_i = (a_i \mod m) \in A_0$ for $1 \le i \le d$. 
The $2$-cocycle $c$ defined by lifts $\{ S(a_i) \}$ and $\{ S(a_ia_j) \}$ satisfies 
$c(\overline{a}_i\otimes \overline{a}_j) = S(a_ia_j)-S(a_i)S(a_j) = \widetilde{a}_i \widetilde{a}_j - \widetilde{a}_i \widetilde{a}_j = 0$. Hence the cohomology class $[ c ]$ is zero in $H^2(A_0, ({\rm M}_n(k)/A_0)\otimes_{k} I)$. 

Conversely, assume that the cohomology class $[c]$ is zero in $H^2(A_0, ({\rm M}_n(k)/A_0)\otimes_{k} I)$. 
For an $R$-basis $a_1, \ldots, a_d$ of $A$, 
choose lifts $S(a_i) \in {\rm M}_n(\widetilde{R})$ of $a_i$ ($1 \le i \le d$) and lifts 
$S(a_ia_j) \in \widetilde{R}S(a_1)\oplus \cdots \oplus \widetilde{R}S(a_d)$ of $a_ia_j$ ($1 \le i, j \le d$), respectively.    
The $R$-linear map $c'_{S} : A\otimes_{R} A \to {\rm M}_n(k)\otimes_{k} I$ defined by $a_i \otimes a_j \mapsto S(a_ia_j)-S(a_i)S(a_j)$ for $1 \le i, j \le d$ induces a $2$-coboundary $c_S : A_0 \otimes A_0 \to ({\rm M}_n(k)/A_0)\otimes_{k}I$ by the assumption.   In other words, there exists a $k$-linear map 
$\theta : A_0 \to ({\rm M}_n(k)/A_0)\otimes_{k}I$ 
such that $c_S = d^1(\theta)$.  
Let us denote by 
$\theta'' : A \to A_0 \stackrel{\theta}{\to} ({\rm M}_n(k)/A_0)\otimes_{k} I$ the composition of 
$\theta$ with the projection $A \to A_0$. Choose an $R$-linear map 
$\theta' : A \to {\rm M}_n(I)$ as a lift of $\theta''$. 
Then there exist  
$t_{ij}^{l} \in I$ such that 
\begin{eqnarray}\label{eq:Saiaj} 
S(a_ia_j) - S(a_i)S(a_j) & = &  S(a_i)\theta'(a_j) - \theta'(a_ia_j) +\theta'(a_i)S(a_j) - \sum_{l=1}^{d} t_{ij}^{l} S(a_l) 
\end{eqnarray}
%(\ref{eq:Saiaj}) 
in ${\rm M}_n(I)$ for $1 \le i, j \le d$. 
Put $\widetilde{a}_i = S(a_i)+\theta'(a_i) \in {\rm M}_n(\widetilde{R})$ and 
$\widetilde{A} = \sum_{i=1}^{d} \widetilde{R}\widetilde{a}_{i} \subset {\rm M}_n(\widetilde{R})$.   
It is easy to see that 
$\widetilde{A}$ and ${\rm M}_n(\widetilde{R})/\widetilde{A}$ are free modules 
over $\widetilde{R}$ and that ${\rm rank}_{\widetilde{R}} \widetilde{A} = d$. 
Let us show that $\widetilde{A}$ is an $\widetilde{R}$-subalgebra of ${\rm M}_n(\widetilde{R})$. 
For $1 \le i, j \le d$, we can write $S(a_ia_j) = \sum_{l=1}^{d} \widetilde{c}_{ij}^{l} S(a_l)$ for some $\widetilde{c}_{ij}^{l} \in \widetilde{R}$. 
By using (\ref{eq:Saiaj}), we have 
\begin{eqnarray*}
\widetilde{a}_i\widetilde{a}_j & = & (S(a_i)+\theta'(a_i))(S(a_j)+\theta'(a_j)) \\ 
 & = & S(a_i)S(a_j)+ S(a_i)\theta'(a_j)+\theta'(a_i)S(a_j) \\ 
 & = & S(a_ia_j) + \theta'(a_ia_j) + \sum_{l=1}^{d} t_{ij}^{l} S(a_l) \\
 & = & \sum_{l=1}^{d} \widetilde{c}_{ij}^{l}S(a_l) +  \sum_{l=1}^{d} \widetilde{c}_{ij}^{l}\theta'(a_l) 
 + \sum_{l=1}^{d} t_{ij}^{l} S(a_l) \\ 
 & = & \sum_{l=1}^{d} (\widetilde{c}_{ij}^{l} + t_{ij}^{l})(S(a_l)+\theta'(a_l)) \\ 
 & = & \sum_{l=1}^{d} (\widetilde{c}_{ij}^{l} + t_{ij}^{l}) \widetilde{a}_l \in \widetilde{A}  
\end{eqnarray*}  
for $1 \le i, j \le d$.
Thus, $\widetilde{A}$ is closed under multiplication. 
Since $1 \in A = \oplus_{i=1}^{d} Ra_i$, we can write $1 = \sum_{i=1}^{d} r_i a_i$ for some 
$r_i \in R$. Take a lift $\widetilde{r}_{i} \in \widetilde{R}$ of $r_{i}$ for $1 \le i \le d$, respectively.   
Put $\widetilde{a} = \sum_{i=1}^{d} \widetilde{r}_i\widetilde{a}_i = \sum_{i=1}^{d} \widetilde{r}_i(S(a_i) + \theta'(a_i)) 
\in \widetilde{A}$. 
Then $\pi(\widetilde{a})=\pi(\sum_{i=1}^{d} \widetilde{r}_i\widetilde{a}_i)=\sum_{i=1}^{d} r_i a_i = 1$ and there exists $x \in{\rm M}_n(I)$ such that 
$\widetilde{a} = 1 + x$.  Hence $2\widetilde{a}-\widetilde{a}^2 = 2(1+x) - (1+x)^2 = 1 \in  \widetilde{A}$. 
Therefore, $\widetilde{A}$ is an $\widetilde{R}$-subalgebra of ${\rm M}_n(\widetilde{R})$. 
Obviously, $\widetilde{A}\otimes_{\widetilde{R}} R = A$. Thus, we have proved the statement. 
\qed 

%{\color{red} 
%\begin{theorem}
%Let $x \in {\rm Mold}_{n, d}$.  
%Let ${\mathcal A}$ be the universal mold on ${\rm Mold}_{n, d}$. 
%Set ${\mathcal A}(x) = {\mathcal A}\otimes_{{\mathcal O}_{{\rm Mold}_{n, d}}} k(x)$. 
%Then $x$ is a smooth point of ${\rm Mold}_{n, d}$ if and only if 
%the cohomology class $c({\mathcal A}(x))$ is zero in 
%$H^2({\mathcal A}(x), {\rm M}_n(k)/{\mathcal A}(x))$. 
%\end{theorem}
%}

\begin{theorem}\label{th:smooth} 
Let $x \in {\rm Mold}_{n, d}$.  
Let ${\mathcal A}$ be the universal mold on ${\rm Mold}_{n, d}$. 
Set ${\mathcal A}(x) = {\mathcal A}\otimes_{{\mathcal O}_{{\rm Mold}_{n, d}}} k(x)$.  
If $H^2({\mathcal A}(x), {\rm M}_n(k(x))/{\mathcal A}(x)) = 0$, then the canonical morphism 
${\rm Mold}_{n, d} \to {\Bbb Z}$ is smooth at $x$. 
\end{theorem} 

\proof
Since ${\rm Mold}_{n, d} \to {\Bbb Z}$ is a morphism of finite type of noetherian schemes, 
it suffices to show that if 
\[
\begin{array}{ccc}
{\rm Spec}\: R/I & \stackrel{f}{\to} & {\rm Mold}_{n, d} \\
\!\!\!\pi \downarrow & & \downarrow \\
{\rm Spec}\: R & \to & {\rm Spec}\: {\Bbb Z}  
\end{array}
\] 
is a commutative diagram such that 
$I$ is an ideal of an Artin local ring $(R, m, k)$ with $mI=0$,  
$f$ maps the special point $m/I$ to $x$, and $k(x)=R/m=k$, then there exists $g : {\rm Spec}\: R \to {\rm Mold}_{n, d}$ such that  
$g\circ \pi = f$ (for details, see \cite[Lemma~02HX]{Stacks}). Let $\overline{A}  \subset {\rm M}_n(R/I)$ be the mold corresponding to $f$. 
The cohomology class $[c]$ defined by $\overline{A}$ and $(R, m, k)$ is zero, 
since $\overline{A}\otimes_{R/I} k  
= {\mathcal A}(x) \subset 
{\rm M}_n(k)={\rm M}_n(k(x))$ and $H^2({\mathcal A}(x), ({\rm M}_n(k(x))/{\mathcal A}(x)) \otimes_{k(x)} I) = H^2({\mathcal A}(x), {\rm M}_n(k(x))/{\mathcal A}(x))\otimes_{k(x)} I  = 0$. 
By Proposition \ref{prop:cohomogyzerolift}, $\overline{A}$ 
has a lift $A \subset {\rm M}_n(R)$, and hence 
we have $g : {\rm Spec}\: R \to {\rm Mold}_{n, d}$ corresponding to $A$ such that $g\circ \pi = f$. 
This completes the proof. 
\qed 

\begin{remark}\rm%[{\cite[Remark~30]{Nakamoto-Torii:51st}}]
Even if $H^2({\mathcal A}(x), {\rm M}_n(k(x))/{\mathcal A}(x)) \neq 0$, 
the morphism ${\rm Mold}_{n, d} \to {\rm Spec}\: {\Bbb Z}$ may be smooth at $x \in {\rm Mold}_{n, d}$. 
Indeed, assume that ${\mathcal A}(x)= {\rm J}_n(k(x))$ for $x \in {\rm Mold}_{n, n}$ with $n \ge 2$, where ${\rm J}_n$ is defined in Definition~\ref{def:xofJn} below. We will see that $H^{2}({\rm J}_n(k(x)), {\rm M}_n(k(x))/{\rm J}_n(k(x))) \neq 0$ by Corollary~\ref{cor:hochschildjn}. However, $x$ is contained in an open subscheme ${\rm Mold}_{n, n}^{\rm reg}$ of ${\rm Mold}_{n, n}$ and ${\rm Mold}_{n, n}^{\rm reg}$ is smooth over ${\Bbb Z}$ 
(for details, see \cite{Nakamoto-Torii:classification}). 
%For example, 
\end{remark} 

\subsection{Smoothness of the morphism 
$\phi_{{\mathcal A}} : {\rm PGL}_{n, S} \to {\rm Mold}_{n, d}\otimes_{{\Bbb Z}} S$}\label{subsection:3.4} 

Let $S$ be a locally noetherian scheme. 
For a rank $d$ mold ${\mathcal A}$ of degree $n$ on $S$, 
we obtain a morphism $\tau_{{\mathcal A}} : S \to {\rm Mold}_{n, d}\otimes_{{\Bbb Z}} S$: 
\[  
\begin{array}{ccccc}
S & \stackrel{\tau_{{\mathcal A}}}{\to} & {\rm Mold}_{n, d}\otimes_{{\Bbb Z}} S & \to & {\rm Mold}_{n, d} \\ 
 & \searrow & \downarrow & & \downarrow \\ 
  & {id_{S}} & S & \to & {\rm Spec}\: {\Bbb Z}.  
\end{array} 
\]  
Let us consider the group scheme ${\rm PGL}_{n, S} = {\rm PGL}_n\otimes_{{\Bbb Z}} S$ over $S$. 
We define the $S$-morphism $\phi_{{\mathcal A}} : {\rm PGL}_{n, S} \to {\rm Mold}_{n, d}\otimes_{{\Bbb Z}} S$ by 
$P \mapsto P^{-1}{\mathcal A}P$. 
For an $S$-scheme $f : X \to S$, set ${\mathcal A}_{X} = f^{\ast}{\mathcal A} \subseteq {\rm M}_n({\mathcal O}_X)$. 
In particular, set $A_{R} = f^{\ast}{\mathcal A} \subseteq {\rm M}_n(R)$ in the case $X = {\rm Spec} R$. 
For an $X$-valued point $P$ of ${\rm PGL}_{n, S}$, $\phi_{{\mathcal A}}(P) = P^{-1}{\mathcal A}_{X}P$. 

Let us consider the question whether $\phi_{{\mathcal A}} : {\rm PGL}_{n, S} \to {\rm Mold}_{n, d}\otimes_{{\Bbb Z}} S$ is  (formally) smooth or not. Let $I$ be an ideal of an Artin local ring $(R, m, k)$ with $mI=0$. 
Assume that 
\begin{eqnarray}\label{diagram:phiasmooth} 
\begin{array}{ccc}
{\rm Spec} R/I & \stackrel{\overline{g}}{\to} & {\rm PGL}_{n, S} \\
\iota \downarrow & & \downarrow {\phi_{\mathcal A}} \\
{\rm Spec} R & \stackrel{\psi}{\to} & {\rm Mold}_{n, d}\otimes_{{\Bbb Z}} S \\
\end{array} 
\end{eqnarray} 
is a commutative diagram. If there exists $g : {\rm Spec} R \to {\rm PGL}_{n, S}$ such that 
$\overline{g} = g \circ \iota$ and $\psi = \phi_{{\mathcal A}} \circ g$ for any commutative diagram above, then 
$\phi_{{\mathcal A}}$ is smooth since $\phi_{{\mathcal A}}$ is locally of finite type (for details, see \cite[Lemma~02HX]{Stacks}).  

Denote by $B' (\subseteq {\rm M}_n(R))$ the mold associated to $\psi$.  We can identify  
${\rm PGL}_n(R) = {\rm GL}_n(R)/(R^{\times}\cdot I_n)$ with the set of $R$-valued points of the group scheme 
${\rm PGL}_n$ for a local ring $R$. Note that there is a commutative diagram consisting of surjective group homomorphisms: 
\[
\begin{array}{ccc}
{\rm GL}_n(R) & \twoheadrightarrow & {\rm GL}_n(R/I) \\
 \rotatebox{270}{\hspace*{-1.5ex}$\twoheadrightarrow$} & & \rotatebox{270}{\hspace*{-1.5ex}$\twoheadrightarrow$} \\
{\rm PGL}_n(R) & \twoheadrightarrow & {\rm PGL}_n(R/I). \\
\end{array} 
\]
Given diagram (\ref{diagram:phiasmooth}), we have $\overline{P} \in {\rm GL}_n(R/I)$ such that 
$\overline{P}^{-1} A_{R/I} \overline{P} = B'\otimes_{R} (R/I)$. 
There exists 
$g : {\rm Spec} R \to {\rm PGL}_{n, S}$ satisfying $\overline{g} = g \circ \iota$ 
and $\psi = \phi_{{\mathcal A}} \circ g$ if and only if there exists 
$P \in {\rm GL}_n(R)$ such that $P^{-1}A_{R} P = B'$ and $(P \mod I) = \overline{P} \in {\rm GL}_n(R/I)$. 
Take a lift $P' \in {\rm GL}_n(R)$ of $\overline{P}$. Set $B = P' B' P'^{-1} \subseteq {\rm M}_n(R)$. 
Then $B \otimes_{R} (R/I) = \overline{P} (B'\otimes_{R} (R/I)) \overline{P}^{-1} = A_{R/I}$. Let us denote by 
$\overline{I}_n \in {\rm GL}_n(R/I)$ the identity matrix. 
There exists $P'' \in {\rm GL}_n(R)$ such that $P''^{-1}A_{R}P'' = B$ and $(P'' \mod I) = \overline{I}_n \in {\rm GL}_n(R/I)$ if and only if $g$ exists. 
Hence, we only need to consider whether there exists $P \in {\rm GL}_n(R)$ such that $P^{-1}A_{R} P = B$ and $(P \mod I) = \overline{I}_n \in {\rm GL}_n(R/I)$ 
when $B \otimes_{R} (R/I) = A_{R/I}$. Here we have: 

\begin{lemma}\label{lemma:pmodiIncase}
Let ${\mathcal A}$ be a rank $d$ mold of degree $n$ on a locally noetherian scheme $S$.  
The morphism $\phi_{{\mathcal A}} : {\rm PGL}_{n, S} \to {\rm Mold}_{n, d}\otimes_{{\Bbb Z}} S$ defined by $P\mapsto P^{-1}{\mathcal A}P$ is smooth if and only if 
for any ideal $I$ of an Artin local ring $(R, m, k)$ over $S$ with $mI=0$, and for any rank $d$ mold $B \subset {\rm M}_n(R)$ with $B \otimes_{R} (R/I) = A_{R/I}$, there exists $P \in {\rm GL}_n(R)$ such that $P^{-1}A_{R}P =B$ and 
$(P \mod I) = \overline{I}_n \in {\rm GL}_n(R/I)$.   

\end{lemma} 

Assume that $B \otimes_{R} (R/I) = A_{R/I}$.  
Let us define $\delta': A_{R/I} = A_{R}\otimes_{R}(R/I) \to ({\rm M}_n(R)/A_{R})\otimes_{R} I$. Take $R$-bases $a_1, a_2, \ldots, a_d \in A_R$ and 
$b_1, b_2, \ldots, b_d \in B$ such that $\overline{a}_i = \overline{b}_i$ in ${\rm M}_n(R/I)$ for $1 \le i  \le d$, where $\overline{a}_i = (a_i \mod I)$ and $\overline{b}_i = (b_i \mod I)$.   This is possible because $A_R$ and $B$ are free modules over a local ring $R$ and $B \otimes_{R} (R/I) = A_{R/I}$. 
Set $c_i = b_i -a_i \in {\rm M}_n(I)$ for $1 \le i \le d$. 
For $\sum_{i=1}^{d} \overline{r}_i \overline{a}_i \in A_{R/I}$, we define  
\[
\delta'\left(\sum_{i=1}^{d} \overline{r}_i \overline{a}_i \right) = \left(\sum_{i=1}^{d} r_i c_i \mod A_{R}\otimes_{R} I \right) \in ({\rm M}_n(R)/A_{R})\otimes_{R} I = {\rm M}_n(I)/(A_{R}\otimes_{R}I), 
\]  
where $r_i \in R$ and $\overline{r}_i = (r_i \mod I) \in R/I$ for $1 \le i \le d$.  
Here note that $\sum_{i=1}^{d} r_i c_i \in {\rm M}_n(I)$.

First, we show that $\delta'$ does not depend on the choice of lifts $r_i$ of $\overline{r}_i$.  
Let $r'_i \in R$ be another lift of $\overline{r}_i$ for $1 \le i \le d$. Set $s_i = r_i - r'_i \in I$. %for $1 \le i \le d$. 
Note that $s_i c_i = 0$ for $1 \le i \le d$ because $I^2=0$. 
Since 
\[
\sum_{i=1}^{d} r_i c_i = \sum_{i=1}^{d} (r'_i +s_i) c_i = \sum_{i=1}^{d} r'_i c_i + \sum_{i=1}^{d} s_i c_i = \sum_{i=1}^{d} r'_i c_i, 
\]
we see that $\delta'$ does not depend on the choice of lifts $r_i$ of $\overline{r}_i$.  

Since $I^2=0$, $({\rm M}_n(R)/A_{R})\otimes_{R} I$ is an $R/I$-module. It is obvious that $\delta' : 
A_{R/I} \to ({\rm M}_n(R)/A_{R})\otimes_{R} I$ is an $R/I$-linear map. 
Second, we show that $\delta'$ does not depend on the choice of $b_i$ for $1 \le i \le d$. 
Let $b'_1, b'_2, \ldots, b'_d$ be another basis of $B \subseteq {\rm M}_n(R)$ such that $b'_i -b_i \in {\rm M}_n(I)$ for $1 \le i \le n$.  We can write $b'_i = b_i + \sum_{j=1}^{d} x_{ij} b_j$ for $x_{ij} \in I$. 
By using $c_i = b_i -a_i \in {\rm M}_n(I)$ and $I^2=0$, we have 
\[ 
b'_i = b_i + \sum_{j=1}^{d} x_{ij} b_j = b_i + \sum_{j=1}^{d} x_{ij} (a_j+c_j) = b_i + \sum_{j=1}^{d} x_{ij} a_j.  
\]
Set $c'_i = b'_i-a_i \in {\rm M}_n(I)$ for $1 \le i \le d$. Since $c'_i = b'_i - a_i = b_i-a_i+\sum_{j=1}^{d}x_{ij}a_{j} = c_i +\sum_{j=1}^{d}x_{ij}a_{j}$,  
\[ 
\sum_{i=1}^{d} r_i c'_i  =  \sum_{i=1}^{d} r_ic_i + \sum_{i=1}^{d} \sum_{j=1}^{d} r_i x_{ij} a_{j}.  
\] 
Hence $(\sum_{i=1}^{d} r_i c_i \mod A_{R}\otimes_{R} I ) = (\sum_{i=1}^{d} r_i c'_i \mod A_{R}\otimes_{R} I )$, which implies that $\delta'$ does not depend on the choice of a basis $\{ b_1, b_2, \ldots, b_d \}$ of $B$ such that  $\overline{a}_i = \overline{b}_i$ in ${\rm M}_n(R/I)$ for $1 \le i  \le d$.  

Third, we show that $\delta'$ does not depend on the choice of a basis $\{ a_1, a_2, \ldots, a_d \}$ of $A_{R}$ over $R$.
Let $\{ a'_1, a'_2, \ldots, a'_d \}$ be another basis of $A_{R}$ over $R$. 
There exists $P=(p_{ij}) \in {\rm GL}_n(R)$ such that $(a_1, a_2, \ldots, a_d) = (a'_1, a'_2, \ldots, a'_d)P$. 
Let $\{ b_1, b_2, \ldots, b_d \}$ and $\{ b'_1, b'_2, \ldots, b'_d \}$ be bases of $B$ over $R$ such that 
$(b_1, b_2, \ldots, b_d)= (b'_1, b'_2, \ldots, b'_d) P$ and  $c_i=b_i - a_i, c'_i=b'_i - a'_i \in {\rm M}_n(I)$ for $1 \le i \le d$. 
Then $(c_1, c_2, \ldots, c_d)= (c'_1, c'_2, \ldots, c'_d) P$. 
By using $a_i = \sum_{j=1}^{d} p_{ji}a'_j$, we have 
$\sum_{i=1}^{d} r_i a_i = \sum_{i=1}^{d} r_i (\sum_{j=1}^{d} p_{ji} a'_j) = \sum_{j=1}^{d} (\sum_{i=1}^{d} r_i p_{ji}) a'_j$.   
Similarly, by using $c_i = \sum_{j=1}^{d} p_{ji}c'_j$, we obtain  $\sum_{i=1}^{d} r_i c_i = \sum_{i=1}^{d} r_i (\sum_{j=1}^{d} p_{ji} c'_j) = \sum_{j=1}^{d} (\sum_{i=1}^{d} r_i p_{ji}) c'_j$. 
Hence we see that $\delta'$ does not depend on the choice of a basis $\{ a_1, a_2, \ldots, a_d \}$ of $A_{R}$ over $R$.  

Since $I^2=0$, we can regard $({\rm M}_n(R)/A_{R})\otimes_{R} I$ as a bimodule over $A_{R/I} = A_{R}\otimes_{R} (R/I)$.  Fourth, we show that the $R/I$-linear  map $\delta' : A_{R/I} \to ({\rm M}_n(R)/A_{R})\otimes_{R} I$ is a derivation. For proving it, we only need to verify that $\delta'(\overline{a}_i\overline{a}_j) = \delta'(\overline{a}_i)\overline{a}_j + \overline{a}_i\delta'(\overline{a}_j)$ for $1 \le i, j \le d$.   Let $a_i a_j = \sum_{k=1}^{d} c_{ij}^{k} a_{k}$. 
Because $\overline{a}_i = \overline{b}_i$ and $\overline{a}_i \overline{a}_j = \overline{b}_i \overline{b}_j$ in ${\rm M}_n(R/I)$, there exist 
$d_{ij}^{k} \in I$ such that $b_i b_j = \sum_{k=1}^{d} (c_{ij}^{k} + d_{ij}^{k})b_{k}$ for $1 \le i, j  \le d$.  
By the definition of $\delta'$, 
\[ 
\delta'(\overline{a}_i)\overline{a}_j + \overline{a}_i\delta'(\overline{a}_j)  =  \left( c_i a_j + a_ic_j \mod A_{R}\otimes_{R} I \right). 
\]   
By using $c_ic_j=0$ for $c_i = b_i -a_i \in {\rm M}_n(I)$, we see that 
\begin{eqnarray*}
c_ia_j+a_ic_j & = & (b_i-a_i) a_j + (b_i - c_i) c_j \\ 
 & = &  b_ia_j-a_ia_j + b_ic_j \\ 
 & = & b_i(b_j-c_j)-a_ia_j + b_ic_j \\ 
 & = & b_ib_j-a_ia_j \\
  & = & \sum_{k=1}^d (c_{ij}^{k}+d_{ij}^{k}) b_k -\sum_{k=1}^{d}c_{ij}^{k} a_k \\
  & = & \sum_{k=1}^d c_{ij}^{k} c_k + \sum_{k=1}^{d}d_{ij}^{k} b_k \\
\end{eqnarray*} 
in ${\rm M}_n(I) \subset {\rm M}_n(R)$. 
Since $d_{ij}^{k} c_k = 0$, we have 
$d_{ij}^{k}b_k = d_{ij}^{k}(a_k+c_k) = 
d_{ij}^{k}a_k \in A_{R}\otimes_{R} I$.  Thereby,  
\[ \delta'(\overline{a}_i)\overline{a}_j + \overline{a}_i\delta'(\overline{a}_j) = 
\left( \sum_{k=1}^{d} c_{ij}^{k} c_k \mod A_{R}\otimes_{R} I \right) 
\] 
in $({\rm M}_n(R)/A_{R}) \otimes_{R} I = {\rm M}_n(I)/(A_{R}\otimes_{R}I)$. 
On the other hand, 
\begin{eqnarray*}
\delta'(\overline{a}_i\overline{a}_j) & = & \delta'(\sum_{k=1}^{d} \overline{c}_{ij}^{k} \overline{a}_k ) \\ 
 & = & \left( \sum_{k=1}^{d} c_{ij}^{k} c_k  \mod A_{R}\otimes_{R} I \right),   \\ 
\end{eqnarray*} 
where $\overline{c}_{ij}^{k} = (c_{ij}^{k} \mod I) \in R/I$. 
Hence $\delta'(\overline{a}_i\overline{a}_j) = \delta'(\overline{a}_i)\overline{a}_j + \overline{a}_i\delta'(\overline{a}_j)$. 

Let us consider the $d$-dimensional subalgebra $A_0 = A_{R} \otimes_{R} k \subseteq {\rm M}_n(k)$. 
Since $I$ is finitely generated over $R$ and $mI=0$, $I$ is a finite-dimensional $k$-vector space. Then  
\begin{eqnarray*} 
({\rm M}_n(R)/A_{R})\otimes_{R} I & = & ({\rm M}_n(R)/A_{R})\otimes_{R}  (R/m) \otimes_{R} I \\ 
 & = & ({\rm M}_n(R)/A_{R})\otimes_{R} k \otimes_{k} k \otimes_{R} I \\
 & = & ({\rm M}_n(k)/A_0) \otimes_{k} I.  
\end{eqnarray*} 
It is easy to see that the derivation $\delta' : A_{R/I} \to ({\rm M}_n(R)/A_{R})\otimes_{R} I = ({\rm M}_n(k)/A_0) \otimes_{k} I$ factors through 
$A_{R/I}\otimes_{R/I} k = A_0$. Hence we obtain a $k$-linear map 
$\delta : A_0 \to ({\rm M}_n(k)/A_0) \otimes_{k} I$. Set $[a_i] = (a_i \mod A_{R}\otimes_{R}m) \in A_{0} = A_{R}\otimes_{R} k = A_{R}/(A_{R}\otimes_{R}m)$. Note that $\delta([a_i]) = (c_i \mod A_{R}\otimes_{R} I) \in {\rm M}_n(I)/(A_{R}\otimes_{R}I) = ({\rm M}_n(k)/A_0)\otimes_{k} I$ for $1 \le i \le d$.  
We regard $({\rm M}_n(k)/A_0) \otimes_{k} I$ as an $A_0$-bimodule by $a(\overline{X}\otimes x) b = \overline{aXb}\otimes x$ for $a, b \in A_0$, $X \in {\rm M}_n(k)$, and $x \in I$, where  $\overline{X} = (X \mod A_0)$ and 
$\overline{aXb} = (aXb \mod A_0)$ in ${\rm M}_n(k)/A_0$.  We easily see that $\delta$ is a derivation, that is, 
$\delta(ab) = \delta(a)b+a\delta(b)$ for $a, b \in A_0$.   Then $\delta$ is a $1$-cocycle in $C^{1}(A_0, ({\rm M}_n(k)/A_0)\otimes_{k} I)$. 

\begin{definition}\rm
Let $\delta : A_0 \to ({\rm M}_n(k)/A_0)\otimes_{k} I$ be as above.  
We say that $\delta$ is the {\it derivation} associated to a rank $d$ mold $B \subset {\rm M}_n(R)$ with $B \otimes_{R} (R/I) = A_{R/I}$. Note that $\delta$ does not depend on the choice of an $R$-basis $\{ a_i \}$ of $A$ and an  $R$-basis $\{ b_i \}$ of $B$ with $a_i - b_i \in {\rm M}_n(I)$ for $1 \le i \le d$.  
\end{definition}

Let $x_1, x_2, \ldots, x_r \in I$ be a basis of $I$ over $k$. 
We can write $\delta = \sum_{i=1}^{r} x_i \delta_i$, where $\delta_i \in C^{1}(A_0, {\rm M}_n(k)/A_0)$ is a derivation for 
$1 \le i \le r$. The cohomology class $[ \delta ]$ is expressed by 
$[ \delta ] = \sum_{i=1}^{r} [\delta_i ] \otimes x_i$ in $H^{1}(A_0, ({\rm M}_n(k)/A_0)\otimes_{k} I) = H^{1}(A_0, {\rm M}_n(k)/A_0)\otimes_{k} I$.  

\begin{lemma}\label{lemma:eqdelta=0} %\ref{diagram:phiasmooth}
In the situation above, there exists $P \in {\rm GL}_n(R)$ such that $P^{-1}A_R P =B$ and $(P \mod I) = \overline{I}_n \in {\rm GL}_n(R/I)$ if and only if $[ \delta ] =0$ in $H^{1}(A_0, ({\rm M}_n(k)/A_0)\otimes_{k} I)$. 
\end{lemma} 

\proof 
Suppose that $[\delta]=0$.  For $1 \le i \le r$, $[\delta_i] = 0$ in $H^{1}(A_0, {\rm M}_n(k)/A_0)$. 
By the exact sequence  
\[ 
0 \to N(A_0)/A_0 \to {\rm M}_n(k)/A_0 \stackrel{\overline{d}}{\to} {\rm Der}_{k}(A_0, {\rm M}_n(k)/A_0) \to H^{1}(A_0, {\rm M}_n(k)/A_0) \to 0 
\]  
in Corollary \ref{cor:dimtangent}, there exists $X_i \in {\rm M}_n(k)$ such that $\delta_i(a) = ([X_i, a] \mod A_0) \in {\rm M}_n(k)/A_0$ 
for $a \in A_0$. Set $X = \sum_{i=1}^{r} x_i X_i \in {\rm M}_n(I)$. Then $\delta([a_i]) = ([X, a_i] \mod A_{R}\otimes_{R} I)$ 
in ${\rm M}_n(I)/(A_{R}\otimes_{R} I) = ({\rm M}_n(k)/A_0) \otimes_{k} I$ for $1 \le i \le d$.  On the other hand, $\delta([a_i]) = (c_i \mod A_{R}\otimes_{R} I)$ by the definition of $\delta$.  
There exists $d_i \in A_{R}\otimes_{R} I$ such that $[X, a_i]=c_i + d_i$ for $1 \le i \le d$. 
Put $P=I_n-X \in {\rm GL}_n(R)$. Note that $(P \mod I) = \overline{I}_n \in {\rm GL}_n(R/I)$. 
Using $P^{-1} = I_n + X$ and $c_i=b_i-a_i \in {\rm M}_n(I)$, we have 
\begin{eqnarray*} 
P^{-1}a_iP & = & (I_n+X)a_i(I_n-X)  =  a_i + [X, a_i] \\ 
 & = & a_i+c_i+d_i = a_i +b_i -a_i +d_i =b_i+d_i.
\end{eqnarray*}  
For $x \in I$, $a_i \otimes x \in A_{R}\otimes_{R} I \subset A_{R} \subseteq {\rm M}_n(R)$.   
Since $a_i\otimes x = a_i x = (b_i-c_i) x = b_i x \in B$ for $1 \le i \le d$, $A_{R}\otimes_{R} I \subseteq B$.  
Hence $P^{-1}a_iP = b_i +d_i \in B$, which implies that $P^{-1}A_RP \subseteq B$. 
By Lemma \ref{lemma:moldinclusionequal} below, $P^{-1}A_RP = B$. 

Conversely, suppose that there exists $P \in {\rm GL}_n(R)$ such that $P^{-1}A_R P =B$ and $(P \mod I) = \overline{I}_n \in {\rm GL}_n(R/I)$. We can write $P = I_n - X$, where $X \in {\rm M}_n(I)$. 
For a basis $a_1, a_2, \ldots, a_d$ of $A_R$ over $R$, set $b_i= P^{-1}a_iP$. 
Note that $b_i = (I_n+X)a_i(I_n-X) = a_i + [X, a_i]$ and that  
$b_1, b_2, \ldots, b_d$ is a basis of $B$ over $R$ such that $(a_i \mod I) = (b_i \mod I) \in {\rm M}_n(R/I)$. 
The derivation $\delta' :A_{R/I} \to {\rm M}_n(I)/(A_{R}\otimes I)$ can be written 
by $\delta'(a_i) = ([X, a_i] \mod A_{R}\otimes_{R} I)$ for $1 \le i \le d$.  
Then $\delta : A_0 \to ({\rm M}_n(k)/A_0)\otimes_{k} I$ is a $1$-coboundary in $C^{1}(A_0, ({\rm M}_n(k)/A_0)\otimes_{k} I)$.  Hence $[\delta] = 0$  in $H^{1}(A_0, ({\rm M}_n(k)/A_0)\otimes_{k} I)$. 
\qed 

\bigskip 

The following lemma has been used in Lemma \ref{lemma:eqdelta=0}. 

\begin{lemma}\label{lemma:moldinclusionequal} 
Let ${\mathcal A}$ and ${\mathcal B}$ be subbundles of rank $d$ of a locally free sheaf ${\mathcal E}$ of rank $m$ on a scheme $X$. If ${\mathcal A} \subseteq {\mathcal B}$, then ${\mathcal A} = {\mathcal B}$. 
In particular, if ${\mathcal A}, {\mathcal B} \subseteq {\rm M}_n({\mathcal O}_{X})$ are rank $d$ molds of degree $n$ on a scheme $X$ and if ${\mathcal A} \subseteq {\mathcal B}$, then ${\mathcal A} = {\mathcal B}$. 
\end{lemma} 

\proof 
By the assumption that ${\mathcal A}$ and ${\mathcal B}$ are subbuldes of ${\mathcal E}$, ${\mathcal E}/{\mathcal A}$ and ${\mathcal E}/{\mathcal B}$ are locally free sheaves of rank $m-d$ on $X$. If ${\mathcal A} \subseteq {\mathcal B}$, then we have the following commutative diagram with rows exact: 
\[
\begin{array}{ccccccccc}
 & & 0 & & & & & & \\ 
 & & \downarrow & & & & & & \\ 
0 & \to & {\mathcal A} & \to & {\mathcal E} & \to & {\mathcal E}/{\mathcal A} & \to & 0 \\  
   &  & \cap  &  & \| &  & \phantom{\psi} \downarrow \psi &  & \\   
0 & \to & {\mathcal B} & \to & {\mathcal E} & \to & {\mathcal E}/{\mathcal B} & \to & 0 \\
 & & & & & & \downarrow & & \\ 
 & & & & & & \phantom{.\quad} 0\quad . &  & \\ 
\end{array} 
\]
For $\psi : {\mathcal E}/{\mathcal A} \to {\mathcal E}/{\mathcal B}$, set ${\mathcal K} = {\rm Ker} \psi$. For proving that ${\mathcal A} = {\mathcal B}$, it suffices to show that ${\mathcal K} = 0$. Let $x \in X$. 
Taking the stalks at $x$, we have an exact sequence   
\[
0 \to {\mathcal K}_{x} \to {\mathcal E}_{x}/{\mathcal A}_{x} \stackrel{\psi_{x}}{\to} {\mathcal E}_{x}/{\mathcal B}_{x} \to 0. 
\] Since ${\mathcal E}_{x}/{\mathcal B}_{x}$ is free over the local ring ${\mathcal O}_{X, x}$, 
${\mathcal E}_{x}/{\mathcal A}_{x} \cong ({\mathcal E}_{x}/{\mathcal B}_{x}) \oplus {\mathcal K}_x$. In particular, ${\mathcal K}_x$ is a finitely generated module over ${\mathcal O}_{X, x}$. By taking the tensor products with the residue field $k(x) = {\mathcal O}_{X, x}/m_{x}$, we obtain 
$({\mathcal E}_{x}/{\mathcal A}_x)\otimes_{{\mathcal O}_{X, x}} k(x) \cong (({\mathcal E}_{x}/{\mathcal B}_x)\otimes_{{\mathcal O}_{X, x}} k(x)) \oplus ({\mathcal K}_{x} \otimes_{{\mathcal O}_{X, x}} k(x))$. 
Because both $({\mathcal E}_{x}/{\mathcal A}_x)\otimes_{{\mathcal O}_{X, x}} k(x)$ and $({\mathcal E}_{x}/{\mathcal B}_x)\otimes_{{\mathcal O}_{X, x}} k(x)$ have dimension $m-d$, ${\mathcal K}_{x} \otimes_{{\mathcal O}_{X, x}} k(x) = {\mathcal K}_x/m_x{\mathcal K}_x = 0$. By Nakayama's lemma, ${\mathcal K}_x =0$, which implies that ${\mathcal K}=0$. Hence ${\mathcal A}={\mathcal B}$. 
\qed 

\bigskip

Let $k$ be a field over the residue field $k(x)$ of a point $x \in S$. For a rank $d$ mold ${\mathcal A}$ on $S$, 
put $A_0 = {\mathcal A}\otimes_{{\mathcal O}_{S}} k \subseteq {\rm M}_n(k)$, which is a $d$-dimensional $k$-subalgebra of ${\rm M}_n(k)$. 
Let $R = k[\epsilon]/(\epsilon^2)$, and let $I = (\epsilon)$. 
We regard ${\rm Spec} R$ and ${\rm Spec} R/I$ as  $S$-schemes by ${\rm Spec} \: k = {\rm Spec} R/I \to {\rm Spec} R \to {\rm Spec} \: k \to {\rm Spec} \: k(x) \to S$ induced by %$k(x)$-algebra homomorphisms 
the canonical homomorphisms $k(x) \to k \to R \to R/I$.  Set $A = A_0 \otimes_{k} R \subseteq {\rm M}_n(k[\epsilon]/(\epsilon^2))$. Then $A$ is the pull-back of ${\mathcal A}$ on the $S$-scheme ${\rm Spec} R$. 
Take a basis $a_1, a_2, \ldots, a_d$ of $A_0$ over $k$. 
We can also regard $a_1, a_2, \ldots, a_d$ as an $R$-basis of the rank $d$ mold $A$ by 
$a_i = a_i + 0\epsilon \in {\rm M}_n(k) \subseteq {\rm M}_n(k)\oplus {\rm M}_n(k)\epsilon = {\rm M}_n(k[\epsilon]/(\epsilon^2)) = {\rm M}_n(R)$. 

Let us show that for any $\delta \in {\rm Der}_{k}(A_0, {\rm M}_n(k)/A_0)$, there exists a rank $d$ mold $B \subseteq {\rm M}_n(R)$ such that 
$B \otimes_{R} (R/I) = A\otimes_{R} (R/I) = A_0$ and $\delta : A_0 \to ({\rm M}_n(k)/A_0)\otimes_{k} I \cong {\rm M}_n(k)/A_0$ is the derivation associated to $B$. 
For a derivation $\delta \in {\rm Der}_{k}(A_0, {\rm M}_n(k)/A_0)$, choose a lift $\widetilde{\delta} : A_0 \to {\rm M}_n(k)$ of $\delta$ as a map. Set $b_i = a_i + \widetilde{\delta}(a_i)\epsilon \in {\rm M}_n(k[\epsilon]/(\epsilon^2))$ and $B = Rb_1 \oplus Rb_2 \oplus \cdots \oplus Rb_d \subseteq {\rm M}_n(k[\epsilon]/(\epsilon^2))$.  Note that $b_i \epsilon = a_i \epsilon \in B$ and that $A_0\otimes_{k} k\epsilon = A_0 \epsilon \subseteq B$ in ${\rm M}_n(k[\epsilon]/(\epsilon^2))$. We claim that the definition of $B$ does not depend on the choice of a lift $\widetilde{\delta}$ of $\delta$. Indeed, let us choose another lift $\widetilde{\delta'} : A_0 \to {\rm M}_n(k)$ of $\delta$. Set $b'_i = a_i + \widetilde{\delta'}(a_i)\epsilon$ and $B' = Rb'_1 \oplus Rb'_2 \oplus \cdots \oplus Rb'_d \subseteq {\rm M}_n(k[\epsilon]/(\epsilon^2))$. 
Since $b'_i - b_i = (\widetilde{\delta'}(a_i) - \widetilde{\delta}(a_i))\epsilon \in A_0 \epsilon \subseteq B$, 
we have $B' \subseteq B$. Similarly, we can verify that $B \subseteq B'$. Hence $B=B'$, which implies that the definition of $B$ does not depend on the choice of  lifts $\widetilde{\delta}$ of $\delta$. 
Note that $a+\widetilde{\delta}(a)\epsilon \in B$ for any $a \in A_0$ and that $B$ is generated by $\{ a+\widetilde{\delta}(a)\epsilon \mid a \in A_0 \}$ as an $R$-module.

Let us prove that $B$ is an $R$-subalgebra of ${\rm M}_n(k[\epsilon]/(\epsilon^2))$. Calculating $b_ib_j$, we have 
\begin{eqnarray*}
b_ib_j & = & (a_i+\widetilde{\delta}(a_i)\epsilon)(a_j+\widetilde{\delta}(a_j)\epsilon) \\  
 & = & a_ia_j + (a_i\widetilde{\delta}(a_j) + \widetilde{\delta}(a_i)a_j)\epsilon \\
 & = & a_ia_j + (\widetilde{\delta}(a_ia_j)+ c)\epsilon  
 \end{eqnarray*} 
 for some $c \in A_0$. 
Since $a+\widetilde{\delta}(a)\epsilon \in B$ for any $a \in A_0$, $a_ia_j + (\widetilde{\delta}(a_ia_j))\epsilon \in B$.    By using $c\epsilon \in A_0\epsilon \subseteq B$, we see that $b_ib_j \in B$. 
We easily see that $1 \in B$. 
Hence $B$ is an $R$-subalgebra of ${\rm M}_n(k[\epsilon]/(\epsilon^2))$. 
We also see that $B$ is a rank $d$ mold on $R$ such that $B\otimes_{R} k = A\otimes_{R} k= A_0$. 

We denote by $\psi : {\rm Spec} R \to {\rm Mold}_{n, d}\otimes_{{\Bbb Z}} S$ 
the morphism induced by the rank $d$ mold $B$. 
We also denote by $\overline{g} : {\rm Spec} R/I = {\rm Spec}\: k \to {\rm PGL}_{n, S}$ the 
morphism given by the identity $[I_n]$. 
Then we obtain commutative diagram (\ref{diagram:phiasmooth}). 
Diagram (\ref{diagram:phiasmooth}) induces a derivation $\delta' : A_{R/I}=A\otimes_{R} (R/I) = A_0 \to ({\rm M}_n(k)/A_0)\otimes_{k} I$, that is, $\delta'(a_i) = b_i - a_i = (\widetilde{\delta}(a_i) \mod A_0)\otimes \epsilon = \delta(a_i)\otimes \epsilon$ for $1 \le i \le d$.  
Hence $\delta : A_0 \to ({\rm M}_n(k)/A_0)\otimes_{k} I \cong  {\rm M}_n(k)/A_0$ is the derivation associated to $B$. 
Therefore, we have the following lemma. 

\begin{lemma}\label{lemma:existencemoldB}
Let $k$ be a field over the residue field $k(x)$ of a point $x \in S$. 
Let $R=k[\epsilon]/(\epsilon^2)$, and let $I = (\epsilon)$. 
Put $A_0 = {\mathcal A}\otimes_{{\mathcal O}_{S}} k \subseteq {\rm M}_n(k)$. 
Set $A_{R} = A_0\otimes_{k}R$ and $A_{R/I} = A_0\otimes_{k}(R/I) = A_0$. 
For any $\delta \in {\rm Der}_{k}(A_0, {\rm M}_n(k)/A_0)$, there exists a rank $d$ mold $B \subseteq {\rm M}_n(R)$ such that $\delta : A_0 \to ({\rm M}_n(k)/A_0)\otimes_{k} I \cong  {\rm M}_n(k)/A_0$ is the derivation associated to $B$ with $B\otimes_{R} (R/I) = A_{R/I}$. 
\end{lemma} 

Now we have: 

\begin{theorem}\label{th:orbitsmooth} 
Let ${\mathcal A}$ be a rank $d$ mold of degree $n$ on a locally noetherian scheme $S$. 
%$\tau_{{\mathcal A}} : S \to {\rm Mold}_{n, d}\otimes_{{\Bbb Z}} S$ the morphism induced by ${\mathcal A}$. 
%For a rank $d$ mold ${\mathcal A}$ of degree $n$ on $S$, denote by 
%$\tau_{{\mathcal A}} : S \to {\rm Mold}_{n, d}\otimes_{{\Bbb Z}} S$ the morphism induced by ${\mathcal A}$. 
Set ${\mathcal A}(x)= {\mathcal A}\otimes_{{\mathcal O}_S} k(x) \subseteq {\rm M}_n(k(x))$, where $k(x)$ is the residue field of a point $x \in S$. 
Put ${\rm PGL}_{n, S} = {\rm PGL}_n\otimes_{{\Bbb Z}} S$. 
Let us define the $S$-morphism $\phi_{{\mathcal A}} : {\rm PGL}_{n, S} \to {\rm Mold}_{n, d}\otimes_{{\Bbb Z}} S$ by 
$P \mapsto P^{-1}{\mathcal A}P$. Then $\phi_{{\mathcal A}}$ is smooth if and only if $H^{1}({\mathcal A}(x), {\rm M}_n(k(x))/{\mathcal A}(x))=0$ for each $x \in S$. 
\end{theorem} 

\proof
Assume that $H^{1}({\mathcal A}(x), {\rm M}_n(k(x))/{\mathcal A}(x))=0$ for each $x \in S$. Let $I$ be an ideal of an Artin local ring $(R, m, k)$ over $S$ with $mI=0$. 
Suppose that $B \subseteq {\rm M}_n(R)$ is a rank $d$ mold with $B \otimes_{R} (R/I) = A_{R/I}$. By Lemma \ref{lemma:pmodiIncase}, it suffices to prove that there exists $P \in {\rm GL}_n(R)$ such that $P^{-1}A_{R}P =B$ and 
$(P \mod I) = \overline{I}_n \in {\rm GL}_n(R/I)$. 
Let $\delta : A_0 \to ({\rm M}_n(k)/A_0)\otimes_{k} I$ be the derivation associated to $B$. 
Denote by $x \in S$ the image of $m$ by the canonical morphism ${\rm Spec} R \to S$. 
Then $k$ is a field over $k(x)$ and $A_0 = {\mathcal A}(x)\otimes_{k(x)} k \subseteq {\rm M}_n(k)$.   
By the assumption that $H^{1}({\mathcal A}(x), {\rm M}_n(k(x))/{\mathcal A}(x))=0$ and Proposition \ref{prop:hochextfield}, 
$H^{1}(A_0, {\rm M}_n(k)/A_0) = H^{1}({\mathcal A}(x)\otimes_{k(x)} k, ({\rm M}_n(k(x))/{\mathcal A}(x))\otimes_{k(x)} k) 
= H^{1}({\mathcal A}(x), {\rm M}_n(k(x))/{\mathcal A}(x))\otimes_{k(x)} k = 0$. 
Hence $H^{1}(A_0, ({\rm M}_n(k)/A_0)\otimes_{k} I) = H^{1}(A_0, {\rm M}_n(k)/A_0)\otimes_{k} I = 0$ and the cohomology class $[\delta]$ is $0$. 
By Lemma \ref{lemma:eqdelta=0}, there exists $P \in {\rm GL}_n(R)$ such that $P^{-1}A_{R}P =B$ and 
$(P \mod I) = \overline{I}_n \in {\rm GL}_n(R/I)$. 

Conversely, assume that $\phi_{{\mathcal A}}$ is smooth. Let $x \in S$. 
By Corollary \ref{cor:dimtangent}, there exists a surjection ${\rm Der}_{k(x)}({\mathcal A}(x), {\rm M}_n(k(x))/{\mathcal A}(x)) \to H^{1}({\mathcal A}(x), {\rm M}_n(k(x))/{\mathcal A}(x)) \to 0$.  It suffices to show that $[\delta]=0$ in $H^{1}({\mathcal A}(x), {\rm M}_n(k(x))/{\mathcal A}(x))$ for any 
$\delta \in {\rm Der}_{k(x)}({\mathcal A}(x), {\rm M}_n(k(x))/{\mathcal A}(x))$. 
Let $k = k(x)$ and  $A_0 = {\mathcal A}\otimes_{{\mathcal O}_S} k = {\mathcal A}(x) \subseteq {\rm M}_n(k)$.  
%for $\delta \in {\rm Der}_{k}({\mathcal A}(x), {\rm M}_n(k(x))/{\mathcal A}(x))$.   
By Lemma \ref{lemma:existencemoldB}, there exists a rank $d$ mold $B \subseteq {\rm M}_n(R)$ such that 
$\delta : A_0 \to ({\rm M}_n(k)/A_0)\otimes_{k} I \cong {\rm M}_n(k)/A_0$ is the derivation associated to $B$ with $B\otimes_{R} (R/I) =A_{R/I}$, where $R = k[\epsilon]/(\epsilon^2)$, $I=(\epsilon)$, $A_{R} = A_0\otimes_{k} R$, and $A_{R/I} = A_0\otimes_{k} (R/I) = A_0$.  
Using Lemma \ref{lemma:pmodiIncase}, we have $P \in {\rm GL}_n(R)$ such that $P^{-1}A_{R}P =B$ and $(P \mod I) = \overline{I}_n \in {\rm GL}_n(R/I)$, because $\phi_{{\mathcal A}}$ is smooth. 
Hence, Lemma \ref{lemma:eqdelta=0} implies that $[\delta]=0$ in $H^{1}(A_0, ({\rm M}_n(k)/A_0)\otimes_{k}I) \cong  
H^{1}({\mathcal A}(x), {\rm M}_n(k(x))/{\mathcal A}(x))$. Thereby, $ H^{1}({\mathcal A}(x), {\rm M}_n(k(x))/{\mathcal A}(x)) = 0$ for each $x \in S$. 
\qed 

\begin{corollary}\label{cor:openorbit} 
In the situation of Theorem~\ref{th:orbitsmooth}, assume that $H^{1}({\mathcal A}(x), {\rm M}_n(k(x))/{\mathcal A}(x))=0$ for each $x \in S$. 
Then ${\rm Im} \phi_{{\mathcal A}}$ is open in ${\rm Mold}_{n, d}\otimes_{{\Bbb Z}} S$. 
\end{corollary} 

\proof 
By Theorem \ref{th:orbitsmooth}, the assumption implies that $\phi_{{\mathcal A}} : {\rm PGL}_{n, S} \to {\rm Mold}_{n, d}\otimes_{{\Bbb Z}} S$ is smooth. In particular, 
$\phi_{{\mathcal A}}$ is flat morphism locally of finite presentation.  
Hence $\phi_{{\mathcal A}}$ is open, which completes the proof. 
\qed 

\section{How to calculate Hochschild cohomology groups} 
In this section, we introduce how to calculate Hochschild cohomology groups.  
By using Cibils's result (Proposition~\ref{prop:cibils}), %in \cite{Cibils}, 
we can calculate Hochschild cohomology for several cases. As a result, we see that if $\Lambda$ is the incidence algebra of an ordered quiver $Q$ with 
$n = |Q_{0}|$, then $H^{i}(\Lambda, {\rm M}_n(R)/\Lambda) = 0$ for $i \ge 0$ (Theorem~\ref{th:incidence}).   
We also explain several techniques and perform several calculations. 

\bigskip 

Let $Q$ be a finite quiver. 
Denote by $Q_0$ the set of vertices of $Q$. 
Let $RQ$ be the path algebra over a commutative ring $R$. 
We define the {\it arrow ideal} $F$ as the two-sided ideal of $RQ$ generated by  
the paths of positive length of $Q$. 
A two-sided ideal of $I$ of $RQ$ is called {\it admissible} if $F^n \subseteq I \subseteq F$ for a positive 
integer $n$ and $F/I$ is an $R$-free module which has an $R$-basis consisting of 
oriented paths.  
For an admissible ideal $I$, set $\Lambda = RQ/I$ and $r=F/I$.  Denote by $E$ 
the $R$-subalgebra of $\Lambda$ generated by $Q_0$.  

\begin{proposition}[\cite{Cibils}, Proposition~1.2]\label{prop:cibils}
Let $M$ be a $\Lambda$-bimodule. 
The Hochschild cohomology $R$-modules $H^{i}(\Lambda, M)$ are 
the cohomology groups of the complex %of $E$-bimodules (i.e. modules over 
%$E^{e} = E\otimes_{R}E^{op}$) 
\begin{multline*} 
0 \to M^{E} \stackrel{\delta^0}{\to} {\rm Hom}_{E^{e}}(r, M) 
\stackrel{\delta^1}{\to} {\rm Hom}_{E^{e}}(r\otimes_{E}r, M) 
\stackrel{\delta^2}{\to} \cdots \\ 
\cdots \stackrel{\delta^{i-1}}{\to}  
{\rm Hom}_{E^{e}}(r^{\otimes i}, M) \stackrel{\delta^{i}}{\to} 
{\rm Hom}_{E^{e}}(r^{\otimes i+1}, M) \stackrel{\delta^{i+1}}{\to}  \cdots ,  
\end{multline*}
where the tensor products are over $E$,  
\begin{eqnarray*} 
M^{E} & = & \{ m \in M \mid sm = ms \mbox{ for each } s \in Q_{0} \} = \oplus_{s \in Q_0} sMs, \\ 
\delta^{0}(m)(x) & = & xm-mx \mbox{ for } m \in M^{E} \mbox{ and } x \in r, 
\end{eqnarray*} 
and 
\begin{multline*} 
\delta^{i}(f)(x_1\otimes \cdots \otimes x_{i+1})  =   x_1f(x_2 \otimes \cdots \otimes x_{i+1})   
+ \sum_{j=1}^{i} (-1)^{j} f(x_1\otimes \cdots \otimes x_{j}x_{j+1} \otimes \cdots \otimes x_{i+1})  \\ 
+(-1)^{i+1}f(x_1\otimes \cdots \otimes x_{i})x_{i+1} 
\end{multline*} 
for $f \in {\rm Hom}_{E^{e}}(r^{\otimes i}, M)$ and $x_1\otimes \cdots \otimes x_{i+1} \in r^{\otimes i+1}$. 
\end{proposition}

\begin{remark}\rm
Set $r^{\otimes 0} = E$. Then ${\rm Hom}_{E^{e}}(r^{\otimes 0}, M) = M^{E}$.   
Hence the complex above can be written by $\{ {\rm Hom}_{E^{e}}(r^{\otimes n}, E),\delta^{n}\}$. 
\end{remark}

\bigskip 

Denote by $Q_1$ the set of arrows of a finite quiver $Q$. 
For each oriented path $\alpha$ of $Q$, 
we denote by $h(\alpha)$ and $t(\alpha)$ the head and the tail of $\alpha$, 
respectively. 

\begin{definition}\rm
Let $Q$ be a finite quiver without oriented cycles. 
We say that $Q$ is {\it ordered}   
if there exists no oriented path other than $\alpha$ joining $t(\alpha)$ to $h(\alpha)$ for 
each arrow $\alpha \in Q_1$.  
\end{definition}

\begin{definition}\rm 
Let $Q$ be an ordered quiver. 
Let $I$ be the two-sided ideal of $RQ$ generated by 
$$\{ \gamma-\delta \in RQ \mid 
\mbox{ $\gamma$ and $\delta$ are oriented paths with } 
h(\gamma)=h(\delta) \mbox{ and }   t(\gamma)=t(\delta) \}.$$ 
We call $\Lambda = RQ/I$ an {\it incidence $R$-algebra}. 
Note that $I$ is an admissible ideal. 
\end{definition}

\bigskip 

For an ordered quiver $Q$, set $n = | Q_0 |$. 
For $a, b \in Q_0$, we define $a \ge b$ if $a = b$ or there exists an 
oriented path $\alpha$ such that $t(\alpha)=a$ and $h(\alpha)=b$. 
Then $(Q_0, \ge)$ is a partially ordered set (i.e. poset).  
Let $\Lambda = RQ/I$ be the incidence algebra associated to $Q$. 
For $a \ge b$, let $e_{ba}$ be the equivalence class of oriented paths from $a$ to $b$ in 
$\Lambda$. We can write $\Lambda = \oplus_{a \ge b} Re_{ba}$. 
Fix a numbering on $Q_0$.  
By regarding $e_{ba}$ as $E_{ba}$, $\Lambda$ can be considered as an $R$-subalgebra of 
${\rm M}_n(R) = \oplus_{a, b\in Q_0} RE_{ba}$, where $E_{ba}$ is the matrix unit.  
We can write $E = \oplus_{a \in Q_0} Re_{aa}$ and 
$E^{e} = E\otimes_{R} E^{op} = \oplus_{a, b \in Q_0} Re_{aa}\otimes e_{bb}$. 
We also have $r = F/I = \oplus_{a>b} Re_{ba}$. (In the sequel, we denote 
$E_{ba} \in {\rm M}_n(R)$ by $e_{ba}$ for simplicity.)  

\begin{lemma}\label{lemma:hom=0}
For $i \ge 0$, ${\rm Hom}_{E^{e}}(r^{\otimes i}, {\rm M}_n(R)/\Lambda) = 0$. 
\end{lemma}

\proof
As $E$-bimodules, $r^{\otimes i}$ is isomorphic to $\oplus_{s_0 > s_1 > \cdots > s_i} Re_{s_is_0}$. 
On the other hand, ${\rm M}_n(R)/\Lambda \cong \oplus_{a \not\ge b} Re_{ba}$. 
Hence we have ${\rm Hom}_{E^{e}}(r^{\otimes i}, {\rm M}_n(R)/\Lambda) \cong 
\oplus_{s_0 > s_1 > \cdots > s_i,\;  a \not\ge b}  \;  
{\rm Hom}_{E^{e}}(Re_{s_is_0}, Re_{ba})$.  
Since ${\rm Hom}_{E^{e}}(Re_{s_is_0}, Re_{ba}) \cong e_{s_is_i} (Re_{ba}) e_{s_0s_0} = 0$, 
${\rm Hom}_{E^{e}}(r^{\otimes i}, {\rm M}_n(R)/\Lambda) = 0$. 
\qed 

\bigskip 

Summarizing the discussion above, we have the following theorem. 

\begin{theorem}\label{th:incidence} 
Let $Q$ be an ordered quiver with $n= | Q_0 |$. Let $\Lambda$ be the incidence algebra 
associated to $Q$. Then  
$H^{i}(\Lambda, {\rm M}_n(R)/\Lambda) = 0$ for $i \ge 0$. 
\end{theorem}

\proof 
By Proposition \ref{prop:cibils} and Lemma \ref{lemma:hom=0},  we can prove the statement. 
\qed

\bigskip 

We show several examples of Hochschild cohomology groups $H^{i}(A, {\rm M}_n(R)/A)$ 
for $R$-subalgebras $A$ of ${\rm M}_n(R)$. We also refer to the moduli ${\rm Mold}_{n, d}$ of 
molds. 
Let ${\mathcal A}$ be the universal mold on ${\rm Mold}_{n, d}$. 
For $x \in {\rm Mold}_{n, d}$, we set 
${\mathcal A}(x) = {\mathcal A}\otimes_{{\mathcal O}_{{\rm Mold}_{n, d}}} k(x)$, 
where $k(x)$ is the residue field of $x$.

\begin{example}\label{ex:borel}\rm 
Let $R$ be a commutative ring, and let us consider the following quiver $Q$:  
\[
1 \longleftarrow 2 \longleftarrow 3 \longleftarrow \cdots \longleftarrow n. 
\]
%\[
% \xymatrix{ 
%   1 & \ar[l]  2 & \ar[l]  3  & \ar[l] \cdots  & \ar[l] n.  
% }
%\]
Let $\Lambda = RQ/I$ be the incidence algebra associated to $Q$ over a commutative ring $R$. 
Then $\Lambda = \oplus_{1 \le i \le j  \le n} Re_{ij}$. 
We can regard $\Lambda$ as the upper triangular matrix ring ${\mathcal B}_n(R) 
= \{  (a_{ij}) \in {\rm M}_n(R) \mid  a_{ij}=0 \mbox{ for } i>j  \}$.  
By Theorem \ref{th:incidence}, $H^{i}({\mathcal B}_n(R), {\rm M}_n(R)/{\mathcal B}_n(R)) = 0$ 
for $i \ge 0$. 

This result is compatible with the fact that the connected component containing ${\mathcal B}_n$ in 
${\rm Mold}_{n, d}$ is isomorphic to ${\rm GL}_n/{\rm B}_n$, where 
$d = n(n+1)/2$ and ${\rm B}_n = \{  (a_{ij}) \in {\rm GL}_n \mid  a_{ij}=0 \mbox{ for } i>j  \}$ 
(For details, see \cite[Theorem~1.1]{Nakamoto2}). 
Indeed, 
the image of the morphism $\phi_{{\mathcal B}_n} : {\rm PGL}_{n} \to {\rm Mold}_{n,d}$ 
associated to the mold ${\mathcal B}_{n}({\mathbb Z})$ on ${\rm Spec}\: {\mathbb Z}$ 
is open by Corollary~\ref{cor:openorbit}, since 
$H^{1}({\mathcal B}_{n}(R), 
{\rm M}_n(R)/{\mathcal B}_{n}(R)) = 0$ for any commutative ring $R$. 
It is easy to see that 
${\rm Im} \phi_{{\mathcal B}_{n}} = {\rm GL}_n/{\rm B}_{n}$ is irreducible, open and closed, and hence that ${\rm GL}_n/{\rm B}_{n}$ is an irreducible component and a connected component.    
We also see that 
$H^{2}({\mathcal B}_{n}(R), 
{\rm M}_n(R)/{\mathcal B}_{n}(R)) = 0$ is compatible with the fact that 
${\rm GL}_n/{\rm B}_{n}$ 
is smooth over ${\Bbb Z}$ (see Theorem~\ref{th:smooth}). 
If ${\mathcal A}(x) = 
{\mathcal B}_n(k(x))$ for a point $x \in {\rm Mold}_{n, d}$, then  
$\dim_{k(x)} T_{{\rm Mold}_{n, d}/{\Bbb Z}, x} = \dim_{k(x)}  
H^{1}({\mathcal B}_n(k(x)), {\rm M}_n(k(x))/{\mathcal B}_n(k(x))) + n^2 
- \dim_{k(x)} N({\mathcal B}_n(k(x))) 
= n^2 - \dim_{k(x)} {\mathcal B}_n(k(x)) = \dim {\rm GL}_n(k(x))/{\rm B}_n(k(x)) = n(n-1)/2$ 
by Corollary \ref{cor:dimtangent}, 
since $N({\mathcal B}_n(k(x))) = {\mathcal B}_n(k(x))$.  
For more general result, see Example \ref{ex:parabolic}. 
\end{example}

\begin{example}\label{ex:scalar}\rm  
Let $R$ be a commutative ring, and  
let $A = RI_n \subset {\rm M}_n(R)$. 
The bar complex $C^{i}(RI_n, {\rm M}_n(R)/RI_n)$ is isomorphic to 
$0 \to {\rm M}_n(R)/RI_n \stackrel{d^0}{\to} {\rm M}_n(R)/RI_n \stackrel{d^1}{\to} {\rm M}_n(R)/RI_n \stackrel{d^2}{\to} \cdots$, where $d^i = 0$ if $i$ is even and 
$d^i = id_{{\rm M}_n(R)/RI_n}$ if $i$ is odd.  Hence we have 
\[
H^{i}(RI_n, {\rm M}_n(R)/RI_n) \cong \left\{  
\begin{array}{cl} 
{\rm M}_n(R)/RI_n & (i=0) \\
0 & (i>0). 
\end{array}
\right. 
\] 
The moduli ${\rm Mold}_{n, 1}$ is smooth over ${\Bbb Z}$, since 
it is isomorphic to ${\rm Spec} \: {\Bbb Z}$. 
This is compatible with the fact that 
$H^{2}(RI_n, {\rm M}_n(R)/RI_n) = 0$ (see Theorem~\ref{th:smooth}).    
Note that ${\mathcal A}(x) = k(x)I_n$ for each point $x \in {\rm Mold}_{n, 1}$. 
Then  $\dim_{k(x)} T_{{\rm Mold}_{n, 1}/{\Bbb Z}, x} = \dim_{k(x)}  
H^{1}(k(x)I_n, {\rm M}_n(k(x))/k(x)I_n) + n^2 
- \dim_{k(x)} N(k(x)I_n) 
= n^2 - \dim_{k(x)} {\rm M}_n(k(x)) = 0$ 
by Corollary \ref{cor:dimtangent}, 
since $N(k(x)I_n) = {\rm M}_n(k(x))$.  
\end{example}

\begin{example}\label{ex:fullmatrix}\rm  
Let $R$ be a commutative ring, and  
let $A = {\rm M}_n(R)$. 
Since ${\rm M}_n(R)/{\rm M}_n(R) = 0$, 
$H^{i}({\rm M}_n(R), {\rm M}_n(R)/{\rm M}_n(R))=0$ for $i \ge 0$. 
The moduli ${\rm Mold}_{n, n^2}$ is isomorphic to ${\rm Spec} \: {\Bbb Z}$ 
(see \cite[Proposition~1.1]{Nakamoto2}), 
and hence it is smooth over ${\Bbb Z}$.  
This is compatible with the fact that 
$H^{2}({\rm M}_n(R), {\rm M}_n(R)/{\rm M}_n(R)) = 0$ (see Theorem~\ref{th:smooth}).     
We see that $\dim_{k(x)} T_{{\rm Mold}_{n, n^2}/{\Bbb Z}, x} = \dim_{k(x)}  
H^{1}({\rm M}_n(k(x)), {\rm M}_n(k(x))/{\rm M}_n(k(x))) + n^2 
- \dim_{k(x)} N({\rm M}_n(k(x))) = 0$ for $x \in {\rm Mold}_{n, n^2}$ 
by Corollary \ref{cor:dimtangent}.   
\end{example} 

\begin{example}\label{ex:dn}\rm  
Let $R$ be a commutative ring, and  
let $A = {\rm D}_n(R) = \{ (a_{ij}) \in {\rm M}_n(R) \mid a_{ij} = 0  \mbox { for } i\neq j \} 
\subset {\rm M}_n(R)$. 
In other words, ${\rm D}_n(R)$ is the $R$-subalgebra of diagonal matrices in 
${\rm M}_n(R)$. 
Let $Q$ be a quiver with $Q_0 = \{ 1, 2, \ldots, n \}$ and $Q_1 = \emptyset$. 
Then ${\rm D}_n(R) = RQ = \oplus_{i=1} Re_{ii} \subset {\rm M}_n(R) = 
\oplus_{i, j=1}^{n} Re_{ij}$ and ${\rm M}_n(R)/{\rm D}_n(R) = \oplus_{i\neq j} Re_{ij}$. 
The arrow ideal $F$ of $RQ$ is $0$. Set $I = F = 0$. 
Then $\Lambda = RQ/I = RQ = {\rm D}_n(R)$,  $r= F/I = 0$, and $E= {\rm D}_n(R)$. 
The complex  
in Proposition \ref{prop:cibils} for $M = {\rm M}_n(R)/{\rm D}_n(R)$ is the zero cochain complex, 
since $r = 0$ and $M^{E} = ({\rm M}_n(R)/{\rm D}_n(R))^{{\rm D}_n(R)} = 0$. 
Hence $H^{i}({\rm D}_n(R), {\rm M}_n(R)/{\rm D}_n(R))=0$ for $i \ge 0$. 
This result also follows from that ${\rm D}_n(R)$ is a separable $R$-algebra.  
\end{example} 

\begin{definition}\label{def:parabolic}\rm
Let $n_1, n_2, \ldots, n_s$ be positive integers with $\sum_{i=1}^{s} n_i = n$. 
We define the $R$-subalgebra ${\mathcal P}_{n_1, n_2, \ldots, n_s}(R)$ of ${\rm M}_n(R)$ 
over a commutative ring $R$ by  
\[
{\mathcal P}_{n_1, n_2, \ldots, n_s}(R) = 
\{
(a_{ij}) \in {\rm M}_n(R) \mid a_{ij}=0 \mbox{ if } \sum_{k=1}^{t} n_k < i \le \sum_{k=1}^{t+1} n_k 
\mbox{ and }  j \le \sum_{k=1}^{t} n_k  
\}.  
\] 
To simplify notation, we write ${\mathcal P}_{\bf n}(R)$ instead of ${\mathcal P}_{n_1, n_2, \ldots, n_s}(R)$ 
for ${\bf n} = (n_1, n_2, \ldots, n_s)$. 
Set  
\[
E = \left\{ 
\begin{array}{c|c}  
\left(
\begin{array}{ccccc} 
X_1 & 0 & 0 & \cdots & 0 \\
0 & X_2 & 0 & \cdots & 0 \\
0 & 0 & X_3 & \cdots & 0 \\
\vdots & \vdots & \ddots & \ddots & \vdots \\
0 & 0 & 0 & \cdots & X_s \\
\end{array} 
\right) \in {\mathcal P}_{\bf n}(R) & X_i \in {\rm M}_{n_i}(R) \mbox{ for } 1 \le i \le s 
\end{array} 
\right\}
\] 
and 
\[
r = \left\{ 
\begin{array}{c|c}  
\left(
\begin{array}{ccccc} 
0 & X_{12} & X_{13} & \cdots & X_{1s}  \\
0 & 0 & X_{23} & \cdots & X_{2s} \\
0 & 0 & 0 & \cdots & X_{3s} \\
\vdots & \vdots & \ddots & \ddots & \vdots \\
0 & 0 & 0 & \cdots & 0 \\
\end{array} 
\right) \in {\mathcal P}_{\bf n}(R) & 
\begin{array}{c} 
X_{ij} \in {\rm M}_{n_i, n_j}(R) \\ 
\mbox{ for } 1 \le i < j  \le s \\ 
\end{array} 
\end{array} 
\right\}, 
\] 
where ${\rm M}_{i, j}(R)$ is the set of $(i\times j)$-matrices over $R$. 
Note that $E$ is an $R$-subalgebra of ${\mathcal P}_{\bf n}(R)$ and that ${\mathcal P}_{\bf n}(R) = E\oplus r$ as $E$-bimodules.  
We also set 
\[
r_{ii} = \left\{ 
\begin{array}{c|c}  
\left(
\begin{array}{ccccc} 
X_1 & 0 & 0 & \cdots & 0 \\
0 & X_2 & 0 & \cdots & 0 \\
0 & 0 & X_3 & \cdots & 0 \\
\vdots & \vdots & \ddots & \ddots & \vdots \\
0 & 0 & 0 & \cdots & X_s \\
\end{array} 
\right) \in E & 
\begin{array}{c}
X_{i} \in {\rm M}_{n_i}(R), \mbox{ but } X_{j} \mbox{ equals  } 0 \\ 
\mbox{ in } {\rm M}_{n_j}(R)  \mbox{ for }  j \neq i   
\end{array} 
\end{array} 
\right\}   
\] 
for $1 \le i \le s$ and  
\[
r_{ij} = \left\{ 
\begin{array}{c|c}  
\left(
\begin{array}{ccccc} 
0 & X_{12} & X_{13} & \cdots & X_{1s}  \\
0 & 0 & X_{23} & \cdots & X_{2s} \\
0 & 0 & 0 & \cdots & X_{3s} \\
\vdots & \vdots & \ddots & \ddots & \vdots \\
0 & 0 & 0 & \cdots & 0 \\
\end{array} 
\right) \in r & 
\begin{array}{c}
X_{ij} \in {\rm M}_{n_i, n_j}(R), \mbox{ but } X_{kl} \mbox{ equals  } 0 \\ 
\mbox{ in } {\rm M}_{n_k, n_l}(R)  \mbox{ for }  (k, l) \neq (i, j)  
\end{array} 
\end{array} 
\right\}     
\] 
for $1 \le i < j \le s$. 
We easily see that $r_{ij}$ is an $E$-bimodule and that 
${\mathcal P}_{\bf n}(R) = \oplus_{1 \le i \le j \le s} r_{ij}$ and 
$r = \oplus_{1 \le i < j \le s} r_{ij}$ as $E$-bimodules. 
\end{definition}

\bigskip 

To calculate $H^{i}({\mathcal P}_{\bf n}(R), 
{\rm M}_n(R)/{\mathcal P}_{\bf n}(R))$, we need to make several preparations.  

\begin{proposition}\label{prop:rijeprojective} 
For $1 \le i \le j  \le s$, $r_{ij}$ is a projective $E$-bimodule. In particular, ${\mathcal P}_{\bf n}(R)$ and $r$ are projective $E$-bimodules. 
\end{proposition} 

\proof 
The $E$-bimodule $E\otimes_{R} E = \oplus_{1 \le i, j  \le s} r_{ii}\otimes_{R} r_{jj}$ is isomorphic to 
$E\otimes E^{op}$ as $E^{e} = E\otimes_{R}E^{op}$-modules. Hence $r_{ii}\otimes_{R}r_{jj}$ is 
a projective $E$-bimodule for $1 \le i, j \le s$. 
For $1 \le i \le j \le s$, we can easily check that $r_{ii}\otimes_{R} r_{jj} \cong r_{ij}^{\oplus n_{i}n_{j}}$ as $E$-bimodules.  Therefore, $r_{ij}$ is a projective $E$-bimodule. The last statement follows from that 
${\mathcal P}_{\bf n}(R) = \oplus_{1 \le i \le j \le s} r_{ij}$ and 
$r = \oplus_{1 \le i < j \le s} r_{ij}$. 
\qed 

\begin{proposition}\label{prop:parabolicprojresol} 
For a commutative ring $R$,  
put $P = {\mathcal P}_{\bf n}(R)$. Let $E$ and $r$ be as in Definition \ref{def:parabolic}. 
The following complex gives a projective resolution of $P$ in the category of $P$-bimodules: 
\[
\cdots \to P\otimes_{E}r^{\otimes i} \otimes_{E} P \stackrel{d_i}{\to} P\otimes_{E}r^{\otimes i-1} \otimes_{E} P  
\to \cdots \to P\otimes_{E}r \otimes_{E} P \stackrel{d_1}{\to} P\otimes_{E}P \stackrel{d_0}{\to} P \to 0, 
\] 
where $r^{\otimes i} = r\otimes_{E} \cdots \otimes_{E} r$ ($i$ times),  
$d_i : P\otimes_{E}r^{\otimes i} \otimes_{E} P \stackrel{d_i}{\to} P\otimes_{E}r^{\otimes i-1}\otimes_{E} P$ is the $P$-bimodule homomorphism determined by 
$d_i(1 \otimes x_1\otimes \cdots \otimes x_i \otimes 1) = x_1\otimes x_2 \otimes \cdots \otimes x_i \otimes 1 + \sum_{j=1}^{i-1} (-1)^j 1\otimes x_1 \otimes \cdots \otimes x_{j}x_{j+1} \otimes \cdots \otimes x_{i} \otimes 1 + (-1)^i 1\otimes x_1 \otimes x_2 \otimes \cdots \otimes x_i$, and $d_0 : P\otimes_{E} P \to P$ is defined by $d_0(a\otimes b) = ab$. 
\end{proposition}

\proof 
As $E$-bimodules, $r_{ij}\otimes_{E} r_{jl} \cong r_{il}$ and 
$r_{ij} \otimes_{E} r_{kl} \cong 0$ for $j \neq k$.  We see that 
$r^{\otimes i} \cong \oplus_{1 \le j_1 < j_2 < \cdots < j_{i+1} \le s} \; r_{j_1j_2} \otimes_{E} \cdots \otimes_{E} r_{j_{i}j_{i+1}} \cong \oplus_{1 \le j_1 < j_2 < \cdots < j_{i+1} \le s} \; r_{j_1j_{i+1}}$.  (For convenience, set $r^{\otimes 0} = E$.)  
By Proposition \ref{prop:rijeprojective}, $r^{\otimes i}$ is a projective $E$-bimodule for $i \ge 0$.  
There exists an $E^e = E\otimes_{R}E^{op}$-module $M$ such that $r^{\otimes i}\oplus M \cong (E\otimes_{R}E^{op})^{\oplus q}$ as $E^{e}$-modules for some $q \in \mathbb{Z}$. 
Then we have 
\begin{eqnarray*}
(P\otimes_{E} r^{\otimes i} \otimes_{E} P) \oplus (P\otimes_{E} M \otimes_{E} P) &\cong & P\otimes_{E} (r^{\otimes i}\oplus M) \otimes_{E} P  \\ 
& \cong &  (P\otimes_{E} (E\otimes_{R}E^{op}) \otimes_{E} P)^{\oplus q} \\ 
& \cong &  (P\otimes_{R}P)^{\oplus q}.  
\end{eqnarray*} 
Since $P\otimes_{R}P$ is a projective $P$-bimodule,  $P\otimes_{E} r^{\otimes i} \otimes_{E} P$ is also a projective $P$-bimodule for $i \ge 0$.   

Let us show that the complex is exact. For $\lambda \in P = E \oplus r$, we write 
$\lambda = \lambda_E + \lambda_r$, where $\lambda_E \in E$ and $\lambda_r \in r$. 
For $i \ge 1$, we define the $R$-homomorphism $t_i : P \otimes_{E} r^{\otimes i-1} \otimes_{E} P \to P \otimes_{E} r^{\otimes i} \otimes_{E} P$ by $t_{i}(\lambda\otimes x_1\otimes \cdots \otimes x_{i-1} \otimes \mu) = 1\otimes \lambda_r \otimes x_1 \otimes \cdots x_{i-1} \otimes \mu$. 
We also define $t_0 : P \to P\otimes_{E} P$ by $\lambda \mapsto 1 \otimes \lambda$. 
By $d_0t_0(\lambda) = d_0(1\otimes \lambda)=\lambda$, we have $d_0t_0 = id$.  
Next, let us check $t_0d_0+d_1t_1 = id$. By $t_0d_0(\lambda\otimes \mu) = t_0(\lambda\mu) = 1\otimes \lambda\mu$ and $d_1t_1(\lambda\otimes \mu) = d_1(1\otimes \lambda_r \otimes \mu) = \lambda_r\otimes \mu - 1\otimes \lambda_r \mu$,  
\begin{eqnarray*} 
(t_0d_0+d_1t_1)(\lambda\otimes \mu) & = & 1\otimes \lambda\mu + \lambda_r\otimes \mu - 1\otimes \lambda_r \mu \\ 
& = & 1\otimes (\lambda_E+\lambda_r)\mu + \lambda_r\otimes \mu - 1\otimes \lambda_r \mu  \\ 
& = & 1\otimes \lambda_E\mu + \lambda_r\otimes \mu \\
& = & \lambda_E\otimes \mu + \lambda_r\otimes \mu \\ 
& = & \lambda\otimes \mu.
\end{eqnarray*}  
This implies that $t_0d_0+d_1t_1 = id$. 

Finally, let us prove that $t_id_i+d_{i+1}t_{i+1}=id$ for $i \ge 1$. Since 
\begin{eqnarray*}
& & d_i(\lambda \otimes x_1 \otimes \cdots \otimes x_i \otimes \mu) \\ 
& = &  \lambda x_1\otimes x_2 \otimes \cdots \otimes x_i \otimes \mu  \\ 
& &  + \sum_{j=1}^{i-1} (-1)^j \lambda \otimes x_1 \otimes \cdots \otimes x_jx_{j+1} \otimes \cdots \otimes x_i \otimes \mu + (-1)^{i} \lambda \otimes x_1 \otimes \cdots \otimes x_{i-1} \otimes x_i\mu, 
\end{eqnarray*}
\begin{eqnarray}\label{eq:td}
 & & t_id_i(\lambda \otimes x_1 \otimes \cdots \otimes x_i \otimes \mu) \\
 &  = & 
1 \otimes (\lambda x_1)_r \otimes x_2 \otimes \cdots \otimes x_i \otimes \mu  
+ \sum_{j=1}^{i-1} (-1)^j 1\otimes \lambda_r \otimes x_1 \otimes \cdots \otimes x_jx_{j+1} \otimes \cdots \otimes x_i \otimes \mu \nonumber \\ 
& &  + (-1)^{i} 1 \otimes \lambda_r \otimes x_1 \otimes \cdots \otimes x_{i-1} \otimes x_i\mu. \nonumber 
\end{eqnarray}  
On the other hand, 
\begin{eqnarray}\label{eq:dt} 
& &  d_{i+1}t_{i+1}(\lambda \otimes x_1 \otimes \cdots \otimes x_i \otimes \mu)  \\ 
& = & 
d_{i+1}(1 \otimes \lambda_r \otimes x_1 \otimes \cdots \otimes x_i \otimes \mu) \nonumber \\  
& = &  \lambda_r \otimes x_1 \otimes \cdots \otimes x_i \otimes \mu - 1\otimes \lambda_rx_1 \otimes x_2 \otimes \cdots \otimes x_i \otimes \mu  \nonumber \\ 
& &  + \sum_{j=1}^{i-1} (-1)^{j+1} 1\otimes \lambda_r \otimes x_1 \otimes \cdots \otimes x_j x_{j+1} \otimes \cdots \otimes x_i \otimes \mu  \nonumber \\
& & +(-1)^{i+1} 1\otimes \lambda_r \otimes x_1 \otimes \cdots \otimes x_{i-1} \otimes x_i\mu.  \nonumber 
\end{eqnarray} 
By (\ref{eq:td}) and (\ref{eq:dt}),  
%\begin{flushleft} 
%$(t_id_i+d_{i+1}t_{1+1})(\lambda\otimes x_1 \otimes \cdots \otimes x_i \otimes \mu)$ 
%\end{flushleft}  
%\vspace*{-1.8ex}
\begin{eqnarray*}
& & (t_id_i+d_{i+1}t_{1+1})(\lambda\otimes x_1 \otimes \cdots \otimes x_i \otimes \mu) \\ 
& = & 1 \otimes (\lambda x_1)_r \otimes x_2 \otimes \cdots \otimes x_i \otimes \mu + \lambda_r \otimes x_1 \otimes \cdots \otimes x_i \otimes \mu  - 1\otimes \lambda_rx_1 \otimes x_2 \otimes \cdots \otimes x_i \otimes \mu \\  
& = & 1 \otimes (\lambda_E+\lambda_r)x_1 \otimes x_2 \otimes \cdots \otimes x_i \otimes \mu + \lambda_r \otimes x_1 \otimes \cdots \otimes x_i \otimes \mu  - 1\otimes \lambda_rx_1 \otimes x_2 \otimes \cdots \otimes x_i \otimes \mu \\ 
& =  & 1 \otimes \lambda_E x_1 \otimes x_2 \otimes \cdots \otimes x_i \otimes \mu + \lambda_r \otimes x_1 \otimes \cdots \otimes x_i \otimes \mu \\
& = & \lambda_E \otimes  x_1 \otimes x_2 \otimes \cdots \otimes x_i \otimes \mu + \lambda_r \otimes x_1 \otimes \cdots \otimes x_i \otimes \mu \\
& = & \lambda\otimes x_1 \otimes \cdots \otimes x_i \otimes \mu. 
\end{eqnarray*} 
Here we used $(\lambda x_1)_r = \lambda x_1 = (\lambda_E+\lambda_r)x_1$. 
Hence $t_id_i+d_{i+1}t_{i+1}=id$ for $i \ge 1$. Thus, we have proved that the complex is exact. 
\qed 

\begin{proposition}\label{prop:parabolichochschild} 
Let $R$ be a commutative ring. Let ${\mathcal P}_{\bf n}(R)$ be as in Definition \ref{def:parabolic}. Then 
$H^{i}({\mathcal P}_{\bf n}(R), 
{\rm M}_n(R)/{\mathcal P}_{\bf n}(R)) = 0$ for $i \ge 0$. 
\end{proposition} 

\proof 
Put $P = {\mathcal P}_{\bf n}(R)$. We use the same notation in the proof of Proposition \ref{prop:parabolicprojresol}. 
For $1 \le j < i \le s$, we set 
\[
r_{ij} = \left\{ 
\begin{array}{c|c}  
\left(
\begin{array}{ccccc} 
0 & 0 & 0 & \cdots & 0  \\
X_{21} & 0 & 0 & \cdots & 0 \\
X_{31} & X_{32} & 0 & \cdots & 0 \\
\vdots & \vdots & \ddots & \ddots & \vdots \\
X_{s1} & X_{s2} & X_{s3} & \cdots & 0 \\
\end{array} 
\right) \in {\rm M}_n(R) & 
\begin{array}{c}
X_{ij} \in {\rm M}_{n_i, n_j}(R), \mbox{ but } X_{kl} \mbox{ equals  } 0 \\ 
\mbox{ in } {\rm M}_{n_k, n_l}(R)  \mbox{ for }  (k, l) \neq (i, j)  
\end{array} 
\end{array} 
\right\}.    
\] 
Note that ${\rm M}_n(R)/P \cong \oplus_{1 \le j < i \le s} \; r_{ij}$ as $E$-bimodules. 

Let us consider the projective resolution of $P$ in Proposition \ref{prop:parabolicprojresol}: 
\[
\cdots \to P\otimes_{E}r^{\otimes i} \otimes_{E} P \stackrel{d_i}{\to} P\otimes_{E}r^{\otimes i-1} \otimes_{E} P  
\to \cdots \to P\otimes_{E}r \otimes_{E} P \stackrel{d_1}{\to} P\otimes_{E}P \stackrel{d_0}{\to} P \to 0. 
\] 
To calculate $H^{i}({\mathcal P}_{\bf n}(R), 
{\rm M}_n(R)/{\mathcal P}_{\bf n}(R))$, it suffices to take the cohomology of the following complex 
\begin{multline*} 
0 \to {\rm Hom}_{P^{e}}(P\otimes_{E}P, {\rm M}_n(R)/P) \to 
 {\rm Hom}_{P^{e}}(P\otimes_{E}r \otimes_{E} P, {\rm M}_n(R)/P)  \to \cdots \\ 
 \to {\rm Hom}_{P^{e}}(P\otimes_{E}r^{\otimes i-1} \otimes_{E} P, {\rm M}_n(R)/P)  \to 
 {\rm Hom}_{P^{e}}(P\otimes_{E}r^{\otimes i} \otimes_{E} P, {\rm M}_n(R)/P)  \to \cdots.   
\end{multline*}  
For $i \ge 1$, we see that 
\begin{eqnarray*}
{\rm Hom}_{P^{e}}(P\otimes_{E}r^{\otimes i} \otimes_{E} P, {\rm M}_n(R)/P) & \cong & 
{\rm Hom}_{E^{e}}(r^{\otimes i}, {\rm M}_n(R)/P) \\ 
 & \cong & {\rm Hom}_{E^{e}}(\oplus_{1\le j_1 < j_2 < \cdots < j_{i+1} \le s} \;r_{j_1j_{i+1}}, \oplus_{1 \le j < i \le s} \; r_{ij}) \\ 
 & \cong & 0. 
\end{eqnarray*} 
We also see that 
\begin{eqnarray*}
{\rm Hom}_{P^{e}}(P\otimes_{E} P, {\rm M}_n(R)/P) & \cong & 
{\rm Hom}_{E^{e}}(E, {\rm M}_n(R)/P) \\ 
 & \cong & {\rm Hom}_{E^{e}}(\oplus_{i=1}^{s} r_{ii}, \oplus_{1 \le j < i \le s} \; r_{ij}) \\ 
 & \cong & 0. 
\end{eqnarray*} 
Hence we have $H^{i}({\mathcal P}_{\bf n}(R), 
{\rm M}_n(R)/{\mathcal P}_{\bf n}(R)) = 0$ for $i \ge 0$. 
\qed 

\begin{example}\label{ex:parabolic}\rm 
Proposition \ref{prop:parabolichochschild} is compatible with the fact that the connected component containing ${\mathcal P}_{\bf n} = {\mathcal P}_{n_1, n_2, \ldots, n_s}$ in 
${\rm Mold}_{n, d}$ is isomorphic to ${\rm GL}_n/{\rm P}_{n_1, n_2, \ldots, n_s} \cong {\rm Flag}_{n_1, n_2, \ldots, n_s}$, where 
$d = \sum_{1 \le i \le j \le s} n_in_j$ and 
${\rm P}_{n_1, n_2, \ldots, n_s} = \{  (a_{ij}) \in {\rm GL}_n \mid  a_{ij}=0 \mbox{ if } \sum_{k=1}^{t} n_k < i \le \sum_{k=1}^{t+1} n_k \mbox{ and }  j \le \sum_{k=1}^{t} n_k   \}$ 
(for details, see \cite[Theorem~1.1]{Nakamoto2}). Indeed, since 
$H^{1}({\mathcal P}_{\bf n}(R), 
{\rm M}_n(R)/{\mathcal P}_{\bf n}(R)) = 0$ for any commutative ring $R$, 
the image of the morphism $\phi_{{\mathcal P}_{\bf n}} : {\rm PGL}_{n} \to {\rm Mold}_{n,d}$ 
associated to the mold ${\mathcal P}_{\bf n}({\mathbb Z})$ on ${\rm Spec}\; {\mathbb Z}$ 
is open by Corollary~\ref{cor:openorbit}. 
It is easy to see that 
${\rm Im} \: \phi_{{\mathcal P}_{\bf n}} = {\rm GL}_n/{\rm P}_{n_1, n_2, \ldots, n_s}$ is irreducible, open and closed, and hence that ${\rm GL}_n/{\rm P}_{n_1, n_2, \ldots, n_s}$ is an irreducible component and a connected component.    
We also see that 
$H^{2}({\mathcal P}_{\bf n}(R), 
{\rm M}_n(R)/{\mathcal P}_{\bf n}(R)) = 0$ is compatible with the fact that 
${\rm GL}_n/{\rm P}_{n_1, n_2, \ldots, n_s}$ 
is smooth over ${\Bbb Z}$ (see Theorem~\ref{th:smooth}). 
If ${\mathcal A}(x) = 
{\mathcal P}_{\bf n}(k(x))$ for a point $x \in {\rm Mold}_{n, d}$, then  
\begin{eqnarray*} 
& & \dim_{k(x)} T_{{\rm Mold}_{n, d}/{\Bbb Z}, x} \\ 
& = & \dim_{k(x)}  
H^{1}({\mathcal P}_{\bf n}(k(x)), 
{\rm M}_n(k(x))/{\mathcal P}_{\bf n}(k(x)))  + n^2 
- \dim_{k(x)} N({\mathcal P}_{\bf n}(k(x))) \\ 
& = &  n^2 - \dim_{k(x)} {\mathcal P}_{\bf n}(k(x)) \\ 
& = & \dim {\rm GL}_n(k(x))/{\rm P}_{n_1, n_2, \ldots, n_s}(k(x))
\end{eqnarray*}  
by Corollary \ref{cor:dimtangent}, 
since $N({\mathcal P}_{\bf n}(k(x))) = {\mathcal P}_{\bf n}(k(x))$.  
\end{example} 

\begin{definition}\label{def:xofJn}\rm 
Let $R$ be a commutative ring. We define $x \in {\rm M}_n(R)$ by 
\[ 
x = \left(
\begin{array}{cccccc} 
0 & 1 & 0 & 0 & \cdots & 0 \\
0 & 0 & 1 & 0 & \cdots & 0 \\
0 & 0 & 0 & 1 & \cdots & 0 \\
\vdots & \vdots & \vdots & \ddots & \ddots & \vdots \\
0 & 0 & 0 & 0 & \ddots & 1 \\
0 & 0 & 0 & 0 & \cdots & 0 \\
\end{array} 
\right). 
\] 
Let ${\rm J}_n(R)$ be the $R$-subalgebra of ${\rm M}_n(R)$ generated by $x$. 
Then ${\rm J}_n(R) \cong R[x]/(x^n)$ as $R$-algebras. 
\end{definition}

\bigskip 

To calculate $H^{i}({\rm J}_n(R), 
{\rm M}_n(R)/{\rm J}_n(R))$, we need to make several preparations.  
Let $A = R[x]/(x^n)$. We introduce the following proposition without proof. 
This gives a projective resolution of $A$ over $A^{e} = A\otimes_{R}A^{op}$. 

\begin{proposition}[{\cite[Proposition~1.3]{Buenos}, \cite[Example~2.6]{Redondo}}]\label{prop:projresonegen}   
The following complex gives a projective resolution of $A$ over $A^{e}$: 
\[ 
\cdots \to A^{e} \stackrel{d_n}{\to} A^{e} \to \cdots \to A^{e} \stackrel{d_1}{\to} A^{e} \stackrel{\mu}{\to} A \to 0,  
\]
where 
\[ d_i(a\otimes b) = \left\{ 
\begin{array}{cc}
\left(\sum_{j=0}^{n-1} x^j \otimes x^{n-1-j} \right) (a\otimes b) & (i : \mbox{ even } ) \\
(1\otimes x - x \otimes 1) (a\otimes b) & (i : \mbox{ odd } ) \\
\end{array} 
\right. 
\]  
and $\mu(a\otimes b) = ab$. 
\end{proposition}  

\bigskip 

Set $M = {\rm M}_n(R)/{\rm J}_n(R)$. For calculating $H^{i}({\rm J}_n(R), M)$, it suffices to 
take the cohomology of the complex 
\[
0 \to {\rm Hom}_{A^{e}}(A^{e}, M) \stackrel{d_1^{\ast}}{\to} {\rm Hom}_{A^{e}}(A^{e}, M) \to \cdots 
\to {\rm Hom}_{A^{e}}(A^{e}, M) \stackrel{d_n^{\ast}}{\to} {\rm Hom}_{A^{e}}(A^{e}, M) \to \cdots,   
\]  
which is isomorphic to 
\[ 
0 \to M \stackrel{b^1}{\to} M \to \cdots \to M \stackrel{b^n}{\to} M \to \cdots,  
\] 
where 
\[ b^i(m) = \left\{ 
\begin{array}{cc}
\sum_{j=0}^{n-1} x^j m x^{n-1-j}  & (i : \mbox{ even}) \\
mx - xm  & (i : \mbox{ odd}). \\
\end{array} 
\right. 
\]  
Let $E_{ij} \in {\rm M}_n(R)$ be the matrix unit. 
We can write $x \in {\rm J}_n(R)$ by $x = E_{12}+E_{23}+ \cdots + E_{n-1, n}$. Note that 
\[ E_{ij}x = \left\{ 
\begin{array}{cc}
E_{i, j+1}  & (j \le n-1 ) \\
0  & (j=n) \\
\end{array} 
\right. 
\]  
and 
\[ xE_{ij} = \left\{ 
\begin{array}{cc}
E_{i-1, j}  & (i \ge 2) \\
0  & (i=1). \\
\end{array} 
\right. 
\]  
First, let us calculate $b^i : M \to M$ for even $i$. 
For $E_{kl} \in M$, $b^i(E_{kl}) = \sum_{j=0}^{n-1} x^{j} E_{kl} x^{n-1-j} = \sum_{j=l-1}^{k-1} E_{k-j, l+n-1-j} = \sum_{j=1}^{k-l+1} E_{j, n-1+l-k+j} = x^{n+l-k-1} = 0$ in $M$ 
(if $k<l$, then there is no term in the sum).  Hence if $i$ is even, then $b^i = 0$. 

Next, let us calculate $b^i : M \to M$ for odd $i$. The rank of the $R$-free module $M$ is $n(n-1)$. 
We can choose an $R$-basis $E_{n1}, E_{n2}, \ldots, E_{nn}, E_{n-1, 1}, E_{n-1, 2}, \ldots, E_{n-1, n}, \ldots, E_{2, 1}, E_{2, 2}, \ldots, E_{2, n}$ of $M$. Set $b = b^1 = b^3 = b^5 = \cdots$. 

\begin{lemma}
With respect to the $R$-basis 
\[
E_{n1}, E_{n2}, \ldots, E_{nn}, E_{n-1, 1}, E_{n-1, 2}, \ldots, E_{n-1, n}, \ldots, E_{2, 1}, E_{2, 2}, \ldots, E_{2, n}
\] 
of $M$, the matrix $B$ representing $b: M \to M$ is given by 
\[ 
B = 
\left(
\begin{array}{cccccc}
J & 0 & 0 &  \cdots & 0 & J^{n-1} \\
-I_n & J & 0 & \cdots & 0 &  J^{n-2} \\ 
0 & -I_n & J & \cdots & 0 &  J^{n-3} \\  
\vdots & \vdots & \ddots & \ddots & \vdots & \vdots \\  
0 & 0 & 0 & \ddots & J & J^2 \\ 
0 & 0 & 0 & \cdots & -I_n & 2J \\ 
\end{array} 
\right) \in {\rm M}_{n(n-1)}(R), 
\]
where 
\[ 
J = \left(
\begin{array}{cccccc} 
0 & 0 & 0 &  \cdots & 0 & 0 \\
1 & 0 & 0 &  \cdots & 0 & 0 \\
0 & 1 & 0 &  \cdots & 0 & 0 \\
\vdots & \vdots & \ddots & \ddots & \vdots & \vdots \\
0 & 0 & 0 &  \ddots & 0 & 0 \\
0 & 0 & 0 &  \cdots & 1 & 0 \\
\end{array} 
\right) \in {\rm M}_n(R). 
\] 
\end{lemma} 

\proof 
For $3 \le i \le n$, 
\[ b(E_{ij}) = \left\{ 
\begin{array}{cc}
E_{i, j+1}-E_{i-1, j}  & (j \le n-1) \\
-E_{i-1, n}  & (j=n). \\
\end{array} 
\right. 
\]  
For $1 \le j \le n-2$, $b(E_{2j}) = E_{2, j+1} - E_{1j} = 2E_{2, j+1} + E_{3, j+2} + \cdots + E_{n-j+1, n}$. 
We also see that $b(E_{2, n-1}) = E_{2, n}-E_{1, n-1} = 2E_{2, n}$ and that $b(E_{2n}) = -E_{1n} = 0$.  
By these results, we can check the statement. 
\qed 

\bigskip 

By multiplying 
\[
\left(
\begin{array}{cccccc}
I_n & J & 0 &  \cdots & 0 & 0 \\
0 & I_n & 0 & \cdots & 0 &  0 \\ 
0 & 0 & I_n & \cdots & 0 &  0 \\  
\vdots & \vdots & \ddots & \ddots & \vdots & \vdots \\  
0 & 0 & 0 & \ddots & I_n & 0 \\ 
0 & 0 & 0 & \cdots & 0 & I_n \\ 
\end{array} 
\right) 
\left(
\begin{array}{cccccc}
I_n & 0 & 0 &  \cdots & 0 & 0 \\
0 & I_n & J & \cdots & 0 &  0 \\ 
0 & 0 & I_n & \cdots & 0 &  0 \\  
\vdots & \vdots & \ddots & \ddots & \vdots & \vdots \\  
0 & 0 & 0 & \ddots & I_n & 0 \\ 
0 & 0 & 0 & \cdots & 0 & I_n \\ 
\end{array} 
\right) \cdots 
\left(
\begin{array}{cccccc}
I_n & 0 & 0 &  \cdots & 0 & 0 \\
0 & I_n & 0 & \cdots & 0 &  0 \\ 
0 & 0 & I_n & \cdots & 0 &  0 \\  
\vdots & \vdots & \ddots & \ddots & \vdots & \vdots \\  
0 & 0 & 0 & \ddots & I_n & J \\ 
0 & 0 & 0 & \cdots & 0 & I_n \\ 
\end{array} 
\right) 
\] 
by $B$,  we have 
\[
\left(
\begin{array}{cccccc}
0 & 0 & 0 &  \cdots & 0 & nJ^{n-1} \\
-I_n & 0 & 0 & \cdots & 0 &  (n-1)J^{n-2} \\ 
0 & -I_n & 0 & \cdots & 0 &  (n-2)J^{n-3} \\  
\vdots & \vdots & \ddots & \ddots & \vdots & \vdots \\  
0 & 0 & 0 & \ddots & 0 & 3J^2 \\ 
0 & 0 & 0 & \cdots & -I_n & 2J \\ 
\end{array} 
\right). 
\]
Furthermore, we can obtain the following Smith normal form of $B$ by multiplying elementary matrices (although $R$ may not be a principal ideal domain, we use the terminology "Smith normal form"):  
\[
\left(
\begin{array}{ccccc}
I_n & 0 &  \cdots & 0 & 0 \\ 
0 & I_n &  \cdots & 0 & 0 \\ 
\vdots & \vdots & \ddots & \vdots & \vdots \\ 
0 & 0 & \cdots & I_n & 0 \\ 
0 & 0 & \cdots & 0 & X \\ 
\end{array}
\right), \mbox{ where } X = 
\left(
\begin{array}{cccc}
n & 0 &  \cdots & 0  \\ 
0 & 0 &  \cdots & 0  \\ 
\vdots & \vdots & \ddots & \vdots  \\ 
0 & 0 & \cdots & 0  \\ 
\end{array}
\right) \in {\rm M}_n(R). 
\]   
By the discussion above, we have the following proposition: 

\begin{proposition}
Let ${\rm J}_n(R)$ be as above. Then 
\[ 
H^{i}({\rm J}_n(R), {\rm M}_n(R)/{\rm J}_n(R)) \cong  
\left\{ 
\begin{array}{cc}
R^{n-1} \oplus {\rm Ann}(n) & (i:\mbox{ even }) \\
R^{n-1} \oplus (R/nR) & (i:\mbox{ odd }), \\
\end{array} 
\right. 
\] 
where ${\rm Ann}(n) = \{ a \in R \mid an = 0 \}$. 
\end{proposition} 

\proof 
Let us consider the complex 
\[ 
0 \to M \stackrel{b}{\to} M \stackrel{0}{\to} M \stackrel{b}{\to} M \stackrel{0}{\to} M \stackrel{b}{\to}  \cdots.   
\] 
By the Smith normal form of $B$, ${\rm Ker} B = R^{n-1} \oplus {\rm Ann}(n)$ and 
${\rm Coker} B = R^{n-1} \oplus (R/nR)$.  The statement follows from this result. 
\qed 

\begin{corollary}\label{cor:hochschildjn}
Let $k$ be a field. For each $i \ge 0$,  
\[ 
H^{i}({\rm J}_n(k), {\rm M}_n(k)/{\rm J}_n(k)) \cong  
\left\{ 
\begin{array}{cc}
k^{n-1} & ({\rm ch}(k) \not\;\mid \; n) \\
k^n & ({\rm ch}(k)\; \mid \;n).  \\
\end{array} 
\right. 
\] 
\end{corollary} 

\bigskip 

For calculating the dimension of the tangent space of ${\rm Mold}_{n, n}$ over ${\Bbb Z}$ 
at ${\rm J}_n$,  
we determine the normalizer $N({\rm J}_n(k)) = \{ z \in {\rm M}_n(k) \mid [z, y] \in {\rm J}_n(k) 
\mbox{ for any } y \in {\rm J}_n(k)  \}$ for a field $k$, where $[z, y] = zy-yz$.   

\begin{proposition}\label{prop:normalizerjn}
Let $k$ be a field. Let $x \in {\rm J}_n(k)$ be as in Definition \ref{def:xofJn}.  
Put $A_0 = I_n$ and $A_i = x^{i}$ for $1 \le i \le n-1$.  
We define $B_i, C \in {\rm M}_n(k)$ by  
\[ 
B_i = \sum_{j=1}^{n-i-1} j E_{j+1, i+j+1}= E_{2, i+2} + 2E_{3, i+3} + 3E_{4, i+4} + \cdots + (n-i-1)E_{n-i, n} \; 
\mbox{ for } 0 \le i \le n-2
\] 
and 
\[
C = \sum_{j=1}^{n-1} j E_{j+1, j} =  E_{2, 1} + 2 E_{3, 2}+3E_{4, 3} + \cdots +(n-1) E_{n, n-1}.  
\] 
Then 
\[ 
N({\rm J}_n(k)) = 
\left\{ 
\begin{array}{cc} 
(\oplus_{i=0}^{n-1} k A_i) \oplus (\oplus_{i=0}^{n-2} k B_i)  & ({\rm ch}(k) \not\;\mid \; n) \\
(\oplus_{i=0}^{n-1} k A_i) \oplus (\oplus_{i=0}^{n-2} k B_i) \oplus kC  & ({\rm ch}(k)\; \mid \;n).  \\
\end{array} 
\right. 
\] 
In particular, 
\[ 
\dim_{k} N({\rm J}_n(k)) = 
\left\{ 
\begin{array}{cc} 
2n-1  & ({\rm ch}(k) \not\;\mid \; n) \\
2n  & ({\rm ch}(k)\; \mid \;n).  \\
\end{array} 
\right. 
\] 
\end{proposition} 

\proof 
Note that $N({\rm J}_n(k)) = \{ z \in {\rm M}_n(k) \mid [z, x] \in {\rm J}_n(k) \}$. 
Set ${\rm M}_{\ell} = \oplus_{j-i=l} kE_{ij} \subset {\rm M}_n(k)$. 
Then ${\rm M}_n(k) = \oplus_{l=-(n-1)}^{n-1} {\rm M}_{l}$.  
Since $[E_{ij}, x] = E_{i, j+1}-E_{i-1, j}$, $[z, x] \in {\rm M}_{l+1}$ if $z \in {\rm M}_l$. 
 For $z \in {\rm M}_n(k)$, we can write $z = z_{-(n-1)} + \cdots + z_0 + \cdots + z_{n-1}$, where $z_i \in {\rm M}_{i}$. 
It is easy to see that $z \in N({\rm J}_n(k))$ if and only if $z_i \in N({\rm J}_n(k))$ for $-(n-1) \le i \le n-1$.  It suffices to determine $N({\rm J}_n(k)) \cap{\rm M}_{i}$. 

For $-(n-1) \le i \le -2$, if $0 \neq z_i \in {\rm M}_{i}$, then $0 \neq [z_i, x] \in M_{i+1}$. 
Hence $N({\rm J}_n(k)) \cap{\rm M}_{i} = 0$ for $-(n-1) \le i \le -2$. 
Let $z_{-1} = a_2E_{21}+a_3E_{32}+\cdots + a_n E_{n, n-1} \in {\rm M}_{-1}$. 
Since $[z_{-1}, x] = -a_2E_{11}+(a_2-a_3)E_{22}+(a_3-a_4)E_{33}+\cdots + (a_{n-1}-a_{n})E_{n-1, n-1}+a_nE_{nn}$, $z_{-1} \in N({\rm J}_n(k))$ if and only if 
\begin{eqnarray}\label{eq:z-1condition} 
-a_2=a_2-a_3=a_3-a_4=\cdots = a_{n-1}-a_{n}=a_n.   
\end{eqnarray} 
Suppose that (\ref{eq:z-1condition}) holds. Putting $a_n = -t$, we have  $a_2 = t, a_3=a_2+t, 
a_4=a_3+t, \ldots, a_{n}=a_{n-1}+t$. 
Hence $a_2=t, a_3=2t, \ldots, a_{n-1} = (n-2)t, a_{n}=(n-1)t$. 
By $a_{n}=-t$, we obtain $nt=0$.  If ${\rm ch}(k) \not\;\mid \; n$, then $t=0$. 
In this case, $a_2=a_3=\cdots = a_{n}=0$ and $z_{-1} = 0$. 
If ${\rm ch}(k)\; \mid \;n$, then $z_{-1} = tC$. Conversely, $z_{-1} = tC \in N({\rm J}_n(k))$. Thus, we have 
\[ 
N({\rm J}_n(k)) \cap{\rm M}_{-1} = 
\left\{ 
\begin{array}{cc} 
0  & ({\rm ch}(k) \not\;\mid \; n) \\
kC  & ({\rm ch}(k)\; \mid \;n).  \\
\end{array} 
\right. 
\] 

Let us investigate $N({\rm J}_n(k)) \cap{\rm M}_{0}$.  Let $z_{0} = a_{1}E_{11}+a_{2}E_{22}+\cdots +a_{n}E_{nn} \in M_{0}$.  
Since $[z_{0}, x] = (a_1-a_2)E_{12}+(a_2-a_3)E_{23}+\cdots 
+ (a_{n-1}-a_{n})E_{n-1, n}$, $z_{0} \in N({\rm J}_n(k))$ if and only if 
\begin{eqnarray}\label{eq:z0condition} 
a_1-a_2 = a_2 - a_3 = \cdots = a_{n-1}-a_{n}. 
\end{eqnarray} 
Suppose that (\ref{eq:z0condition}) holds. 
Putting $a_1=s$ and $a_1-a_2=-t$, we have 
$a_1=s, a_2=s+t, a_3=s+2t, \ldots, a_n=s+(n-1)t$. 
Then $z_0 = sI_n+tB_{0} = sA_{0}+tB_{0}$.  Conversely, if $z_0 = sA_{0}+tB_{0}$, then 
$z_0 \in N({\rm J}_n(k))$. Hence $N({\rm J}_n(k)) \cap{\rm M}_{-1} =kA_0 \oplus kB_0$.  
Similarily, we can show that $N({\rm J}_n(k)) \cap{\rm M}_{i} =kA_{i} \oplus kB_{i}$ for $1 \le i \le n-2$ and that 
$N({\rm J}_n(k)) \cap{\rm M}_{n-1} =kA_{n-1}$.   
Therefore, we have proved the statement. 
\qed

\begin{example}\rm 
If ${\mathcal A}(x) = 
{\rm J}_{n}(k(x))$ for a point $x \in {\rm Mold}_{n, n}$, then  
\begin{eqnarray*} 
& & \dim_{k(x)} T_{{\rm Mold}_{n, d}/{\Bbb Z}, x}   \\ 
& = & \dim_{k(x)}  
H^{1}({\rm J}_{n}(k(x)), 
{\rm M}_n(k(x))/{\rm J}_{n}(k(x)))  + n^2 
- \dim_{k(x)} N({\rm J}_{n}(k(x)))  \\ 
& = & n^2 - n
\end{eqnarray*}  
by Corollary \ref{cor:dimtangent}, 
since $\dim_{k}N({\rm J}_{n}(k))  -  \dim_{k}  
H^{1}({\rm J}_{n}(k), 
{\rm M}_n(k)/{\rm J}_{n}(k)) = n$ for any field $k$ by 
Corollary \ref{cor:hochschildjn} and Proposition \ref{prop:normalizerjn}.  
\end{example} 

\bigskip

To calculate several Hochschild cohomology groups,  
%$H^{i}(A, {\rm M}_n(R)/A)$ for $R$-subalgebra $A$ of ${\rm M}_n(R)$ over a commutative ring $R$,  
we introduce several propositions.

\begin{proposition}\label{prop:producthochschild} 
Let $R$ be a commutative ring. 
Let $A$ and $B$ be $R$-subalgebras of ${\rm M}_m(R)$ and ${\rm M}_n(R)$, respectively. 
Assume that $A$ and $B$ are projective $R$-modules. 
We regard the product $A\times B$ as the $R$-subalgebra 
$\left\{ 
\begin{array}{c|} 
\left(\begin{array}{cc} 
X & 0 \\
0 & Y \\
\end{array}
\right) \in  {\rm M}_{m+n}(R) 
\end{array}  \;  X \in A, Y \in B  
\right\}$ of ${\rm M}_{m+n}(R)$.  
Then $H^{i}(A\times B, {\rm M}_{m+n}(R)/(A\times B)) \cong 
H^{i}(A, {\rm M}_m(R)/A) \oplus H^{i}(B, {\rm M}_n(R)/B)$ as $R$-modules for each $i$. 
\end{proposition} 

\proof 
Set $\Lambda = A\times B$. Put $e_A = \left(\begin{array}{cc} 
I_m & 0 \\
0 & 0 \\
\end{array}
\right)$ and 
$e_B = \left(\begin{array}{cc} 
0 & 0 \\
0 & I_n \\
\end{array}
\right)$. 
Note that $e_A$ and $e_B$ are contained in the center of $\Lambda$ and that 
$1= e_A+e_B$, $e_A^2=e_A$, $e_B^2 = e_B$, and $e_Ae_B=e_Be_A=0$. 
There is a projective resolution $P_{\ast} \to A \to 0$ of $A$ in the category of $A$-bimodules. 
There is also a projective resolution $Q_{\ast} \to B \to 0$ of $B$ in the category of $B$-bimodules. 
Then we obtain a projective resolution $P_{\ast}\oplus Q_{\ast} \to 
\Lambda \cong A\oplus B \to 0$ of $\Lambda$ in the category of $\Lambda$-bimodules, 
since $P_{\ast}\oplus Q_{\ast}$ is a projective $\Lambda$-bimodule for each $\ast$.  
Putting $M = {\rm M}_{m+n}(R)/\Lambda = {\rm M}_{m+n}(R)/(A\times B)$, we see that  
\begin{eqnarray*}
M & = & e_A M e_A \oplus e_A M e_B \oplus e_B M e_A \oplus e_B M e_B \\  
   & \cong & {\rm M}_{m}(R)/A \oplus M_{m, n}(R) \oplus M_{n, m}(R) \oplus {\rm M}_n(R)/B,  
\end{eqnarray*} 
where ${\rm M}_{m, n}(R)$ and ${\rm M}_{n, m}(R)$ are the $R$-modules of 
$(m\times n)$-matrices and $(n\times m)$-matrices, respectively. 
By the isomorphism 
\[
{\rm Hom}_{\Lambda^{e}}(P_{\ast}\oplus Q_{\ast}, M) 
\cong {\rm Hom}_{A^{e}}(P_{\ast}, {\rm M}_{m}(R)/A) \oplus 
{\rm Hom}_{B^{e}}(Q_{\ast}, {\rm M}_{n}(R)/B),  
\] 
we have 
$H^{i}(A\times B, {\rm M}_{m+n}(R)/(A\times B)) \cong 
H^{i}(A, {\rm M}_m(R)/A) \oplus H^{i}(B, {\rm M}_n(R)/B)$ for each $i$. 
\qed

\begin{proposition}\label{prop:conjequlhochschild}
Let $A$ be an $R$-subalgebra of ${\rm M}_n(R)$ over a commutative ring $R$. 
Assume that $A$ is a projective module over $R$.  
For $P \in {\rm GL}_n(R)$, set $B = P^{-1}AP$. Then $H^{i}(A, {\rm M}_n(R)/A) \cong 
H^{i}(B, {\rm M}_n(R)/B)$ as $R$-modules for each $i$. 
\end{proposition} 

\proof 
Let $\phi : {\rm M}_n(R) \to {\rm M}_n(R)$ be the isomorphism defined by 
$X \mapsto P^{-1}XP$. 
The commutative diagram 
\[ 
\begin{array}{ccc}
A & \stackrel{\cong}{\to} & B \\ 
\downarrow & & \downarrow \\
{\rm M}_n(R) & \stackrel{\phi}{\to} & {\rm M}_n(R) \\
\end{array} 
\]
is induced by $\phi$. Then we obtain an isomorphism 
$C^{\ast}(A, {\rm M}_n(R)/A) \cong C^{\ast}(B, {\rm M}_n(R)/B)$ of complexes.  
This implies the statement. 
\qed 

\begin{proposition}\label{prop:transequlhochschild}
Let $A$ be an $R$-subalgebra of ${\rm M}_n(R)$ over a commutative ring $R$. 
Assume that $A$ is a projective module over $R$.   
Set ${}^t A = \{ {}^t X \mid X \in A \} \subseteq {\rm M}_n(R)$. 
Then $H^{i}(A, {\rm M}_n(R)/A) \cong 
H^{i}({}^t A, {\rm M}_n(R)/{}^t A)$ as $R$-modules for each $i$. 
\end{proposition} 

\proof 
Let $A^{op}$ be the opposite $R$-algebra of $A$. In other words, $A^{op} = \{ 
a^{op} \mid a \in A \}$ and $a^{op}b^{op} = (ba)^{op}$ for $a, b \in A$. 
For an $A$-bimodule $M$, we define the $A^{op}$-bimodule $M^{op} = \{ m^{op} \mid m \in M \}$ 
by $a^{op} m^{op} b^{op} = (bma)^{op}$ for $a^{op}, b^{op} \in A^{op}$ and $m^{op} \in M^{op}$. 
Let us choose a projective resolution $\cdots \to P_1 \to P_0 \to A \to 0$ of $A$ in the category of $A$-bimodules.  We canonically obtain a projective resolution $\cdots \to P_1^{op} \to P_0^{op} \to A^{op} \to 0$ of $A^{op}$ 
in the category of $A^{op}$-bimodules. Then 
${\rm Hom}_{A-bimod}(P_{\ast}, M)$ and ${\rm Hom}_{A^{op}-bimod}(P_{\ast}^{op}, M^{op})$ are isomorphic as complexes of $R$-modules. 
Hence $H^{i}(A, M) \cong H^{i}(A^{op}, M^{op})$ for each $i$.    

We define a canonical $R$-algebra isomorphism $\phi : A^{op} \to {}^t A$ by $a^{op} \mapsto {}^t a$. 
Note that $A^{op}$ and ${}^t A$ are projective modules over $R$. 
The ${}^t A$-bimodule ${\rm M}_n(R)/{}^t A$ can be regarded as an $A^{op}$-bimodule through $\phi$, which is isomorphic to $({\rm M}_n(R)/A)^{op}$. This implies that 
$H^{i}(A, {\rm M}_n(R)/A) \cong H^{i}(A^{op}, ({\rm M}_n(R)/A)^{op}) \cong 
H^{i}({}^t A, {\rm M}_n(R)/{}^t A)$ for each $i$. 
\qed 

\begin{definition}[{\it cf}. Definitions~\ref{def:locallyequivalent} and \ref{def:equivalentsubalgebra}]\rm
Let $R$ be a commutative ring. Let $A, B \subseteq {\rm M}_n(R)$ be $R$-subalgebras.  
We say that $A$ and $B$ are (globally) {\it equivalent} (or $A \sim B$) 
if there exists $P \in {\rm GL}_n(R)$ such that 
$P^{-1}AP = B$.  
\end{definition} 

\begin{corollary}\label{cor:samefortransposed} 
Let $A$ and $B$ be $R$-subalgebras of ${\rm M}_n(R)$ over a commutative ring $R$. 
Assume that $A$ and $B$ are projective modules over $R$. 
If $A \sim B$ or $A \sim {}^{t}B$, then $H^{i}(A, {\rm M}_n(R)/A) \cong H^{i}(B, {\rm M}_n(R)/B)$ as 
$R$-modules for each $i$. 
\end{corollary} 

\proof 
The statement follows from Propositions \ref{prop:conjequlhochschild} and \ref{prop:transequlhochschild}. 
\qed 

\section{The calculation of  $H^{i}(A, {\rm M}_n(k)/A)$ for $n=2, 3$}  
Let $A$ be a $k$-subalgebra of ${\rm M}_n(k)$ over a field $k$. 
In this section, we discuss the Hochschild cohomology $H^{i}(A, {\rm M}_n(k)/A)$ 
for $n=2$ and $3$. 
For an algebraic closure $\overline{k}$ of $k$, $H^{i}(A\otimes_{k}\overline{k}, 
M_n(\overline{k})/(A\otimes_{k}\overline{k})) \cong 
H^{i}(A, M_{n}(k)/A)\otimes_{k}\overline{k}$ for $i \ge 0$ 
by Proposition \ref{prop:hochextfield}. 
Thereby, 
for studying $H^{i}(A, M_{n}(k)/A)$, it only suffices to investigate 
the case that $k$ is an algebraically closed field. 
For an algebraically closed field $k$, we have the classification of equivalence classes of 
$k$-subalgebras of ${\rm M}_n(k)$ for $n=2, 3$. 
Here we calculate all cases of $k$-subalgebras for $n=2$ and $3$.  

\subsection{The case $n=2$}

In this subsection, we calculate $H^{i}(A, {\rm M}_2(k)/A)$ for 
$k$-subalgebras $A$ of ${\rm M}_2(k)$ over a field $k$. 
In the case $n=2$, we have the following classification. 

\begin{proposition}[{\cite[Proposition~35]{Nakamoto-Torii:50th} and \cite[Proposition~2.2]{Nakamoto-Torii:classification}}]\label{prop:deg2} 
Let $k$ be an algebraically closed field. 
Any subalgebras of ${\rm M}_2(k)$ are equivalent to  
one of the following: 
\begin{enumerate}
\item\label{item:2-1} ${\rm M}_2(k)$ 
\item\label{item:2-2} ${\rm B}_2(k) = \left\{ \left( 
\begin{array} {cc}
\ast & \ast \\
0 & \ast \\
\end{array} 
\right)
\right\}$ 
\item\label{item:2-3} ${\rm D}_2(k) = \left\{ \left( 
\begin{array} {cc}
\ast & 0 \\
0 & \ast \\
\end{array} 
\right)
\right\}$
\item\label{item:2-4} ${\rm N}_2(k) = \left\{ 
\begin{array}{c|c} 
\left( 
\begin{array} {cc}
a & b \\
0 & a \\
\end{array} 
\right) & a, b \in k 
\end{array}  
\right\}$
\item\label{item:2-5} ${\rm C}_2(k) = \left\{ 
\begin{array}{c|c} 
\left( 
\begin{array} {cc}
a & 0 \\
0 & a \\
\end{array} 
\right) & a \in k 
\end{array}  
\right\}$. 
\end{enumerate} 
\end{proposition}

\bigskip 

Let $k$ be a (not necessarily algebraically closed) field. 
We summarize the results on $H^{i}(A, {\rm M}_2(k)/A)$ in the cases (\ref{item:2-1})--(\ref{item:2-5}) in Proposition \ref{prop:deg2}. For details, see Table \ref{table:deg2} in Section 6. 

\begin{enumerate} 
\item[(\ref{item:2-1})]  For $A={\rm M}_2(k)$, we have $H^{i}(A, {\rm M}_2(k)/A)=0$ for $i \ge 0$ by 
Example \ref{ex:fullmatrix}. 
\item[(\ref{item:2-2})]  For $A = {\rm B}_2(k)$, we have $H^{i}(A, {\rm M}_2(k)/A)=0$  
for $i \ge 0$ by Example \ref{ex:borel}. 
\item[(\ref{item:2-3})] For $A = {\rm D}_2(k)$, we have $H^{i}(A, {\rm M}_2(k)/A)=0$  
for $i \ge 0$ by Example \ref{ex:dn}. 
\item[(\ref{item:2-4})]  For $A = {\rm N}_2(k)$, $A$ coincides with ${\rm J}_2(k)$ in 
Definition \ref{def:xofJn}. Then we have 
\[ 
H^{i}(A, {\rm M}_2(k)/A) \cong  
\left\{ 
\begin{array}{cc}
k & ({\rm ch}(k) \neq 2) \\
k^2 & ({\rm ch}(k) = 2).  \\
\end{array} 
\right. 
\] for $i \ge 0$ by Corollary \ref{cor:hochschildjn}. 
\item[(\ref{item:2-5})]  For $A={\rm C}_2(k)$, we have 
\[
H^{i}(A, {\rm M}_2(k)/A) \cong \left\{  
\begin{array}{cl} 
{\rm M}_2(k)/{\rm C}_2(k) \cong k^3 & (i=0) \\
0 & (i>0). 
\end{array}
\right. 
\] 
by Example \ref{ex:scalar}. 
\end{enumerate}

\subsection{The case $n=3$}
In this subsection, we calculate $H^{i}(A, {\rm M}_3(k)/A)$ for 
$k$-subalgebras $A$ of ${\rm M}_3(k)$ over a field $k$. 
In the case $n=3$, we have the following classification. 

\begin{theorem}[{\cite[Theorem~2]{Nakamoto-Torii:50th} and \cite[Theorem~2.1]{Nakamoto-Torii:classification}}]\label{th:deg3} 
Let $k$ be an algebraically closed field. 
Any subalgebras of ${\rm M}_3(k)$ are equivalent to  
one of the following: 
\begin{enumerate}
\item\label{item:3-1} ${\rm M}_3(k)$ \\
\item\label{item:3-2} ${\rm P}_{2, 1}(k) = \left\{ \left(
\begin{array}{ccc}
\ast & \ast & \ast \\
\ast & \ast & \ast \\
0 & 0 & \ast \\ 
\end{array}
\right)  \in {\rm M}_3(k) \right\}$ \\
\item\label{item:3-3} ${\rm P}_{1, 2}(k) = \left\{ \left(
\begin{array}{ccc}
\ast & \ast & \ast \\
0 & \ast & \ast \\
0 & \ast & \ast \\ 
\end{array}
\right)  \in {\rm M}_3(k) \right\}$ \\
\item\label{item:3-4} ${\rm B}_{3}(k) = \left\{ \left(
\begin{array}{ccc}
\ast & \ast & \ast \\
0 & \ast & \ast \\
0 & 0 & \ast \\ 
\end{array}
\right)  \in {\rm M}_3(k) \right\}$ \\
\item\label{item:3-5} ${\rm C}_3(k) = \left\{  \begin{array}{c|c}
\left(
\begin{array}{ccc}
a & 0 & 0 \\
0 & a & 0 \\
0 & 0 & a \\ 
\end{array}
\right)  & a \in k 
\end{array} 
\right\}$ \\
\item\label{item:3-6} ${\rm D}_{3}(k) = \left\{ \left(
\begin{array}{ccc}
\ast & 0 & 0 \\
0 & \ast & 0 \\
0 & 0 & \ast \\ 
\end{array}
\right)  \in {\rm M}_3(k) \right\}$ \\
\item\label{item:3-7} $({\rm C}_{2} \times {\rm D}_{1})(k) = \left\{ 
\begin{array}{c|c} 
\left(
\begin{array}{ccc}
a & 0 & 0 \\
0 & a & 0 \\
0 & 0 & b \\ 
\end{array}
\right) 
& a, b \in k 
\end{array} \right\}$ \\
\item\label{item:3-8} $({\rm N}_{2} \times {\rm D}_{1})(k) = \left\{ 
\begin{array}{c|c} 
\left(
\begin{array}{ccc}
a & c & 0 \\
0 & a & 0 \\
0 & 0 & b \\ 
\end{array}
\right) 
& a, b, c \in k 
\end{array} \right\}$ \\
\item\label{item:3-9} $({\rm B}_{2} \times {\rm D}_{1})(k) = \left\{  
\left(
\begin{array}{ccc}
* & * & 0 \\
0 & * & 0 \\
0 & 0 & * \\ 
\end{array}
\right)  \in {\rm M}_3(k) \right\}$ \\
\item\label{item:3-10} $({\rm M}_{2} \times {\rm D}_{1})(k) = \left\{  
\left(
\begin{array}{ccc}
* & * & 0 \\
* & * & 0 \\
0 & 0 & * \\ 
\end{array}
\right)  \in {\rm M}_3(k) \right\}$ \\
\item\label{item:3-11} ${\rm J}_{3}(k) = \left\{ 
\begin{array}{c|c} 
\left(
\begin{array}{ccc}
a & b & c \\
0 & a & b \\
0 & 0 & a \\ 
\end{array}
\right) 
& a, b, c \in k 
\end{array} \right\}$ \\
\item\label{item:3-12} ${\rm N}_3(k)  = \left\{ 
\begin{array}{c|c} 
\left(
\begin{array}{ccc}
a & b & c \\
0 & a & d \\
0 & 0 & a \\ 
\end{array}
\right) 
& a, b, c, d \in k 
\end{array} \right\}$ \\
\item\label{item:3-13} $S_{1}(k)  = \left\{ 
\begin{array}{c|c} 
\left(
\begin{array}{ccc}
a & b & 0 \\
0 & a & 0 \\
0 & 0 & a \\ 
\end{array}
\right) 
& a, b \in k 
\end{array} \right\}$ \\
\item\label{item:3-14} $S_{2}(k)  = \left\{ 
\begin{array}{c|c} 
\left(
\begin{array}{ccc}
a & 0 & 0 \\
0 & a & c \\
0 & 0 & b \\ 
\end{array}
\right) 
& a, b, c \in k 
\end{array} \right\}$ \\
\item\label{item:3-15} $S_{3}(k)  = \left\{ 
\begin{array}{c|c} 
\left(
\begin{array}{ccc}
a & 0 & c \\
0 & b & 0 \\
0 & 0 & b \\ 
\end{array}
\right) 
& a, b, c \in k 
\end{array} \right\}$ \\
\item\label{item:3-16} $S_{4}(k)  = \left\{ 
\begin{array}{c|c} 
\left(
\begin{array}{ccc}
a & b & c \\
0 & a & 0 \\
0 & 0 & a \\ 
\end{array}
\right) 
& a, b, c \in k 
\end{array} \right\}$ \\
\item\label{item:3-17} $S_{5}(k)  = \left\{ 
\begin{array}{c|c} 
\left(
\begin{array}{ccc}
a & 0 & b \\
0 & a & c \\
0 & 0 & a \\ 
\end{array}
\right) 
& a, b, c \in k 
\end{array} \right\}$ \\
%\item $S_{5}(u, v)  = \left\{ 
%\begin{array}{c|c} 
%\left(
%\begin{array}{ccc}
%a & cu & 0 \\
%0 & a+cv & 0 \\
%0 & 0 & b \\ 
%\end{array}
%\right) 
%& a, b, c \in k 
%\end{array} \right\}$ \\
\item\label{item:3-18} $S_{6}(k)  = \left\{ 
\begin{array}{c|c} 
\left(
\begin{array}{ccc}
a & c & d \\
0 & a & 0 \\
0 & 0 & b \\ 
\end{array}
\right) 
& a, b, c, d \in k 
\end{array} \right\}$ \\
%%\item $S_{6}  = \left\{ 
%%\begin{array}{c|c} 
%%\left(
%%\begin{array}{ccc}
%%a & c & 0 \\
%%0 & a & 0 \\
%%0 & d & b \\ 
%%\end{array}
%%\right) 
%%& a, b, c, d \in k 
%%\end{array} \right\}$ \\
%
%  dim 4 case 4 
%
%%\item $S_{8}  = \left\{ 
%%\begin{array}{c|c} 
%%\left(
%%\begin{array}{ccc}
%%a & c & d \\
%%0 & a & 0 \\
%%0 & 0 & b \\ 
%%\end{array}
%%\right) 
%%& a, b, c, d \in k 
%%\end{array} \right\}$ \\
\item\label{item:3-19} $S_{7}(k)  = \left\{ 
\begin{array}{c|c} 
\left(
\begin{array}{ccc}
a & 0 & c \\
0 & a & d \\
0 & 0 & b \\ 
\end{array}
\right) 
& a, b, c, d \in k 
\end{array} \right\}$ \\
\item\label{item:3-20}$S_{8}(k) = \left\{ 
\begin{array}{c|c} 
\left(
\begin{array}{ccc}
a & c & d \\
0 & b & 0 \\
0 & 0 & b \\ 
\end{array}
\right) 
& a, b, c, d \in k 
\end{array} \right\}$ \\
\item\label{item:3-21} $S_{9}(k)  = \left\{ 
\begin{array}{c|c} 
\left(
\begin{array}{ccc}
a & 0 & c \\
0 & b & d \\
0 & 0 & b \\ 
\end{array}
\right) 
& a, b, c, d \in k 
\end{array} \right\}$ \\
%
%
%  dim 5 case
%
\item\label{item:3-22} $S_{10}(k)  = \left\{ 
\begin{array}{c|c} 
\left(
\begin{array}{ccc}
a & b & c \\
0 & a & d \\
0 & 0 & e \\ 
\end{array}
\right) 
& a, b, c, d, e \in k 
\end{array} \right\}$ \\
\item\label{item:3-23} $S_{11}(k) = \left\{ 
\begin{array}{c|c} 
\left(
\begin{array}{ccc}
a & b & c \\
0 & e & d \\
0 & 0 & a \\ 
\end{array}
\right) 
& a, b, c, d, e \in k 
\end{array} \right\}$ \\
\item\label{item:3-24} $S_{12}(k)  = \left\{ 
\begin{array}{c|c} 
\left(
\begin{array}{ccc}
a & b & c \\
0 & e & d \\
0 & 0 & e \\ 
\end{array}
\right) 
& a, b, c, d, e \in k 
\end{array} \right\}$ \\
\item\label{item:3-25} $S_{13}(k)  = \left\{  
\left(
\begin{array}{ccc}
\ast & \ast & \ast \\
0 & \ast & 0 \\
0 & 0 & \ast \\ 
\end{array}
\right) \in {\rm M}_3(k) 
 \right\}$ \\
\item\label{item:3-26} $S_{14}(k)  = \left\{  
\left(
\begin{array}{ccc}
\ast & 0 & \ast \\
0 & \ast & \ast \\
0 & 0 & \ast \\ 
\end{array}
\right) 
\in {\rm M}_3(k) 
 \right\}$ \\
\end{enumerate}
\end{theorem}

Let $k$ be a (not necessarily algebraically closed) field. 
We summarize the results on $H^{i}(A, {\rm M}_2(k)/A)$ in the cases (\ref{item:3-1})--(\ref{item:3-26}) in Theorem  \ref{th:deg3}. For details, see Table \ref{table:deg3} in Section 6. 

\begin{enumerate} 
\item[(\ref{item:3-1})]  For $A={\rm M}_3(k)$, we have $H^{i}(A, {\rm M}_3(k)/A)=0$ for $i \ge 0$ by Example \ref{ex:fullmatrix}. 
\item[(\ref{item:3-2})]  For $A={\rm P}_{2, 1}(k)$, we have $H^{i}(A, {\rm M}_3(k)/A)=0$ for $i \ge 0$ by 
Proposition \ref{prop:parabolichochschild}.  
\item[(\ref{item:3-3})]  For $A={\rm P}_{1, 2}(k)$, we have $H^{i}(A, {\rm M}_3(k)/A)=0$ for $i \ge 0$ by 
Proposition \ref{prop:parabolichochschild}.  
\item[(\ref{item:3-4})]  For $A={\rm B}_{3}(k)$, we have $H^{i}(A, {\rm M}_3(k)/A)=0$ for $i \ge 0$ by 
Example \ref{ex:borel}.  
\item[(\ref{item:3-5})]  For $A={\rm C}_3(k)$, we have 
\[
H^{i}(A, {\rm M}_3(k)/A) \cong \left\{  
\begin{array}{cl} 
{\rm M}_3(k)/{\rm C}_3(k) \cong k^8 & (i=0) \\
0 & (i>0). 
\end{array}
\right. 
\] 
by Example \ref{ex:scalar}. 
\item[(\ref{item:3-6})] For $A = {\rm D}_3(k)$, we have $H^{i}(A, {\rm M}_3(k)/A)=0$ for    
$i \ge 0$ by Example \ref{ex:dn}.  
\item[(\ref{item:3-7})] For $A = ({\rm C}_2 \times {\rm D}_1)(k)$, we have 
\[
H^{i}(A, {\rm M}_3(k)/A) \cong \left\{  
\begin{array}{cl} 
{\rm M}_2(k)/{\rm C}_2(k) \cong k^3 & (i=0) \\
0 & (i>0). 
\end{array}
\right. 
\] 
Indeed,  
$H^{i}(A, {\rm M}_3(k)/A) \cong H^{i}({\rm C}_2(k), {\rm M}_2(k)/{\rm C}_2(k))\oplus 
H^{i}({\rm D}_1(k), {\rm M}_1(k)/{\rm D}_1(k))$ 
by Proposition \ref{prop:producthochschild}. 
By Example \ref{ex:fullmatrix} (or by Example \ref{ex:dn}), $H^{i}({\rm D}_1(k), {\rm M}_1(k)/{\rm D}_1(k)) = 0$. 
Hence we can calculate $H^{i}(A, {\rm M}_3(k)/A)$ by using the result on  $H^{i}({\rm C}_2(k), {\rm M}_2(k)/{\rm C}_2(k))$.  
\item[(\ref{item:3-8})] For $A=({\rm N}_{2} \times {\rm D}_{1})(k)$, we have 
\[ 
H^{i}(A, {\rm M}_3(k)/A) \cong 
\left\{ 
\begin{array}{cc}
k & ({\rm ch}(k) \neq 2) \\
k^2 & ({\rm ch}(k) = 2).  \\
\end{array} 
\right. 
\]
for each $i$. 
Indeed,  
$H^{i}(A, {\rm M}_3(k)/A) \cong H^{i}({\rm N}_2(k), {\rm M}_2(k)/{\rm N}_2(k))\oplus 
H^{i}({\rm D}_1(k), {\rm M}_1(k)/{\rm D}_1(k))$ 
by Proposition \ref{prop:producthochschild}. 
By Example \ref{ex:fullmatrix} (or by Example \ref{ex:dn}), $H^{i}({\rm D}_1(k), {\rm M}_1(k)/{\rm D}_1(k)) = 0$. 
Hence we can calculate $H^{i}(A, {\rm M}_3(k)/A)$ by using the result on  $H^{i}({\rm N}_2(k), {\rm M}_2(k)/{\rm N}_2(k))$.  
\item[(\ref{item:3-9})] For $A=({\rm B}_{2} \times {\rm D}_{1})(k)$, we have 
$H^{i}(A, {\rm M}_3(k)/A)=0$ for $i \ge 0$. 
Indeed, $H^{i}(A, {\rm M}_3(k)/A) \cong H^{i}({\rm B}_2(k), {\rm M}_2(k)/{\rm B}_2(k))\oplus 
H^{i}({\rm D}_1(k), {\rm M}_1(k)/{\rm D}_1(k)) = 0$  
by Proposition \ref{prop:producthochschild}, Examples \ref{ex:borel} and \ref{ex:fullmatrix} (or Example \ref{ex:dn}). 
\item[(\ref{item:3-10})] For $A=({\rm M}_{2} \times {\rm D}_{1})(k)$, we have 
$H^{i}(A, {\rm M}_3(k)/A)=0$ for $i \ge 0$. 
Indeed, $H^{i}(A, {\rm M}_3(k)/A) \cong H^{i}({\rm M}_2(k), {\rm M}_2(k)/{\rm M}_2(k))\oplus 
H^{i}({\rm D}_1(k), {\rm M}_1(k)/{\rm D}_1(k)) = 0$  
by Proposition \ref{prop:producthochschild} and Example \ref{ex:fullmatrix} (or Example \ref{ex:dn}). 
\item[(\ref{item:3-11})] For $A = {\rm J}_3(k)$, we have 
\[ 
H^{i}({\rm J}_3(k), {\rm M}_3(k)/{\rm J}_3(k)) \cong  
\left\{ 
\begin{array}{cc}
k^{2} & ({\rm ch}(k) \neq 3) \\
k^3 & ({\rm ch}(k) = 3)  \\
\end{array} 
\right. 
\] 
for $i \ge 0$ by Corollary \ref{cor:hochschildjn}. 
\item[(\ref{item:3-12})] For $A = {\rm N}_3(k)$, we have 
\[
H^{i}(A, {\rm M}_3(k)/A) \cong \left\{  
\begin{array}{cl} 
k^2 & (i=0) \\
k^{i+1} & (i>0). 
\end{array}
\right. 
\]
For details, see Section \ref{subsection:N3}. 
\item[(\ref{item:3-13})] For $A = {\rm S}_1(k)$, we have 
\[
H^{i}(A, {\rm M}_3(k)/A) \cong \left\{  
\begin{array}{cl} 
k^4 & (i=0) \\
k & (i>0). 
\end{array}
\right. 
\]
For details, see Section \ref{subsection:S1}. 
\item[(\ref{item:3-14})] For $A = {\rm S}_2(k)$, we have 
\[
H^{i}(A, {\rm M}_3(k)/A) \cong \left\{  
\begin{array}{cl} 
k^2 & (i=0) \\
0 & (i>0). 
\end{array}
\right. 
\] 
For details, see Section \ref{subsection:S2}. 
\item[(\ref{item:3-15})] For $A = {\rm S}_3(k)$, we have 
\[
H^{i}(A, {\rm M}_3(k)/A) \cong \left\{  
\begin{array}{cl} 
k^2 & (i=0) \\
0 & (i>0)  
\end{array}
\right. 
\] 
by the result on $H^{i}({\rm S}_2(k), {\rm M}_3(k)/{\rm S}_{2}(k))$ and Corollary \ref{cor:samefortransposed}, since ${\rm S}_3(k) \sim {}^{t}{\rm S}_{2}(k)$. 
\item[(\ref{item:3-16})] For $A = {\rm S}_4(k)$, we have 
\[
H^{i}(A, {\rm M}_3(k)/A) \cong \left\{  
\begin{array}{cl} 
k^4 & (i=0) \\
k^{3\times 2^i} & (i>0). 
\end{array}
\right. 
\] 
For details, see Section \ref{subsection:S4}. 
\item[(\ref{item:3-17})] For $A = {\rm S}_5(k)$, we have 
\[
H^{i}(A, {\rm M}_3(k)/A) \cong \left\{  
\begin{array}{cl} 
k^4 & (i=0) \\
k^{3\times 2^i} & (i>0)  
\end{array}
\right. 
\] 
by the result on $H^{i}({\rm S}_4(k), {\rm M}_3(k)/{\rm S}_{4}(k))$ and Corollary \ref{cor:samefortransposed}, since ${\rm S}_5(k) \sim {}^{t}{\rm S}_{4}(k)$. 
\item[(\ref{item:3-18})] For $A = {\rm S}_6(k)$, we have $H^{i}(A, {\rm M}_3(k)/A)\cong k$ for $i \ge 0$. 
For details, see Section \ref{subsection:S6}. 
\item[(\ref{item:3-19})] For $A = {\rm S}_7(k)$, we have 
\[
H^{i}(A, {\rm M}_3(k)/A) \cong \left\{  
\begin{array}{cl} 
k^3 & (i=0) \\
0 & (i>0). 
\end{array}
\right. 
\] 
For details, see Section \ref{subsection:S7}. 
\item[(\ref{item:3-20})] For $A = {\rm S}_8(k)$, we have 
\[
H^{i}(A, {\rm M}_3(k)/A) \cong \left\{  
\begin{array}{cl} 
k^3 & (i=0) \\
0 & (i>0) 
\end{array}
\right. 
\] 
by the result on $H^{i}({\rm S}_7(k), {\rm M}_3(k)/{\rm S}_{7}(k))$ and Corollary \ref{cor:samefortransposed}, since ${\rm S}_8(k) \sim {}^{t}{\rm S}_{7}(k)$. 
\item[(\ref{item:3-21})] For $A = {\rm S}_9(k)$, we have $H^{i}(A, {\rm M}_3(k)/A)\cong k$ for $i \ge 0$. Indeed,  
this follows from the result on $H^{i}({\rm S}_6(k), {\rm M}_3(k)/{\rm S}_{6}(k))$ and Corollary \ref{cor:samefortransposed}, since ${\rm S}_9(k) \sim {}^{t}{\rm S}_{6}(k)$. 
\item[(\ref{item:3-22})] For $A = {\rm S}_{10}(k)$, we have 
\[ 
H^{i}(A, {\rm M}_3(k)/A) \cong  
\left\{ 
\begin{array}{cc}
k & ({\rm ch}(k) \neq 2) \\
k^2 & ({\rm ch}(k) = 2)  \\
\end{array} 
\right. 
\] 
for $i \ge 0$. 
For details, see Section \ref{subsection:S10}. 
\item[(\ref{item:3-23})] For $A = {\rm S}_{11}(k)$, we have 
\[
H^{i}(A, {\rm M}_3(k)/A) \cong \left\{  
\begin{array}{cl} 
k & (i=0, 1) \\
0 & (i\ge 2). 
\end{array}
\right. 
\] 
For details, see Section \ref{subsection:S11}. 
\item[(\ref{item:3-24})] For $A = {\rm S}_{12}(k)$, we have 
\[ 
H^{i}(A, {\rm M}_3(k)/A) \cong  
\left\{ 
\begin{array}{cc}
k & ({\rm ch}(k) \neq 2) \\
k^2 & ({\rm ch}(k) = 2)  \\
\end{array} 
\right. 
\] 
for $i \ge 0$. 
Indeed, this follows from the result on $H^{i}({\rm S}_{10}(k), {\rm M}_3(k)/{\rm S}_{10}(k))$ and Corollary \ref{cor:samefortransposed}, since ${\rm S}_{12}(k) \sim {}^{t}{\rm S}_{10}(k)$. 
\item[(\ref{item:3-25})]  For $A={\rm S}_{13}(k)$, we have $H^{i}(A, {\rm M}_3(k)/A)=0$ for $i \ge 0$. 
For details, see Section \ref{subsection:S13}. 
\item[(\ref{item:3-26})]  For $A={\rm S}_{14}(k)$, we have $H^{i}(A, {\rm M}_3(k)/A)=0$ for $i \ge 0$. 
Indeed, this follows from the result on $H^{i}({\rm S}_{13}(k), {\rm M}_3(k)/{\rm S}_{13}(k))$ and Corollary \ref{cor:samefortransposed}, since ${\rm S}_{14}(k) \sim {}^{t}{\rm S}_{13}(k)$. 
\end{enumerate} 

\subsection{The case $A = {\rm N}_{3}(k)$}\label{subsection:N3}
Set 
${\rm N}_3(R)  = \left\{ 
\begin{array}{c|c} 
\left(
\begin{array}{ccc}
a & b & c \\
0 & a & d \\
0 & 0 & a \\ 
\end{array}
\right) 
& a, b, c, d \in R 
\end{array} \right\}$ for a commutative ring $R$.  
We denote ${\rm N}_3(R)$ by $N$ for simplicity.
Let $J$ be the two-sided ideal of $N$ given by 
\[
J  = \left\{ 
\begin{array}{c|c} 
\left(
\begin{array}{ccc}
0 & b & c \\
0 & 0 & d \\
0 & 0 & 0 \\ 
\end{array}
\right) \in N 
& b, c, d \in R 
\end{array} \right\}.
\] 
We set $T=N/J$,
which is an $N$-bimodule over $R$.
First, we calculate
the Hochschild cohomology 
$H^\ast(N, T)$ of $N$ with coefficients in $T$. 
We note that
there is an isomorphism
$T\otimes_N T\cong T$ of $N$-bimodules over $R$.
This implies that
$T$ is a monoid object in the category of $N$-bimodules
over $R$.
The unit $u: N\to T=N/J$ is given by the projection.
Hence $H^\ast(N, T)$ has the structure of
a graded associative algebra over $R$.

We set $I=E_{11}+E_{22}+E_{33}$, 
$U=E_{12},V=E_{23}$ and $W=E_{13}$.
We let $\overline{N}=N/RI$.
The set $\{U,V,W\}$ forms a basis of the free $R$-module 
$\overline{N}$.
Let $\overline{B}_{\ast}(N,N,N)$ be 
the reduced bar resolution of $N$ as $N$-bimodules
over $R$.
We have
\[ \overline{B}_p(N,N,N)\cong
   N\otimes_R \overbrace{\overline{N}\otimes_R\cdots\otimes_R
   \overline{N}}^p\otimes_R N\]
for $p\ge 0$.
We denote  the cochain complex
${\rm Hom}_{N^e}(\overline{B}_{\ast}(N,N,N),T)$
by $\overline{C}^{\ast}(N, T)$. 
%Since the reduced bar resolution $\overline{B}_{\ast}(N,N,N)$
%is chain homotopy equivalent to
%the bar resolution $B_{\ast}(N,N,N)$,
The cohomology of $\overline{C}^{\ast}(N, T)$
is isomorphic to the Hochschild cohomology
$H^{\ast}(N, T)$ (see Section 2). 

Let $e\in T$ be the image of $I$ under the unit $u: N\to T$.
We denote by $U^{\ast},V^{\ast}\in \overline{C}^1(N, T)$ the maps $\overline{N} \to T$
of $R$-modules given by
\[ \begin{array}{rcl}
   U^{\ast}(n)&=&
   \left\{\begin{array}{ll}
     e & \mbox{\rm if $n=U$}\\
     0 & \mbox{\rm if $n=V,W$} 
     \end{array}\right.\\[3mm]
   V^{\ast}(n)&=& 
   \left\{\begin{array}{ll}
     e & \mbox{\rm if $n=V$}\\
     0 & \mbox{\rm if $n=U,W$},
     \end{array}\right.
   \end{array}\]
respectively.
The maps $U^{\ast}$ and $V^{\ast}$
are $1$-cocycles in the cochain complex
$\overline{C}^{\ast}(N, T)$.
We denote by
\[ \alpha,\beta\in H^1(N, T)\] 
the cohomology classes represented by the $1$-cocycles $U^{\ast},V^{\ast}$,
respectively.

Let $W^{\ast}\in \overline{C}^1(N, T)$ be the map $\overline{N}\to T$ 
of $R$-modules given by
\[  W^{\ast}(n)= 
   \left\{\begin{array}{ll}
     e & \mbox{\rm if $n=W$}\\
     0 & \mbox{\rm if $n=U,V$}.
     \end{array}\right.\]
We observe that
\[ U^{\ast}\cup V^{\ast}=-\delta^1(W^{\ast}),\]
where $\delta^1: \overline{C}^1(N, T) \to \overline{C}^2(N, T)$ 
is the coboundary map.
Thus, we obtain that
\[ \alpha\beta=0 \]
in $H^2(N, T)$.

Let $R\langle\alpha,\beta\rangle$
be the free associative algebra over $R$
generated by $\alpha$ and $\beta$.
There is a map 
\[ R\langle\alpha,\beta\rangle/(\alpha\beta)\longrightarrow
   H^\ast(N, T) \]
of graded associative algebras over $R$,
where $(\alpha\beta)$
is the two-sided ideal of $R\langle\alpha,\beta\rangle$
generated by $\alpha\beta$.

\begin{lemma}
We have an isomorphism
$H^\ast(N, T)\cong R\langle\alpha,\beta\rangle/(\alpha\beta)$
of graded associative algebras over $R$.
\end{lemma}

\proof
%{\color{red} We have to change this proof by
%using cochain complexes.}
We observe that the cochain complex
$\overline{C}^{\ast}(N, T)$
is isomorphic to 
the differential graded algebra which is
the free associative algera
\[ R\langle U^{\ast},V^{\ast},W^{\ast}\rangle \]
generated by $U^{\ast},V^{\ast},W^{\ast}$
with differential
\[ \delta(U^{\ast})=\delta(V^{\ast})=0,\quad \delta(W^*)=-U^*V^* .\]
We let $C_{VU}^{\ast}$ be the subcomplex of $R\langle U^{\ast},V^{\ast},W^{\ast}\rangle$
given by 
\[ C_{VU}^{\ast}=\bigoplus_{n=0}^{\infty} \bigoplus_{\tiny 
\begin{array}{c} i+j=n \\ i, j \ge 0 \end{array} }
   R\overbrace{V^{\ast}\cdots V^{\ast}}^i\overbrace{U^{\ast}\cdots U^{\ast}}^j  \]
with trivial differential.
We also let $C_W^{\ast}$ be the subcomplex of $R\langle U^{\ast},V^{\ast},W^{\ast}\rangle$
given by
\[ \begin{array}{rcl}
     C_W^i&=&\left\{ \begin{array}{ll}
                     RW^{\ast} & (i=1),\\
                     RU^{\ast}V^{\ast} & (i=2),\\
                     0             & (i\neq 1,2),\\
                     \end{array}\right.
   \end{array}\]
with differential $\delta(W^{\ast})=-U^{\ast}V^{\ast}$.
We observe that there is an isomorphism of cochain complexes
between $R\langle U^{\ast},V^{\ast},W^{\ast}\rangle$
and 
\[ C_{VU}^{\ast}\otimes_R\bigoplus_{r\ge 0}\overbrace{(C_W^{\ast}\otimes_R C_{VU}^{\ast})\otimes_R\cdots \otimes_R(C_W^{\ast}\otimes_R C_{VU}^{\ast})}^r,  \] 
where  we set $\overbrace{(C_W^{\ast}\otimes_R C_{VU}^{\ast})\otimes_R\cdots \otimes_R(C_W^{\ast}\otimes_R C_{VU}^{\ast})}^r = (0 \to R \stackrel{\delta^{0}}{\to} 0 \stackrel{\delta^{1}}{\to} 0 \stackrel{\delta^{2}}{\to} \cdots )$ if $r=0$. 
Since $C_W^{\ast}$ is acyclic 
and $C_{VU}^*$ has a trivial differential,
we obtain an isomorphism of $R$-modules
\[ H^{\ast}(\overline{C}(N, T))\cong  
\bigoplus_{n=0}^{\infty} \bigoplus_{\tiny 
\begin{array}{c} i+j=n \\ i, j \ge 0 \end{array} }
   R\overbrace{V^{\ast}\cdots V^{\ast}}^i\overbrace{U^{\ast}\cdots U^{\ast}}^j     
%   \bigoplus_{i+j\ge 0}R\overbrace{V^{\ast}\cdots V^{\ast}}^i
%                          \overbrace{U^{\ast}\cdots U^{\ast}}^j. 
\]
This implies that the $R$-algebra homomorphism
$R\langle\alpha,\beta\rangle\to H^{\ast}(N, T)$
induces an isomorphism 
$H^{\ast}(N, T)\cong R\langle\alpha,\beta\rangle/(\alpha\beta)$
of graded associative algebras over $R$. 
\qed

\bigskip 

We set $M={\rm M}_3(R)$.
Let us calculate the Hochschild cohomology
$H^{\ast}(N, M/N)$ of $N$ with coefficients in $M/N$.
For this purpose, we construct a spectral sequence
which converges to $H^{\ast}(N, M/N)$.
We show that the spectral sequence
collapses at the $E_2$-page
and there is no extension problem.
%See, for example, \cite[\S2.2]{McCleary} for construction
%of spectral sequences.

In order to construct the spectral sequence, 
we introduce a filtration on $M/N$.
We set $F^0=M/N$.
Let $L$ be the $R$-submodule of $M={\rm M}_3(R)$ 
consisting of matrices in which the $(3,1)$-entry is $0$.
We set $F^1=L/N$ and $F^2=B/N$,
where $B={\rm B}_3(R) = \{ (a_{ij} ) \in {\rm M}_3(R) \mid a_{ij}=0 \mbox{ for } i>j  \}$.
We have obtained a filtration
\[ 0=F^3\subset F^2\subset F^1\subset F^0=M/N\]
of $N$-bimodules over $R$.
We denote by
${\rm Gr}^p(M/N)$
the $p$-th associated graded module 
$F^p/F^{p+1}$. 
%The filtration induces a long exact sequence
%\[ \cdots \to H^{p+q}(N, F^{p+1})\to
%   H^{p+q}(N, F^p)\to H^{p+q}(N, {\rm Gr}^p(M/N))
%   \to H^{p+q+1}(N, F^{p+1})\to\cdots. \]
%
%We set 
%\[ \begin{array}{rcl}
%     D^{p,q}&=&H^{p+q}(N, F^p),\\
%     E^{p,q}&=&H^{p+q}(N, {\rm Gr}^p(M/N)).\\
%   \end{array} \]
%There is an exact couple
%\[
%\xymatrix{
%D\ar[rr]^{i}&\ar@{}[d]|{}&D\ar[dl]^{j}\\
%&E\ar[ul]^{k}&
%}
%\]
%where $D=\oplus_{p,q}D^{p,q}$ and $E=\oplus_{p,q}E^{p,q}$.
%By standard construction
%(cf.~\cite[\S2.2]{McCleary}), 
By Proposition \ref{prop:spectralsequence}, 
we obtain a spectral sequence
\[ E_1^{p,q}=H^{p+q}(N, {\rm Gr}^p(M/N))\Longrightarrow
   H^{p+q}(N, M/N)\]
with
\[ d_r: E_r^{p,q}\longrightarrow E_r^{p+r,q-r+1} \]
for $r\ge 1$.
Note that 
$E_1^{p,q}=0$ unless $0\le p\le 2$ and 
$p+q\ge 0$.
Thus,
the spectral sequence collapses at the $E_3$-page. 
Since $H^{\ast}(N, T)\cong R\langle\alpha,\beta\rangle/(\alpha\beta)$
and the $N$-bimodule ${\rm Gr}^p(M/N)$
is isomorphic to the direct sum of finitely many
copies of $T$,
we obtain that
\[ E_1^{p,\ast-p}\cong
   R\langle\alpha,\beta\rangle/(\alpha\beta)\otimes_R
   {\rm Gr}^p(M/N).\]

First, we calculate
$d_1:E_1^{0,q}\to E_1^{1,q}$ for $q\ge 0$.
We have
$E_1^{0,\ast}=H^{\ast}(N, M/L)$
which is isomorphic to
$R\langle\alpha,\beta\rangle/(\alpha\beta)\otimes_RRE_{31}$,
and
$E_1^{1,\ast-1}=H^{\ast}(N, L/B)$
which is isomorphic to 
$R\langle\alpha,\beta\rangle/(\alpha\beta)\otimes_R(RE_{21} \oplus RE_{32})$. 
We set 
\[ c(i,j)=\overbrace{\beta\cdots\beta}^i
   \overbrace{\alpha\cdots\alpha}^j
   \in H^{i+j}(N, T) \]
for $i,j\ge 0$. 
The set $\{c(i,j)|\  i, j \ge 0, \; i+j=n\}$
forms a basis of the free $R$-module 
$H^n(N, T)$ for all $n\ge 0$.
Since the differential $d_1: E_1^{0,\ast}\to E_1^{1,\ast}$
can be identified with the connecting homomorphism
$\delta:
H^{\ast}(N, M/L)\to H^{\ast+1}(N, L/B)$,
we obtain
\[ d_1(c(i,j)\otimes E_{31})=
   c(i+1,j)\otimes E_{21}+ (-1)^{i+j+1}c(i,j+1)\otimes E_{32} . 
\]   
%Thus, $E_2^{0,p}=0$ for all $p\ge 0$. 

Next, we calculate
$d_1:E_1^{1,q-1}\to E_1^{2,q-1}$ for $q\ge 0$.
We have $E_1^{1,\ast-1}=H^{\ast}(N, L/B)$,
which is isomorphic to 
$R\langle\alpha,\beta\rangle/(\alpha\beta)\otimes_R(RE_{21} \oplus RE_{32})$.
We have
$E_1^{2,\ast-2}=H^{\ast}(N, B/N)$
which is isomorphic to
$R\langle\alpha,\beta\rangle/(\alpha\beta) \otimes_RF^2$.
Since the differential $d_1:E_1^{1,q-1}\to E_1^{2,q-1}$
can be identified with 
the connecting homomorphism
$\delta: H^q(N, L/B)\to H^{q+1}(N, B/N)$,
we obtain
\[ \begin{array}{ccc}
    d_1(c(i,j)\otimes E_{21})&=&
    \left\{\begin{array}{ll}
            (-1)^{i+j+1}c(i,j+1)\otimes E_{22}&(i>0),\\
            c(0,j+1)\otimes (E_{11}+(-1)^{j+1}E_{22})&(i=0),\\
           \end{array}\right.\\[3mm]
    d_1(c(i,j)\otimes E_{32})&=&
    \left\{\begin{array}{ll}
            c(i+1,j)\otimes E_{22}&(j>0),\\
            c(i+1,0)\otimes (E_{22}+(-1)^{i+1}E_{33})&(j=0).\\
           \end{array}\right.
   \end{array}\]

By the above calculation of $d_1$,
$E_2^{0,q}=E_2^{1,q-1}=0$ and 
$E_2^{2,q-1}$ is a free $R$-module of rank $q+2$ 
for all $q\ge 0$.
Hence the spectral sequence
collapses at the $E_2$-page.
Since $E_{\infty}^{p,q}$ is a free $R$-module 
for all $p,q$, 
there is no extension problem.
Hence we obtain the following theorem.

\begin{theorem}
The $R$-module $H^n(N, M/N)$ is free for all $n\ge 0$. 
The rank of $H^n(N, M/N)$ over $R$ is given by
\[ {\rm rank}_R\,H^n(N, M/N)=
   \left\{\begin{array}{ll}
           2 & (n=0),\\[2mm]
           n+1 & (n\ge 1).\\
          \end{array}\right. \]
\end{theorem}

\subsection{The case $A = {\rm S}_{1}(k)$}\label{subsection:S1}
Let $A = R[x]/(x^2) \cong 
S_{1}(R)  = \left\{ 
\begin{array}{c|c} 
\left(
\begin{array}{ccc}
a & b & 0 \\
0 & a & 0 \\
0 & 0 & a \\ 
\end{array}
\right) 
& a, b \in R  
\end{array} \right\}$. Here  
$R$ is a commutative ring and $x$ corresponds to  $E_{12} \in {\rm S}_1(R)$.   
By Proposition \ref{prop:projresonegen}, there exists a projective resolution of $A$ as $A^{e}$-modules: 
\[ 
\cdots \to A^{e} \stackrel{d_n}{\to} A^{e} \to \cdots \to A^{e} \stackrel{d_1}{\to} A^{e} \stackrel{\mu}{\to} A \to 0,  
\]
where 
\[ d_i(a\otimes b) = \left\{ 
\begin{array}{cc}
(1\otimes x + x \otimes 1)  (a\otimes b) & (i : \mbox{ even}) \\
(1\otimes x - x \otimes 1) (a\otimes b) & (i : \mbox{ odd}) \\
\end{array} 
\right. 
\]  
and $\mu(a\otimes b) = ab$. 
Set $M = {\rm M}_3(R)/{\rm S}_1(R)$. 
By applying ${\rm Hom}_{A^{e}}(-, M)$ to the projective resolution above, 
we have 
\[
0 \to {\rm Hom}_{A^{e}}(A^{e}, M) \stackrel{d_1^{\ast}}{\to} {\rm Hom}_{A^{e}}(A^{e}, M) \to \cdots 
\to {\rm Hom}_{A^{e}}(A^{e}, M) \stackrel{d_n^{\ast}}{\to} {\rm Hom}_{A^{e}}(A^{e}, M) \to \cdots,   
\]  
which is isomorphic to 
\[ 
0 \to M \stackrel{b^1}{\to} M \to \cdots \to M \stackrel{b^n}{\to} M \to \cdots,  
\] 
where 
\[ b^i(m) = \left\{ 
\begin{array}{cc}
mx + xm & (i : \mbox{ even}) \\
mx - xm  & (i : \mbox{ odd}). \\
\end{array} 
\right. 
\]  
We can choose a basis $\{ E_{13}, E_{21}, E_{22}, E_{23}, E_{31}, E_{32}, E_{33} \}$ 
of the $R$-free module $M$. 
With respect to this basis, we have 
\[
b^i = \left(
\begin{array}{ccccccc}
0 & 0 & 0 & -1 & 0 & 0 & 0 \\ 
0 & 0 & 0 & 0 & 0 & 0 & 0 \\ 
0 & 2 & 0 & 0 & 0 & 0 & 0 \\ 
0 & 0 & 0 & 0 & 0 & 0 & 0 \\ 
0 & 0 & 0 & 0 & 0 & 0 & 0 \\ 
0 & 0 & 0 & 0 & 1 & 0 & 0 \\ 
0 & 1 & 0 & 0 & 0 & 0 & 0 \\ 
\end{array} 
\right) (i : \mbox{ odd}) \mbox{ and }
b^i = \left(
\begin{array}{ccccccc}
0 & 0 & 0 & 1 & 0 & 0 & 0 \\ 
0 & 0 & 0 & 0 & 0 & 0 & 0 \\ 
0 & 0 & 0 & 0 & 0 & 0 & 0 \\ 
0 & 0 & 0 & 0 & 0 & 0 & 0 \\ 
0 & 0 & 0 & 0 & 0 & 0 & 0 \\ 
0 & 0 & 0 & 0 & 1 & 0 & 0 \\ 
0 & -1 & 0 & 0 & 0 & 0 & 0 \\ 
\end{array} 
\right) (i : \mbox{ even}).  
\]
Thereby, we easily verify that 
\[
H^{i}({\rm S}_1(R), {\rm M}_3(k)/{\rm S}_1(R)) \cong \left\{  
\begin{array}{ll} 
RE_{13}\oplus RE_{22} \oplus RE_{32} \oplus RE_{33} \cong R^4 & (i=0) \\
(RE_{22}\oplus RE_{33})/R(2E_{22}+E_{33}) \cong R & (i \mbox{ : odd } > 0 ) \\ 
RE_{22} \cong R & (i \mbox{ : even } > 0 ). 
\end{array}
\right. 
\] 

\subsection{The case $A = {\rm S}_{2}(k)$}\label{subsection:S2}
Let us consider the quiver $Q = \quad 1 \longleftarrow 2$.  
Let $\Lambda = RQ/I$ be the incidence algebra associated to $Q$ over a commutative ring $R$. 
Then $\Lambda = Re_{11} \oplus R e_{22} \oplus Re_{12}$. 
We can regard $\Lambda$ as ${\rm S}_2(R) = \left\{ 
\begin{array}{c|c} 
\left(
\begin{array}{ccc}
a & 0 & 0 \\
0 & a & c \\
0 & 0 & b \\ 
\end{array}
\right) 
& a, b, c \in R 
\end{array} \right\}$ by $e_{11} \mapsto E_{11}+E_{22}, e_{22} \mapsto E_{33}$, and $e_{12} \mapsto E_{23}$.  
Set $M = {\rm M}_3(R)/{\rm S}_2(R)$. Then $M$ is a $\Lambda$-bimodule by identifying $\Lambda$ with 
the subalgebra ${\rm S}_2(R)$ of ${\rm M}_3(R)$. 
The free $R$-module $M$ has a basis $\{ E_{11}, E_{12}, E_{13}, E_{21}, E_{31}, E_{32} \}$. 
By Proposition \ref{prop:cibils}, it suffices to calculate the cohomology of 
the complex $\{ {\rm Hom}_{E^e}(r^{\otimes n}, M), \delta^n \}$, where 
$E = Re_{11}\oplus Re_{22} = R(E_{11}+E_{22})\oplus RE_{33}$ and $r = Re_{12} = RE_{23}$. 
Since $r^{\otimes n} = 0$ for $n \ge 2$, the complex is isomorphic to $0 \to M^{E} \stackrel{\delta^0}{\to} {\rm Hom}_{E^e}(r, M) \to 0 \to 0 \to \cdots$. It is easy to see that $M^{E} = RE_{11}\oplus RE_{12} \oplus RE_{21}$ and that 
${\rm Hom}_{E^e}(r, M) \cong RE_{13}$ by ${\rm Hom}_{E^e}(r, M) \ni f \mapsto f(E_{23}) \in RE_{13}$. 
By direct calculation, $\delta^{0}(E_{11})=\delta^{0}(E_{21})=0$ and $\delta^{0}(E_{12})=-E_{13}$.  
Hence we have 
\[
H^{i}({\rm S}_2(R), {\rm M}_3(R)/{\rm S}_2(R)) \cong \left\{  
\begin{array}{cl} 
RE_{11}\oplus RE_{21} \cong R^2 & (i=0) \\
0 & (i>0). 
\end{array}
\right. 
\]

\subsection{The case $A = {\rm S}_{4}(k)$}\label{subsection:S4}
Set ${\rm S}_4(R) = \left\{ 
\begin{array}{c|c} 
\left(
\begin{array}{ccc}
a & b & c \\
0 & a & 0 \\
0 & 0 & a \\ 
\end{array}
\right) 
& a, b, c \in R 
\end{array} \right\}$ for a commutative ring $R$.  
We set $A={\rm S}_4(R)$ and $M={\rm M}_3(R)$. 
In this subsection we calculate the Hochschild cohomology
$H^{\ast}(A, M/A)$ of $A$ with coefficients in $M/A$. 
For this purpose, we construct a spectral sequence
which converges to $H^{\ast}(A, M/A)$. 
We show that the spectral sequence
collapses at the $E_2$-page
and there is no extension problem.

Let $J$ be the two-sided ideal of $A$ given by 
$J = \left\{ 
\begin{array}{c|c} 
\left(
\begin{array}{ccc}
0 & b & c \\
0 & 0 & 0 \\
0 & 0 & 0 \\ 
\end{array}
\right) \in A
& b, c \in R 
\end{array} \right\}$. 
We set $T=A/J$,
which is an $A$-bimodule over $R$.

First, we describe
the Hochschild cohomology
$H^{\ast}(A, T)$ of $A$
with coefficients in $T$. 
We have an isomorphism
$T\otimes_A T\cong T$ of $A$-bimodules over $R$.
This implies that
$T$ is a monoid object in the category of $A$-bimodules
over $R$.
The unit $u: A\to T=A/J$ is given by the projection.
Hence $H^\ast(A, T)$ has the structure of
a graded associative algebra over $R$.

%We recall that $E_{ij}\in {\rm M}_3(R)$ is
%the matrix with entry $1$ in the $(i,j)$-component and
%$0$ the other components.
We set $I=E_{11}+E_{22}+E_{33}$,
$U=E_{12}$ and $V=E_{13}$.
The set $\{U,V\}$ forms a basis of the free $R$-module 
$\overline{A}=A/RI$.
Let $e\in T$ be the image of $I$ under the unit $u: A\to T$.
We denote by $U^{\ast},V^{\ast}\in \overline{C}^1(A, T)$ 
the maps $\overline{A} \to T$
of $R$-modules given by
\[ \begin{array}{rcl}
   U^{\ast}(n)&=&
   \left\{\begin{array}{ll}
     e & \mbox{\rm if $n=U$}\\
     0 & \mbox{\rm if $n=V$}.
     \end{array}\right.\\[3mm]
   V^{\ast}(n)&=& 
   \left\{\begin{array}{ll}
     e & \mbox{\rm if $n=V$}\\
     0 & \mbox{\rm if $n=U$}.
     \end{array}\right.
   \end{array}\]
We see that $U^{\ast}$ and $V^{\ast}$
are $1$-cocycles in the cochain complex
$\overline{C}^{\ast}(A, T)$.
We denote by
\[ \alpha,\beta\in H^1(A,T)\] 
the elements represented by the $1$-cocycles $U^{\ast},V^{\ast}$,
respectively.
It is easy to calculate the cohomology
$H^{\ast}(A,T)$ and we obtain the following lemma.

\begin{lemma}
There is an isomorphism
$H^{\ast}(A,T)\cong R\langle\alpha,\beta\rangle$
of graded associative algebras over $R$,
where $R\langle\alpha,\beta\rangle$
is the free graded associative algebra over $R$
generated by $\alpha$ and $\beta$.
\end{lemma}

\proof
The lemma follows from the observation that
$\overline{C}^{\ast}(A, T)$ is isomorphic to
a differential graded algebra which is the free graded associative
$R$-algebra $R\langle U^{\ast},V^{\ast}\rangle$ generated by 
$U^{\ast},V^{\ast}$ with trivial differential.
\qed

\bigskip 

In order to construct a spectral sequence, 
we introduce a filtration on $M/A$.
We set $F^0=M/A$.
Let $L$ be the $R$-submodule of $M={\rm M}_3(R)$ 
consisting of matrices in which the $(3,1)$-entry is $0$.
We set $F^1=L/A$ and $F^2=B/A$,
where $B={\rm B}_3(R) = \{ (a_{ij} ) \in {\rm M}_3(R) \mid a_{ij}=0 \mbox{ for } i>j  \}$.
We have obtained a filtration
\[ 0=F^3\subset F^2\subset F^1\subset F^0=M/A\]
of $A$-bimodules.
We denote by
${\rm Gr}^p(M/A)$
the $p$-th associated graded module
$F^p/F^{p+1}$.
%The filtration induces a long exact sequence
%\[ \cdots \to H^{p+q}(A, F^{p+1})\to
%   H^{p+q}(A, F^p)\to H^{p+q}(A, {\rm Gr}^p(M/A))
%   \to H^{p+q+1}(A, F^{p+1})\to\cdots. \]
%We set 
%\[ \begin{array}{rcl}
%     D^{p,q}&=&H^{p+q}(A, F^p)\\
%     E^{p,q}&=&H^{p+q}(A, {\rm Gr}^p(M/A))\\
%   \end{array} \]
%We obtain an exact couple
%\[
%\xymatrix{
%D\ar[rr]^{i}&\ar@{}[d]|{}&D\ar[dl]^{j}\\
%&E\ar[ul]^{k}&
%}
%\]
%where $D=\oplus_{p,q}D^{p,q}$ and $E=\oplus_{p,q}E^{p,q}$.
%By standard construction,
By Proposition \ref{prop:spectralsequence}, 
we obtain a spectral sequence
\[ E_1^{p,q}=H^{p+q}(A, {\rm Gr}^p(M/A))\Longrightarrow
   H^{p+q}(A, M/A)\]
with
\[ d_r: E_r^{p,q}\longrightarrow E_r^{p+r,q-r+1} \]
for $r\ge 1$.
Note that 
$E_1^{p,q}=0$ unless $0\le p\le 2$ and 
$p+q\ge 0$.
Thus,
the spectral sequence collapses at the $E_3$-page. 

The $A$-bimodule ${\rm Gr}^p(M/A)$
is isomorphic to the direct sum of finitely many
copies of $T$.
Since
$H^\ast(A, T)\cong R\langle\alpha,\beta\rangle$,
we obtain that
\[ E_1^{p,\ast-p}\cong
   R\langle\alpha,\beta\rangle\otimes_R
   {\rm Gr}^p(M/A).\]

First, we calculate
$d_1:E_1^{0,q}\to E_1^{1,q}$ for $q\ge 0$.
We have
$E_1^{0,\ast}=H^{\ast}(A, M/L)$
which is isomorphic to
$R\langle\alpha,\beta\rangle\otimes_RRE_{31}$,
and
$E_1^{1,\ast-1}=H^{\ast}(A, L/B)$
which is isomorphic to 
$R\langle\alpha,\beta\rangle\otimes_R (RE_{21} \oplus RE_{32})$. 
Let $\delta$ be the connecting homomorphism
$H^{q}(A, M/L)\to H^{q+1}(A, L/B)$.
We can identify $d_1: E_1^{0,q}\to E_1^{1,q}$
with $\delta$ 
and we obtain that
\begin{equation}\label{eq:differentiald1-(0,q)to(1,q)} 
   d_1(z\otimes E_{31})=
    (-1)^{q+1}z\alpha\otimes E_{32},
\end{equation}   
where $z\in H^q(A, T)$
is a monomial of $\alpha$ and $\beta$.

Next, we calculate
$d_1:E_1^{1,q-1}\to E_1^{2,q-1}$ for $q\ge 0$.
We have $E_1^{1,\ast-1}=H^{\ast}(A, L/B)$,
which is isomorphic to 
$R\langle\alpha,\beta\rangle\otimes_R(RE_{21} \oplus RE_{32})$.
We have
$E_1^{2,\ast-2}=H^{\ast}(A, B/A)$
which is isomorphic to
$R\langle\alpha,\beta\rangle\otimes_RF^2$.
Since the differential $d_1:E_1^{1,q-1}\to E_1^{2,q-1}$
can be identified with 
the connecting homomorphism
$\delta: H^q(A, L/B)\to H^{q+1}(A, B/A)$,
we obtain that
\begin{equation}\label{eq:differentiald1-(1,q-1)to(2,q-1)} 
  \begin{array}{rcl}
    d_1(z\otimes E_{21})&=&
    \alpha z\otimes E_{11}+(-1)^{q+1}z\alpha\otimes E_{22}+
    (-1)^{q+1}z\beta \otimes E_{23},\\
    d_1(z\otimes E_{32})&=& 0, \\
   \end{array}
\end{equation}
where $z\in H^q(A, T)$ is a monomial
of $\alpha$ and $\beta$.

By (\ref{eq:differentiald1-(0,q)to(1,q)}),
$E_2^{0,q}=0$ for all $q$.
Furthermore, 
by (\ref{eq:differentiald1-(0,q)to(1,q)}) 
and (\ref{eq:differentiald1-(1,q-1)to(2,q-1)}), 
we see that $E_2^{1,q-1}$ and 
$E_2^{2,q-2}$ are free $R$-modules for all $q\ge 0$. 
Thus, the spectral sequence
collapses at the $E_2$-page
and there is no extension problem.

\begin{theorem}
The $R$-module $H^n(A, M/A)$ is free for all $n\ge 0$. 
The rank of $H^n(A, M/A)$ over $R$ is given by
\[ {\rm rank}_R\, H^n(A, M/A)=
   \left\{\begin{array}{ll}
           4 & (n=0),\\
           3\cdot 2^n & (n\ge 1).\\
          \end{array}\right. \]
\end{theorem}

\proof

We know that the spectral sequence
collapses at the $E_2$-page and
that $E_2^{p,q}$ is a free $R$-module for all $p$ and $q$.
In what follows we calculate the rank of $E_2^{p,q}$ over $R$.

Since $E_1^{2,q-2}$ is isomorphic
to $H^q(A, B/A)$,
it is a free $R$-module of rank $3\cdot 2^q$
for all $q\ge 0$.
In particular, since $E_2^{2,-2}\cong E_1^{2,-2}$,
we obtain that 
\[ {\rm rank}_R\, E_2^{2,-2}=3. \]
By (\ref{eq:differentiald1-(1,q-1)to(2,q-1)}),
the image of $d_1^{1,q-1}: E_1^{1,q-1}\to E_1^{2,q-1}$
is a direct summand of $E_1^{2,q-1}$ and 
a free $R$-module of rank $2^q$ for all $q\ge 0$.
Since ${\rm rank}_R\,E_2^{2,q-1}=
{\rm rank}_R\,E_1^{2,q-1}-{\rm rank}_R\,{\rm Im}\,d_1^{1,q-1}$,
we obtain that
\[ {\rm rank}_R\, E_2^{2,q-1}=5\cdot 2^q\]
for all $q\ge 0$.

Since $E_1^{1,q-1}$ is isomorphic to
$H^q(A, L/B)$, it is a free $R$-module
of rank $2^{q+1}$ for all $q\ge 0$.
From the fact that 
${\rm rank}_R\,{\rm Im}\,d_1^{1,q-1}=2^q$,
we see that ${\rm rank}_R\,{\rm Ker}\,d_1^{1,q-1}=2^q$
for all $q\ge 0$. 
Since $E_1^{0,q}$ is isomorphic to
$H^q(A, M/L)$, 
we have ${\rm rank}_R\, E_1^{0,q}=2^q$
for all $q\ge 0$.
By (\ref{eq:differentiald1-(0,q)to(1,q)}),
we see that 
${\rm rank}_R\,{\rm Im}\,d_1^{0,q}=2^q$
for all $q\ge 0$.
Since
${\rm rank}_R\,E_2^{1,q-1}=
{\rm rank}_R\,{\rm Ker}\,d_1^{1,q-1}-
{\rm rank}_R\,{\rm Im}\,d_1^{0,q-1}$,
we obtain that 
\[ {\rm rank}_R\,E_2^{1,q-1}=\left\{
    \begin{array}{ll}
      2^{q-1}& (q>0),\\
      1     & (q=0).
    \end{array}\right.\]

The theorem follows from the fact
that 
\[ {\rm rank}_R\,H^n(A, M/A)=
   {\rm rank}_R\,E_2^{1,n-1}+{\rm rank}_R\,E_2^{2,n-2}\]
for all $n\ge 0$.
\qed

\subsection{The case $A = {\rm S}_{6}(k)$}\label{subsection:S6}
Let us consider the quiver 
\[
Q = 
\begin{xy}
{  \ar @(lu,ld)@{{*}->}  _{\alpha} }, 
{  \ar @{-<} (4, 0) }, 
{  (4, 0) \ar  @{-{*}} (10, 0) ^{\beta}  }, 
{  \put(0, 5){$e_1$} }, 
{  \put(30, 5){$e_2$} } 
\end{xy} \hspace*{3ex} . 
\] 
Let $RQ$ be the path algebra of $Q$ over a commutative ring $R$. 
Set $I = \langle \alpha^2, \alpha\beta \rangle \subset RQ$. 
Then we can regard $\Lambda = RQ/I = Re_1 \oplus Re_2 \oplus R\alpha \oplus R\beta$ as $S_6(R)= \left\{ 
\begin{array}{c|c} 
\left(
\begin{array}{ccc}
a & c & d \\
0 & a & 0 \\
0 & 0 & b \\ 
\end{array}
\right) 
& a, b, c, d \in R 
\end{array} \right\}$ by $e_1 \mapsto E_{11}+E_{22}$, $e_2 \mapsto E_{33}$, $\alpha \mapsto E_{12}$, and $\beta \mapsto E_{13}$.  
Set $M = {\rm M}_3(R)/{\rm S}_6(R)$. Then $M$ is a $\Lambda$-bimodule by identifying $\Lambda$ with the subalgebra ${\rm S}_6(R)$ of ${\rm M}_3(R)$. 
The free $R$-module $M$ has a basis $\{ E_{21}, E_{22}, E_{23}, E_{31}, E_{32} \}$. 
By Proposition \ref{prop:cibils}, it suffices to calculate the cohomology of 
the complex $\{ {\rm Hom}_{E^e}(r^{\otimes n}, M), \delta^n \}$, where 
$E = Re_1 \oplus Re_2 = R(E_{11}+E_{22})\oplus RE_{33}$ and $r = R\alpha \oplus R\beta = RE_{12}\oplus RE_{13}$. 
Obviously, $M^{E} = RE_{21}\oplus RE_{22}$. 
Since $\beta\otimes \alpha = \beta e_2 \otimes \alpha = \beta \otimes e_2 \alpha = 0$ and $\beta\otimes \beta = \beta e_2\otimes \beta = \beta\otimes e_2\beta = 0$, 
$r^{\otimes n} = R \alpha^{\otimes n} \oplus R(\alpha^{\otimes (n-1)}\otimes \beta)$ 
and ${\rm rank}_{R} r^{\otimes n} = 2$ for 
$n \ge 1$. 
By using $\alpha^{\otimes n} = e_1\alpha^{\otimes n}  e_1$ and 
$\alpha^{\otimes (n-1)}\otimes \beta = e_1 (\alpha^{\otimes (n-1)}\otimes \beta) e_2$,  
we see that $f(\alpha^{\otimes n}) \in e_1Me_1 = RE_{21} \oplus RE_{22}$ 
and $f(\alpha^{\otimes (n-1)}\otimes \beta) \in e_1 M e_2 = RE_{23}$ for 
$f \in {\rm Hom}_{E^{e}}(r^{\otimes n}, M)$. 
Let $(\alpha^{\otimes n})^{\ast}, (\alpha^{\otimes (n-1)}\otimes \beta)^{\ast} 
\in {\rm Hom}_{R}(r^{\otimes n}, R)$ be the dual basis 
of  $\alpha^{\otimes n}, \alpha^{\otimes (n-1)}\otimes \beta \in r^{\otimes n}$. 
Then we can write ${\rm Hom}_{E^{e}}(r^{\otimes n}, M) = (R(\alpha^{\otimes n})^{\ast}\otimes E_{21}) \oplus (R(\alpha^{\otimes n})^{\ast}\otimes E_{22}) \oplus (R(\alpha^{\otimes (n-1)}\otimes \beta)^{\ast}\otimes E_{23})$. In particular, ${\rm rank}_{R} {\rm Hom}_{E^{e}}(r^{\otimes n}, M) = 3$ for $n \ge 1$. 

First, let us calculate $\delta^{0} : M^{E}= RE_{21}\oplus RE_{22} \to {\rm Hom}_{E^{e}}(r, M)$. 
By direct calculation, 
\begin{eqnarray*}
\delta^{0}(E_{21})(\alpha) & = & \alpha E_{21} - E_{21}\alpha = E_{12}E_{21} - E_{21}E_{12} = E_{11}-E_{22} \equiv -2E_{22} \\ 
\delta^{0}(E_{21})(\beta) & = & \beta E_{21} - E_{21}\beta = E_{13}E_{21} - E_{21}E_{13} = - E_{23} \\ 
\delta^{0}(E_{22})(\alpha) & = & \alpha E_{22} - E_{22}\alpha = E_{12}E_{22} - E_{22}E_{12} = E_{12} \equiv 0 \\ 
\delta^{0}(E_{22})(\beta) & = & \beta E_{22} - E_{22}\beta = E_{13}E_{22} - E_{22}E_{13} =  0.  
\end{eqnarray*} 
With respect to the bases $\{ E_{21}, E_{22} \}$ and 
$\{ \alpha^{\ast} \otimes E_{21}, \alpha^{\ast}\otimes E_{22}, \beta^{\ast}\otimes E_{23} \}$, 
$\delta^{0}$ can be described as 
\[
\delta^{0} = \left(
\begin{array}{cc}
0 & 0 \\
-2 & 0 \\
-1 & 0 \\
\end{array} 
\right). 
\] 
Hence $H^{0}({\rm S}_{6}(R), {\rm M}_3(R)/{\rm S}_{6}(R)) \cong R$. 

Next, let us calculate $\delta^{n} : {\rm Hom}_{E^{e}}(r^{\otimes n}, M) \to {\rm Hom}_{E^{e}}(r^{\otimes (n+1)}, M)$. By direct calculation, 
\begin{eqnarray*}
\delta^{n}((\alpha^{\otimes n})^{\ast}\otimes E_{21})(\alpha^{\otimes (n+1)}) & = & \alpha E_{21} +(-1)^{n+1}E_{21}\alpha  \equiv -(1+(-1)^n) E_{22} \\ 
\delta^{n}((\alpha^{\otimes n})^{\ast}\otimes E_{21})(\alpha^{\otimes n}\otimes \beta) & = & (-1)^{n+1}E_{21}\beta = (-1)^{n+1} E_{23} \\ 
\delta^{n}((\alpha^{\otimes n})^{\ast}\otimes E_{22})(\alpha^{\otimes (n+1)}) & = & \alpha E_{22} +(-1)^{n+1}E_{22}\alpha  = E_{12} \equiv 0 \\ 
\delta^{n}((\alpha^{\otimes n})^{\ast}\otimes E_{22})(\alpha^{\otimes n}\otimes \beta) & = & (-1)^{n+1}E_{22}\beta = 0 \\ 
\delta^{n}((\alpha^{\otimes (n-1)}\otimes \beta)^{\ast}\otimes E_{23})(\alpha^{\otimes (n+1)}) & = & 0 \\ 
\delta^{n}((\alpha^{\otimes (n-1)}\otimes \beta)^{\ast}\otimes E_{23})(\alpha^{\otimes n}\otimes \beta) & = & \alpha E_{23} = E_{13} \equiv 0. 
\end{eqnarray*} 
With respect to the bases $\{ (\alpha^{\otimes n})^{\ast}\otimes E_{21}, (\alpha^{\otimes n})^{\ast}\otimes E_{22}, (\alpha^{\otimes (n-1)}\otimes \beta)^{\ast}\otimes E_{23} \}$ and 
$\{ (\alpha^{\otimes (n+1)})^{\ast}\otimes E_{21}, (\alpha^{\otimes (n+1)})^{\ast}\otimes E_{22}, (\alpha^{\otimes n}\otimes \beta)^{\ast}\otimes E_{23} \}$, 
$\delta^{n}$ can be described as 
\[
\delta^{n} = \left(
\begin{array}{ccc}
0 & 0 & 0 \\
0 & 0 & 0 \\
1 & 0 & 0 \\
\end{array} 
\right) (n : \mbox{ odd } ) \mbox{ and } 
\delta^{n} = \left(
\begin{array}{ccc}
0 & 0 & 0 \\
-2 & 0 & 0 \\
-1 & 0 & 0 \\
\end{array} 
\right) (n : \mbox{ even } ). 
\] 

Finally, let us calculate $H^{n}({\rm S}_{6}(R), {\rm M}_3(R)/{\rm S}_{6}(R))$. 
It is easy to see that 
\[ 
H^{n}({\rm S}_{6}(R), {\rm M}_3(R)/{\rm S}_{6}(R)) \cong   
((R(\alpha^{\otimes n})^{\ast}\otimes E_{22})  \oplus (R (\alpha^{\otimes (n-1)}\otimes \beta)^{\ast}\otimes E_{23}))/
{\rm Im} \:\delta^{n-1} 
\] 
for any $n > 0$, where ${\rm Im} \:\delta^{n-1} =  R((1+(-1)^{n+1})(\alpha^{\otimes n})^{\ast}\otimes E_{22}+(-1)^{n+1}(\alpha^{\otimes (n-1)}\otimes \beta)^{\ast}\otimes E_{23})$.  Summarizing the results, we have $H^{n}({\rm S}_{6}(R), {\rm M}_3(R)/{\rm S}_{6}(R)) \cong R$ for $n \ge 0$. 

\subsection{The case $A = {\rm S}_{7}(k)$}\label{subsection:S7}
Let us consider the quiver 
\[
Q = 
\begin{xy}
{  (25, 0) \ar@/^10pt/^{\beta} @{{*}->}  (10, 0) }, 
{  (25, 0) \ar@/_10pt/_{\alpha} @{-{*}} (10, 0) },
{  (25, 0) \ar@/_10pt/ @{->} (10, 0) },
%{  \ar @(lu,ld)@{{*}->}  _{\alpha} }, 
%{  \ar @{-<} (4, 0) }, 
%{  (4, 0) \ar  @{-{*}} (10, 0) ^{\beta}  }, 
{  \put(25, -10){$e_1$} }, 
{  \put(70, -10){$e_2$} } 
\end{xy} \hspace*{3ex} . 
\] 
Let $RQ$ be the path algebra of $Q$ over a commutative ring $R$. 
Then we can regard $\Lambda = RQ = Re_1 \oplus Re_2 \oplus R\alpha \oplus R\beta$ as $S_7(R)= \left\{ 
\begin{array}{c|c} 
\left(
\begin{array}{ccc}
a & 0 & c \\
0 & a & d \\
0 & 0 & b \\ 
\end{array}
\right) 
& a, b, c, d \in R 
\end{array} \right\}$ by $e_1 \mapsto E_{11}+E_{22}$, $e_2 \mapsto E_{33}$, $\alpha \mapsto E_{13}$, and $\beta \mapsto E_{23}$.  
Set $M = {\rm M}_3(R)/{\rm S}_7(R)$. Then $M$ is a $\Lambda$-bimodule by identifying $\Lambda$ with 
the subalgebra ${\rm S}_7(R)$ of ${\rm M}_3(R)$. 
The free $R$-module $M$ has a basis $\{ E_{12}, E_{21}, E_{22}, E_{31}, E_{32} \}$. 
By Proposition \ref{prop:cibils}, it suffices to calculate the cohomology of 
the complex $\{ {\rm Hom}_{E^e}(r^{\otimes n}, M), \delta^n \}$, where 
$E = Re_{1}\oplus Re_{2} = R(E_{11}+E_{22})\oplus RE_{33}$ and $r = R\alpha \oplus R\beta = RE_{13}\oplus RE_{23}$.  
Since $r^{\otimes n} = 0$ for $n \ge 2$, the complex is isomorphic to $0 \to M^{E} \stackrel{\delta^0}{\to} {\rm Hom}_{E^e}(r, M) \to 0 \to 0 \to \cdots$. It is easy to see that $M^{E} = RE_{12}\oplus RE_{21} \oplus RE_{22}$ and 
${\rm Hom}_{E^e}(r, M)=0$ by $r = e_1r e_2$ and $e_1Me_2 = 0$. 
Hence we have 
\[
H^{i}({\rm S}_7(R), {\rm M}_3(R)/{\rm S}_7(R)) \cong \left\{  
\begin{array}{cl} 
RE_{12}\oplus RE_{21} \oplus RE_{22} \cong R^3 & (i=0) \\
0 & (i>0). 
\end{array}
\right. 
\]

\subsection{The case $A = {\rm S}_{10}(k)$}\label{subsection:S10}
Let us consider the quiver 
\[
Q = 
\begin{xy}
{  \ar @(lu,ld)@{{*}->}  _{\alpha} }, 
{  \ar @{-<} (4, 0) }, 
{  (4, 0) \ar  @{-{*}} (10, 0) ^{\beta}  }, 
{  \put(0, 5){$e_1$} }, 
{  \put(30, 5){$e_2$} } 
\end{xy} \hspace*{3ex} . 
\] 
Let $RQ$ be the path algebra of $Q$ over a commutative ring $R$. 
Set $I = \langle \alpha^2 \rangle \subset RQ$ and $\gamma = \alpha\beta$. 
Then we can regard $\Lambda = RQ/I = Re_1 \oplus Re_2 \oplus R\alpha \oplus R\beta \oplus R\gamma$ as $S_{10}(R) = \left\{ 
\begin{array}{c|c} 
\left(
\begin{array}{ccc}
a & b & c \\
0 & a & d \\
0 & 0 & e \\ 
\end{array}
\right) 
& a, b, c, d \in R 
\end{array} \right\}$ by $e_1 \mapsto E_{11}+E_{22}$, $e_2 \mapsto E_{33}$, $\alpha \mapsto E_{12}$, $\beta \mapsto E_{23}$, and $\gamma \mapsto E_{13}$.  
Set $M = {\rm M}_3(R)/{\rm S}_{10}(R)$. Then $M$ is a $\Lambda$-bimodule by identifying $\Lambda$ with 
the subalgebra ${\rm S}_{10}(R)$ of ${\rm M}_3(R)$. 
The free $R$-module $M$ has a basis $\{ E_{21}, E_{22}, E_{31}, E_{32} \}$. 
By Proposition \ref{prop:cibils}, it suffices to calculate the cohomology of 
the complex $\{ {\rm Hom}_{E^e}(r^{\otimes n}, M), \delta^n \}$, where 
$E = Re_{1}\oplus Re_{2} = R(E_{11}+E_{22})\oplus RE_{33}$ and $r = R\alpha \oplus R\beta \oplus R\gamma = RE_{12}\oplus RE_{23} \oplus RE_{13}$.  
Since $e_1r e_1 =R\alpha$ and $e_1r e_2=R\beta \oplus R\gamma$, 
$r\otimes_{E} r = (R\alpha^{\otimes 2}) \oplus (R\alpha\otimes \beta) \oplus (R\alpha \otimes \gamma)$. 
Similarly, $r^{\otimes n} = (R\alpha^{\otimes n}) \oplus (R \alpha^{\otimes (n-1)} \otimes \beta) \oplus (R \alpha^{\otimes (n-1)} \otimes \gamma)$ for $n \ge 2$. 
Let $\{ (\alpha^{\otimes n})^{\ast}, (\alpha^{\otimes (n-1)} \otimes \beta)^{\ast}, (\alpha^{\otimes (n-1)} \otimes \gamma)^{\ast} \} \subset {\rm Hom}_{R}(r^{\otimes n}, R)$ be the dual basis of $\{ \alpha^{\otimes n}, \alpha^{\otimes (n-1)} \otimes \beta, \alpha^{\otimes (n-1)} \otimes \gamma \} \subset r^{\otimes n}$. 
Note that $M^{E} = RE_{21} \oplus RE_{22}$, $e_1Me_1 = RE_{21} \oplus RE_{22}$, and $e_1Me_2 = 0$. 
It is easy to see that ${\rm Hom}_{E^{e}}(r^{\otimes n}, M) = (R(\alpha^{\otimes n})^{\ast}\otimes E_{21}) \oplus (R(\alpha^{\otimes n})^{\ast}\otimes E_{22})$ for $n \ge 1$. 

First, let us calculate $\delta^{0} : M^{E} = RE_{21} \oplus RE_{22} \to {\rm Hom}_{E^{e}}(r, M) = (R\alpha^{\ast}\otimes E_{21}) \oplus (R\alpha^{\ast}\otimes E_{22})$. 
By direct calculation, 
\begin{eqnarray*}
\delta^{0}(E_{21})(\alpha) & = & \alpha E_{21} - E_{21}\alpha = E_{12}E_{21} - E_{21}E_{12} = E_{11}-E_{22} \equiv -2E_{22} \\ 
\delta^{0}(E_{22})(\alpha) & = & \alpha E_{22} - E_{22}\alpha = E_{12}E_{22} - E_{22}E_{12} = E_{12} \equiv 0. 
\end{eqnarray*} 
With respect to the bases $\{ E_{21}, E_{22} \}$ and 
$\{ \alpha^{\ast} \otimes E_{21}, \alpha^{\ast}\otimes E_{22} \}$, 
$\delta^{0}$ can be described as 
\[
\delta^{0} = \left(
\begin{array}{cc}
0 & 0 \\
-2 & 0 \\
\end{array} 
\right). 
\] 
Hence $H^{0}({\rm S}_{10}(R), {\rm M}_3(R)/{\rm S}_{10}(R)) \cong R\oplus {\rm Ann}(2)$, 
where ${\rm Ann}(2) = \{ a \in R \mid 2a=0 \}$.  

Next, let us calculate $\delta^{n} : {\rm Hom}_{E^{e}}(r^{\otimes n}, M) \to {\rm Hom}_{E^{e}}(r^{\otimes (n+1)}, M)$ for $n \ge 1$. By direct calculation, 
\begin{eqnarray*}
\delta^{n}((\alpha^{\otimes n})^{\ast}\otimes E_{21})(\alpha^{\otimes (n+1)}) & = & \alpha E_{21} +(-1)^{n+1}E_{21}\alpha  \equiv -(1+(-1)^n) E_{22} \\ 
\delta^{n}((\alpha^{\otimes n})^{\ast}\otimes E_{22})(\alpha^{\otimes (n+1)}) & = & \alpha E_{22} +(-1)^{n+1}E_{22}\alpha  = E_{12} \equiv 0.  
\end{eqnarray*} 
With respect to the bases $\{ (\alpha^{\otimes n})^{\ast}\otimes E_{21}, (\alpha^{\otimes n})^{\ast}\otimes E_{22} \}$ and 
$\{ (\alpha^{\otimes (n+1)})^{\ast}\otimes E_{21}, (\alpha^{\otimes (n+1)})^{\ast}\otimes E_{22} \}$, 
$\delta^{n}$ can be described as 
\[
\delta^{n} = \left(
\begin{array}{cc}
0 & 0  \\
0 & 0  \\
\end{array} 
\right) (n : \mbox{ odd}) \mbox{ and } 
\delta^{n} = \left(
\begin{array}{cc}
0 & 0  \\
-2 & 0  \\
\end{array} 
\right) (n : \mbox{ even}). 
\] 

Finally, let us calculate $H^{n}({\rm S}_{10}(R), {\rm M}_3(R)/{\rm S}_{10}(R))$. 
For $n \ge 0$, we easily see that 
\[
H^{n}({\rm S}_{10}(R), {\rm M}_3(R)/{\rm S}_{10}(R)) \cong \left\{  
\begin{array}{cl} 
R\oplus {\rm Ann}(2) & (n : \mbox{ even}) \\
R \oplus (R/2R) & (n : \mbox{ odd}). 
\end{array}
\right. 
\] 

\subsection{The case $A = {\rm S}_{11}(k)$}\label{subsection:S11} 
In this subsection, we calculate $H^{i}({\rm S}_{11}(R), {\rm M}_3(R)/{\rm S}_{11}(R))$ for a commutative ring $R$. 
In the following long proof of Theorem~\ref{th:S11}, the Fibonacci numbers appear, which seems strange to us. 
For another proof using spectral sequence without the Fibonacci numbers, see \cite{Nakamoto-Torii:52nd}.  

\bigskip 

Let us consider the quiver 
\[
Q = 
\begin{xy}
{  (25, 0) \ar@/^10pt/^{\beta} @{<-{*}}  (10, 0) }, 
{  (25, 0) \ar@/_10pt/_{\alpha} @{{*}-} (10, 0) },
{  (25, 0) \ar@/_10pt/ @{->} (10, 0) },
{  \put(25, -10){$e_1$} }, 
{  \put(70, -10){$e_2$} } 
\end{xy} \hspace*{3ex} . 
\] 
Let $RQ$ be the path algebra of $Q$ over a commutative ring $R$. 
Set $I = \langle \beta\alpha \rangle \subset RQ$ and $\gamma = \alpha\beta$. 
Then we can regard $\Lambda = RQ/I = Re_1 \oplus Re_2 \oplus R\alpha \oplus R\beta \oplus R\gamma$ as $S_{11}(R)= \left\{ 
\begin{array}{c|c} 
\left(
\begin{array}{ccc}
a & b & c \\
0 & e & d \\
0 & 0 & a \\ 
\end{array}
\right) 
& a, b, c, d, e \in R 
\end{array} \right\}$ by $e_1 \mapsto E_{11}+E_{33}$, $e_2 \mapsto E_{22}$, $\alpha \mapsto E_{12}$, $\beta \mapsto E_{23}$, and $\gamma \mapsto E_{13}$.  
Set $M = {\rm M}_3(R)/{\rm S}_{11}(R)$. Then $M$ is a $\Lambda$-bimodule by identifying $\Lambda$ with 
the subalgebra ${\rm S}_{11}(R)$ of ${\rm M}_3(R)$. 
The free $R$-module $M$ has a basis $\{ E_{11}, E_{21}, E_{31}, E_{32} \}$. 
%By Proposition \ref{prop:cibils}, it suffices to calculate the cohomology of 
%the complex $\{ {\rm Hom}_{E^e}(r^{\otimes n}, M), \delta^n \}$, where 
Set $E = Re_1 \oplus Re_2 = R(E_{11}+E_{33}) \oplus RE_{22}$ and $r = R\alpha \oplus R\beta \oplus R\gamma = RE_{12}\oplus RE_{23} \oplus RE_{13}$. 
Let ${\rm B}_3(R) = \left\{ 
\left(
\begin{array}{ccc}
\ast & \ast & \ast \\
0 & \ast & \ast \\
0 & 0 & \ast \\ 
\end{array}
\right)  
\right\}$ and $M' = {\rm B}_3(R)/{\rm S}_{11}(R)$. 
Then $M'$ is an ${\rm S}_{11}(R)$-bimodule and there exists an exact sequence of 
${\rm S}_{11}(R)$-bimodules (that is, an exact sequence of ${\rm S}_{11}(R)^{e}$-modules): 
\begin{eqnarray}
0 \to M' \to M \to M''=R\overline{E}_{21} \oplus R\overline{E}_{31} \oplus R\overline{E}_{32} \to 0. 
\end{eqnarray} 
Let us define the ${\rm S}_{11}(R)^{e}$-submodules $M_{21}$ and $M_{32}$ of $M''$ by 
$M_{21} = R\overline{E}_{21}$ and 
$M_{32} = R\overline{E}_{32}$, respectively. 
Put $M_{31} = M''/(M_{21} \oplus M_{32}) = R\overline{E}_{31}$.  
Note that the ${\rm S}_{11}(R)^{e}$-module $M_{31}$ is isomorphic to 
$M' = R\overline{E}_{11}$. 

Let ${\mathcal M}_n = \{ x_1 \otimes x_2 \otimes \cdots \otimes x_n \in r^{\otimes n} \mid x_1 \otimes x_2 \otimes \cdots \otimes x_n \neq 0, \mbox{where } x_i = \alpha, 
\beta, \mbox{ or } \gamma \}$ be the set of non-zero monomials of length $n$ in $\alpha$, $\beta$ and $\gamma$. 
For $n \ge 1$, $r^{\otimes n} = \oplus_{m \in {\mathcal M}_n} R m$ and ${\mathcal M}_n$ is a basis of 
the free module $r^{\otimes n}$ over $R$, where the tensor products are over $E$. 
Because $\alpha \otimes \alpha = \beta \otimes \beta = \alpha \otimes \gamma = 0$ and so on,  
${\rm rank}_{R} r^{\otimes n} < 3^n$ for $n >1$.  
Let ${\mathcal M}_n^{\ast} = \{ m^{\ast} \in {\rm Hom}_{R}(r^{\otimes n}, R) 
\mid m \in {\mathcal M}_n \}$ be the dual basis of ${\mathcal M}_n$.

\bigskip 

Let us introduce the following lemmas. 

\begin{lemma}\label{lemma:S11M11} 
Let $M' = R\overline{E}_{11}$ be as above. Then 
\[
H^{n}({\rm S}_{11}(R), M') \cong \left\{ 
\begin{array}{cc} 
R & ( n = 0 )  \\
0 & ( n > 0 ). \\
\end{array} 
\right. 
\] 
\end{lemma} 

\proof 
Note that 
\[
\left(
\begin{array}{ccc}
a & b & c \\
0 & e & d \\
0 & 0 & a \\
\end{array} 
\right) \overline{E}_{11} 
= a \overline{E}_{11} 
= \overline{E}_{11}  \left(
\begin{array}{ccc}
a & b & c \\
0 & e & d \\
0 & 0 & a \\
\end{array} 
\right).  
\] 
It is easy to see that $M' = e_1M'e_1$ and that $M'^{E} = M'$. 
For $f \in {\rm Hom}_{E^{e}}(r, M')$, $f(\alpha) = f(\alpha e_2) = f(\alpha)e_2 = 0$ and 
$f(\beta) = f(e_2 \beta) = e_2 f(\beta) = 0$. Thus, we have an isomorphism 
${\rm Hom}_{E^{e}}(r, M') \stackrel{\cong}{\to} M' =  R\overline{E}_{11}$ by $f \mapsto f(\gamma)$. 
First, let us consider $\delta^{0} : M'^{E} =M' \to {\rm Hom}_{E^{e}}(r, M') \cong R\overline{E}_{11}$. 
Since 
\[
\delta^{0}(\overline{E}_{11})(\gamma) = \gamma \overline{E}_{11} - \overline{E}_{11} \gamma = 0, 
\] 
we have $\delta^{0} = 0$. Hence $H^{0}(S_{11}(R), M') \cong R$. 

Next, let us consider $\delta^{1} : {\rm Hom}_{E^{e}}(r, M') \cong R\overline{E}_{11} \to {\rm Hom}_{E^{e}}(r\otimes_{E} r, M')$. Let $(\alpha\otimes \beta)^{\ast}, \ldots, (\gamma\otimes \gamma)^{\ast} \in {\rm Hom}_{R}(r\otimes_{E}r , R)$ 
be the dual basis of the basis $\alpha\otimes \beta, \ldots, \gamma\otimes \gamma$ of  $r\otimes_{E}r$ over $R$.  
By using $M'=e_1M'e_1$, we see that ${\rm Hom}_{E^{e}}(r\otimes_{E} r, M') = R((\alpha\otimes \beta)^{\ast}\otimes \overline{E}_{11}) \oplus R((\gamma\otimes \gamma)^{\ast}\otimes \overline{E}_{11})$.  
In a similar way, we can write ${\rm Hom}_{E^{e}}(r, M') = R(\gamma^{\ast}\otimes \overline{E}_{11})$.  
Since 
\begin{eqnarray*}
\delta^{1}(\gamma^{\ast}\otimes \overline{E}_{11})(\alpha \otimes \beta) & = & -\overline{E}_{11} \\ 
\delta^{1}(\gamma^{\ast}\otimes \overline{E}_{11})(\gamma \otimes \gamma) & = & \gamma \overline{E}_{11} + \overline{E}_{11}\gamma = 0, 
\end{eqnarray*} 
$\delta^{1}(\gamma^{\ast}\otimes \overline{E}_{11}) = -(\alpha\otimes \beta)^{\ast}\otimes \overline{E}_{11}$ and 
${\rm Ker} \: \delta^{1} = 0$. Hence $H^{1}(S_{11}(R), M') = 0$. 

We claim that ${\rm Hom}_{E^{e}}(r^{\otimes n}, M')$ is a free $R$-module. 
Set $F_n = {\rm rank}_{R} {\rm Hom}_{E^{e}}(r^{\otimes n}, M')$ under this claim.  
Note that $F_1=1$ and $F_2=2$. Since ${\rm Hom}_{E^{e}}(E, M') = M'^{E} = M' = R\overline{E}_{11}$, set $F_0=1$. 
By $M'=e_1M'e_1$ and 
$r^{\otimes n} = (e_1r^{\otimes n} e_1) \oplus (e_1r^{\otimes n} e_2) \oplus (e_2r^{\otimes n} e_1)  \oplus (e_2r^{\otimes n} e_2)$, 
we have 
${\rm Hom}_{E^{e}}(r^{\otimes n}, M') =  {\rm Hom}_{E^{e}}(e_1r^{\otimes n}e_1, R\overline{E}_{11}) \cong {\rm Hom}_{R}(e_1r^{\otimes n}e_1, R)$.   
Set ${\mathcal L}_n = \{ x_1 \otimes x_2 \otimes \cdots \otimes x_n \in e_1r^{\otimes n}e_1 \mid x_1 \otimes x_2 \otimes \cdots \otimes x_n \neq 0, \mbox{where } x_i = \alpha, 
\beta, \mbox{ or } \gamma \} \subseteq {\mathcal M}_n$. 
Then ${\rm Hom}_{E^{e}}(r^{\otimes n}, M')$ is a free $R$-module of rank $F_n = \sharp {\mathcal L}_n$. 
Indeed, ${\rm Hom}_{E^{e}}(r^{\otimes n}, M')$ has an $R$-basis 
${\mathcal L}_n^{\ast} \otimes \overline{E}_{11} = \{ m^{\ast} \otimes \overline{E}_{11} \mid m \in {\mathcal L}_n \}$,  where ${\mathcal L}_n^{\ast} = \{ m^{\ast} \in {\rm Hom}_{R}(e_1r^{\otimes n}e_1, R) \mid m \in {\mathcal L}_n \}$ is the dual basis of ${\mathcal L}_n \subset e_1r^{\otimes n}e_1$. 

Note that 
${\mathcal L}_1 = \{ \gamma \}$ and  ${\mathcal L}_2 = \{ \alpha\otimes \beta, \gamma\otimes \gamma \}$. 
For $n \ge 3$, ${\mathcal L}_n = (\gamma \otimes {\mathcal L}_{n-1}) \cup (\alpha\otimes \beta\otimes {\mathcal L}_{n-2})$, where $\gamma \otimes {\mathcal L}_{n-1}= \{ \gamma \otimes m \mid m \in {\mathcal L}_{n-1} \}$ and $\alpha\otimes \beta\otimes {\mathcal L}_{n-2} = \{ \alpha\otimes \beta \otimes m' \mid m' 
\in {\mathcal L}_{n-2} \}$. 
Thus, we have $F_{n} = F_{n-1} + F_{n-2}$. 
We call a monomial  in $\alpha\otimes \beta\otimes {\mathcal L}_{n-2}$ and 
in $\gamma \otimes {\mathcal L}_{n-1}$ 
type I and type II, respectively. 
Let us consider the lexicographic order on both ${\mathcal M}_n$ and ${\mathcal L}_n$ such that 
$\alpha > \beta > \gamma$.   
For example, $\alpha\otimes \beta \otimes \gamma > \gamma\otimes \alpha \otimes \beta > \gamma\otimes \gamma \otimes \gamma$ with respect to the lexicographic order on ${\mathcal L}_3 = \{ \alpha\otimes \beta \otimes \gamma, \gamma\otimes \alpha \otimes \beta, \gamma\otimes \gamma \otimes \gamma \}$. 
If $m > m'$ in ${\mathcal M}_n$, then $\alpha\otimes m>\alpha\otimes m'$, $\beta\otimes m>\beta\otimes m'$, 
and $\gamma\otimes m>\gamma\otimes m'$ unless they are zero. 
We also define the lexicographic order on ${\mathcal L}_n^{\ast} \otimes \overline{E}_{11}$ 
such that $m^{\ast}\otimes \overline{E}_{11} > m'^{\ast}\otimes \overline{E}_{11}$ if and only if 
$m > m'$ in ${\mathcal L}_n$. 
If $m \in {\mathcal L}_n$ is of type I or II, then we call $m^{\ast} \otimes \overline{E}_{11}$ 
type I or II, respectively.  
When $m_1$ and $m_2$ in ${\mathcal L}_n$ are of type I and II, respectively, 
$m_1^{\ast} \otimes \overline{E}_{11} > m_2^{\ast} \otimes \overline{E}_{11}$. 

For $n\ge 1$, let us describe 
$\delta^{n} : {\rm Hom}_{E^{e}}(r^{\otimes n}, M') \to {\rm Hom}_{E^{e}}(r^{\otimes (n+1)}, M')$ 
with respect to the ordered bases ${\mathcal L}_n^{\ast}\otimes \overline{E}_{11}$ and 
${\mathcal L}_{n+1}^{\ast}\otimes \overline{E}_{11}$. 
For $n=1$, $\delta^{1}(\gamma^{\ast}\otimes \overline{E}_{11}) = 
-(\alpha\otimes \beta)^{\ast}\otimes \overline{E}_{11}$ and 
\[
\delta^{1} = 
\left(
\begin{array}{c}
-1  \\
0 
\end{array} 
\right) 
\] 
with respect to ${\mathcal L}_{1}^{\ast}\otimes \overline{E}_{11} = \{ \gamma^{\ast} \otimes \overline{E}_{11} \}$ and 
${\mathcal L}_2^{\ast}\otimes \overline{E}_{11} = \{ (\alpha\otimes \beta)^{\ast}\otimes \overline{E}_{11}, (\gamma\otimes \gamma)^{\ast}\otimes \overline{E}_{11} \}$. 
For $n=2$, we have 
\begin{eqnarray*}
\delta^{2}((\alpha\otimes \beta)^{\ast}\otimes \overline{E}_{11})(\alpha\otimes \beta \otimes \gamma) & = & 
-\overline{E}_{11}\gamma = 0   \\ 
\delta^{2}((\alpha\otimes \beta)^{\ast}\otimes \overline{E}_{11})(\gamma \otimes \alpha\otimes \beta ) & = & 
\gamma \overline{E}_{11} = 0   \\
\delta^{2}((\alpha\otimes \beta)^{\ast}\otimes \overline{E}_{11})(\gamma\otimes \gamma \otimes \gamma) & = & 
0   \\
\delta^{2}((\gamma\otimes \gamma)^{\ast}\otimes \overline{E}_{11})(\alpha\otimes \beta \otimes \gamma) & = & 
-\overline{E}_{11}   \\ 
\delta^{2}((\gamma\otimes \gamma)^{\ast}\otimes \overline{E}_{11})(\gamma \otimes \alpha\otimes \beta ) & = & 
\overline{E}_{11}   \\
\delta^{2}((\gamma\otimes \gamma)^{\ast}\otimes \overline{E}_{11})(\gamma\otimes \gamma \otimes \gamma) & = & \gamma \overline{E}_{11} - \overline{E}_{11}\gamma = 0.    
\end{eqnarray*} 
Hence we can write 
\[
\delta^{2} = 
\left(
\begin{array}{cc}
0 & -1  \\
0 & 1  \\
0 & 0   
\end{array} 
\right) 
\] 
with respect to ${\mathcal L}_2^{\ast}\otimes \overline{E}_{11} = \{ (\alpha\otimes \beta)^{\ast}\otimes \overline{E}_{11}, (\gamma\otimes \gamma)^{\ast}\otimes \overline{E}_{11} \}$ and 
${\mathcal L}_{3}^{\ast}\otimes \overline{E}_{11} = \{ (\alpha\otimes \beta\otimes \gamma)^{\ast}\otimes \overline{E}_{11}, (\gamma\otimes \alpha\otimes \beta)^{\ast}\otimes \overline{E}_{11}, (\gamma\otimes \gamma\otimes \gamma)^{\ast}\otimes \overline{E}_{11} \}$. We easily see that $H^{2}(S_{11}(R), M') = 0$. 

For $n\ge 3$, we put   
\[
\delta^{n} = 
\left(
\begin{array}{cc}
A_n & B_n  \\
C_n & D_n  \\
\end{array} 
\right) 
\] 
with respect to ${\mathcal L}_{n}^{\ast}\otimes \overline{E}_{11} = \{ \mbox{ type I } \}\cup \{ \mbox{ type II } \}$ and  
${\mathcal L}_{n+1}^{\ast}\otimes \overline{E}_{11} = \{ \mbox{ type I } \}\cup \{ \mbox{ type II } \}$, where  
$A_n, B_n, C_n, D_n$ are matrices of size $F_{n-1}\times F_{n-2}$, $F_{n-1}\times F_{n-1}$, $F_{n}\times F_{n-2}$, and 
$F_{n}\times F_{n-1}$, respectively. For $m \in {\mathcal L}_{n-2}$, $m', m'' \in {\mathcal L}_{n-1}$ and $l \in {\mathcal L}_{n}$, we have 
\begin{eqnarray*}
\delta^{n}((\alpha\otimes \beta \otimes m)^{\ast}\otimes \overline{E}_{11})(\alpha\otimes \beta \otimes m') & = & \delta^{n-2}(m^{\ast}\otimes \overline{E}_{11})(m') \\ 
\delta^{n}((\alpha\otimes \beta \otimes m)^{\ast}\otimes \overline{E}_{11})(\gamma \otimes l) & = & 0 \\ 
\delta^{n}((\gamma \otimes m')^{\ast}\otimes \overline{E}_{11})(\alpha\otimes \beta \otimes m'') & = & 
\left\{ 
\begin{array}{cc}
-\overline{E}_{11} & (m'=m'') \\
0 & (m'\neq m'') \\ 
\end{array} 
\right. \\ 
\delta^{n}((\gamma \otimes m')^{\ast}\otimes \overline{E}_{11})(\gamma \otimes l) & = & -\delta^{n-1}(m'^{\ast}\otimes \overline{E}_{11})(l). 
\end{eqnarray*} 
Thus, we obtain $A_n = \delta^{n-2}$, $B_n= -I_{F_{n-1}}$, $C_n = 0$, and $D_n = -\delta^{n-1}$.  
Hence 
\[
\delta^{n} = 
\left(
\begin{array}{cc}
\delta^{n-2} & -I_{F_{n-1}}  \\
0 & -\delta^{n-1}  \\
\end{array} 
\right).  
\] 
Multiplying $\delta^{n}$ by invertible matrices, we have 
\[
\left(
\begin{array}{cc}
I_{F_{n-1}} & 0 \\
-\delta^{n-1} & I_{F_{n}} \\ 
\end{array}
\right) 
\delta^{n} 
\left(
\begin{array}{cc}
0 & I_{F_{n-2}}  \\
-I_{F_{n-1}} & \delta^{n-2}  \\ 
\end{array}
\right) = 
\left(
\begin{array}{cc}
I_{F_{n-1}} & 0  \\
0 & 0  \\ 
\end{array}
\right). 
\]
Thereby, ${\rm Ker} \: \delta^{n}$ and ${\rm Im} \: \delta^{n}$ are $R$-free modules of rank $F_{n-2}$ and 
$F_{n-1}$, respectively. We also see that 
the induced surjection ${\rm Hom}_{E^{e}}(r^{\otimes n}, M')/{\rm Im} \: \delta^{n-1} \to {\rm Hom}_{E^{e}}(r^{\otimes n}, M')/{\rm Ker} \: \delta^{n}$ is an $R$-homomorphism of free $R$-modules of the same rank $F_{n-1}$.  
Hence, the surjection is an isomorphism and ${\rm Ker} \: \delta^n = {\rm Im} \: \delta^{n-1}$. 
Therefore, $H^{n}({\rm S}_{11}(R), M')=0$ for $n\ge 3$. This completes the proof. 
\qed 

\begin{lemma}\label{lemma:S11M21} 
Let $M_{21} = R\overline{E}_{21}$ be as above. Then 
\[
H^{n}({\rm S}_{11}(R), M_{21}) \cong \left\{ 
\begin{array}{cc} 
R & ( n = 1 )  \\
0 & ( n \neq 1 ). \\
\end{array} 
\right. 
\] 
\end{lemma} 

\proof 
Note that 
\[
\left(
\begin{array}{ccc}
a & b & c \\
0 & e & d \\
0 & 0 & a \\
\end{array} 
\right) \overline{E}_{21} 
= e \overline{E}_{21}, \quad %\mbox{ and }  
\overline{E}_{21}  \left(
\begin{array}{ccc}
a & b & c \\
0 & e & d \\
0 & 0 & a \\
\end{array} 
\right) = 
a \overline{E}_{21}. 
\] 
It is easy to see that $M_{21} = e_2M_{21}e_1$ and that $M_{21}^{E} = 0$. 
For $f \in {\rm Hom}_{E^{e}}(r, M_{21})$, $f(\alpha) = f(\alpha e_2) = f(\alpha)e_2 = 0$ and 
$f(\gamma) = f(e_1 \gamma) = e_1 f(\gamma) = 0$. Thus, we have an isomorphism 
${\rm Hom}_{E^{e}}(r, M_{21}) \stackrel{\cong}{\to} M_{21} =  
R\overline{E}_{21}$ by $f \mapsto f(\beta)$. 
First, let us consider $\delta^{0} : M_{21}^{E} =0 \to {\rm Hom}_{E^{e}}(r, M_{21}) 
\cong R\overline{E}_{21}$. 
Since $\delta^{0} = 0$, $H^{0}(S_{11}(R), M_{21}) = 0$. 

Next, let us consider $\delta^{1} : {\rm Hom}_{E^{e}}(r, M_{21}) \cong R\overline{E}_{21} \to {\rm Hom}_{E^{e}}(r\otimes_{E} r, M_{21})$. Let $(\alpha\otimes \beta)^{\ast}, \ldots, (\gamma\otimes \gamma)^{\ast} \in {\rm Hom}_{R}(r\otimes_{E}r , R)$ 
be the dual basis of the basis $\alpha\otimes \beta, \ldots, \gamma\otimes \gamma$ of  $r\otimes_{E}r$ over $R$.  
By using $M_{21}=e_2M_{21}e_1$, we see that ${\rm Hom}_{E^{e}}(r\otimes_{E} r, M_{21}) = R((\beta\otimes \gamma)^{\ast}\otimes \overline{E}_{21})$.  
In a similar way, we can write ${\rm Hom}_{E^{e}}(r, M_{21}) = 
R(\beta^{\ast}\otimes \overline{E}_{21})$.  
Since 
\begin{eqnarray*}
\delta^{1}(\beta^{\ast}\otimes \overline{E}_{21})(\beta \otimes \gamma) & = & 
\overline{E}_{21}\gamma = 0, 
\end{eqnarray*} 
$\delta^{1}=0$ and $H^{1}(S_{11}(R), M_{21}) = R\overline{E}_{21} \cong R$. 

We claim that ${\rm Hom}_{E^{e}}(r^{\otimes n}, M_{21})$ is a free $R$-module. 
Set $F'_n = {\rm rank}_{R} {\rm Hom}_{E^{e}}(r^{\otimes n}, M_{21})$ under this claim.  
Note that $F'_1=1$ and $F'_2=1$. Since ${\rm Hom}_{E^{e}}(E, M_{21}) = M_{21}^{E} = 0$, set $F'_0=0$. 
By $M_{21}=e_2M_{21}e_1$ and 
$r^{\otimes n} = (e_1r^{\otimes n} e_1) \oplus (e_1r^{\otimes n} e_2) \oplus (e_2r^{\otimes n} e_1)  \oplus (e_2r^{\otimes n} e_2)$, 
we have 
${\rm Hom}_{E^{e}}(r^{\otimes n}, M_{21}) =  {\rm Hom}_{E^{e}}(e_2r^{\otimes n}e_1, R\overline{E}_{21}) \cong {\rm Hom}_{R}(e_2r^{\otimes n}e_1, R)$.   
Set ${\mathcal L}'_n = \{ x_1 \otimes x_2 \otimes \cdots \otimes x_n \in e_2r^{\otimes n}e_1 \mid x_1 \otimes x_2 \otimes \cdots \otimes x_n \neq 0, \mbox{where } x_i = \alpha, 
\beta, \mbox{ or } \gamma \} \subseteq {\mathcal M}_n$. 
Then ${\rm Hom}_{E^{e}}(r^{\otimes n}, M_{21})$ is a free $R$-module of rank $F'_n = \sharp {\mathcal L}'_n$. 
Indeed, ${\rm Hom}_{E^{e}}(r^{\otimes n}, M_{21})$ has an $R$-basis 
${\mathcal L'}_n^{\ast} \otimes \overline{E}_{21} = \{ m^{\ast} \otimes \overline{E}_{21} \mid m \in {\mathcal L}'_n \}$,  where ${\mathcal L'}_n^{\ast} = \{ m^{\ast} \in {\rm Hom}_{R}(e_2r^{\otimes n}e_1, R) \mid m \in {\mathcal L}'_n \}$ is the dual basis of ${\mathcal L}'_n \subset e_2r^{\otimes n}e_1$. 

Note that 
${\mathcal L}'_1 = \{ \beta \}$, ${\mathcal L}'_2 = \{ \beta\otimes \gamma \}$, and ${\mathcal L}'_3 = \{ 
\beta\otimes \alpha \otimes \beta, \beta\otimes \gamma \otimes \gamma \}$. 
For $n \ge 3$, ${\mathcal L}'_n = ( {\mathcal L}'_{n-1}\otimes \gamma) \cup ({\mathcal L}'_{n-2}\otimes \alpha\otimes \beta)$, where ${\mathcal L}'_{n-1}\otimes \gamma = \{ m\otimes \gamma \mid m \in {\mathcal L}'_{n-1} \}$ and ${\mathcal L}'_{n-2}\otimes \alpha\otimes \beta = \{ m'\otimes \alpha \otimes \beta \mid m' \in {\mathcal L}'_{n-2} \}$. Thus, we have $F'_{n} = F'_{n-1} + F'_{n-2}$. 
We call a monomial  in ${\mathcal L}_{n-2} \otimes \alpha\otimes \beta$ and 
in ${\mathcal L}_{n-1}\otimes \gamma$ 
type I and type II, respectively. 
Set $\alpha > \beta > \gamma$. For $x_1\otimes x_2 \otimes \cdots \otimes x_n, y_1 \otimes y_2 \otimes \cdots \otimes y_n \in {\mathcal L}'_n$ or ${\mathcal M}_n$, we say that $x_1\otimes x_2 \otimes \cdots \otimes x_n > y_1 \otimes y_2 \otimes \cdots \otimes y_n$ if there exists $k$ such that $x_n = y_n, x_{n-1} = y_{n-1}, \ldots, x_{k+1} = y_{k+1}$, and $x_k > y_k$.    For example, $\beta\otimes \gamma \otimes \alpha \otimes \beta >  
\beta\otimes \alpha \otimes \beta \otimes \gamma > \beta\otimes \gamma \otimes \gamma \otimes \gamma$ 
on ${\mathcal L}'_4 = \{ \beta\otimes \gamma \otimes \alpha \otimes \beta, 
\beta\otimes \alpha \otimes \beta \otimes \gamma, \beta\otimes \gamma \otimes \gamma \otimes \gamma \}$. 
If $m > m'$ in ${\mathcal M}_n$, then $m\otimes \alpha >m' \otimes \alpha$, $m\otimes \beta >m' \otimes \beta$, 
and $m\otimes \gamma>m'\otimes \gamma$ unless they are zero. 
We define an order on ${\mathcal L'}_n^{\ast} \otimes \overline{E}_{21}$ 
such that $m^{\ast}\otimes \overline{E}_{21} > m'^{\ast}\otimes \overline{E}_{21}$ if and only if 
$m > m'$ in ${\mathcal L}'_n$. 
If $m \in {\mathcal L}'_n$ is of type I or II, then we call $m^{\ast} \otimes \overline{E}_{21}$ 
type I or II, respectively.  
When $m_1$ and $m_2$ in ${\mathcal L}'_n$ are of type I and II, respectively, 
$m_1^{\ast} \otimes \overline{E}_{21} > m_2^{\ast} \otimes \overline{E}_{21}$. 

For $n\ge 1$, let us describe 
$\delta^{n} : {\rm Hom}_{E^{e}}(r^{\otimes n}, M_{21}) \to {\rm Hom}_{E^{e}}(r^{\otimes (n+1)}, M_{21})$ 
with respect to the ordered bases ${\mathcal L'}_n^{\ast}\otimes \overline{E}_{21}$ and 
${\mathcal L'}_{n+1}^{\ast}\otimes \overline{E}_{21}$. 
For $n=1$, $\delta^{1} = (0)$ 
with respect to ${\mathcal L'}_{1}^{\ast}\otimes \overline{E}_{21} = \{ \beta^{\ast} \otimes \overline{E}_{21} \}$ and 
${\mathcal L'}_2^{\ast}\otimes \overline{E}_{21} = \{ (\beta\otimes \gamma)^{\ast}\otimes \overline{E}_{21} \}$. 
For $n=2$, we have 
\begin{eqnarray*}
\delta^{2}((\beta\otimes \gamma)^{\ast}\otimes \overline{E}_{21})(\beta\otimes \alpha \otimes \beta) & = & 
\overline{E}_{21}  \\ 
\delta^{2}((\beta\otimes \gamma)^{\ast}\otimes \overline{E}_{21})(\beta\otimes \gamma \otimes \gamma) & = & 
-\overline{E}_{21}\gamma = 0.   
\end{eqnarray*} 
We can write 
\[
\delta^{2} = 
\left(
\begin{array}{c}
1  \\
0  \\
\end{array} 
\right) 
\] 
with respect to ${\mathcal L'}_2^{\ast}\otimes \overline{E}_{21} = \{ (\beta\otimes \gamma)^{\ast}\otimes \overline{E}_{21} \}$ and 
${\mathcal L'}_{3}^{\ast}\otimes \overline{E}_{21} = \{ (\beta\otimes \alpha\otimes \beta)^{\ast}\otimes \overline{E}_{21}, (\beta\otimes \gamma\otimes \gamma)^{\ast}\otimes \overline{E}_{21} \}$. Hence $H^{2}(S_{11}(R), M_{21}) = 0$. 

For $n\ge 3$, we put   
\[
\delta^{n} = 
\left(
\begin{array}{cc}
A_n & B_n  \\
C_n & D_n  \\
\end{array} 
\right) 
\] 
with respect to ${\mathcal L'}_{n}^{\ast}\otimes \overline{E}_{21} = \{ \mbox{ type I } \}\cup \{ \mbox{ type II } \}$ and  
${\mathcal L'}_{n+1}^{\ast}\otimes \overline{E}_{21} = \{ \mbox{ type I } \}\cup \{ \mbox{ type II } \}$, where  
$A_n, B_n, C_n, D_n$ are matrices of size $F'_{n-1}\times F'_{n-2}$, $F'_{n-1}\times F'_{n-1}$, $F'_{n}\times F'_{n-2}$, and 
$F'_{n}\times F'_{n-1}$, respectively. For $m \in {\mathcal L}'_{n-2}$, $m', m'' \in {\mathcal L}'_{n-1}$ and $l \in {\mathcal L}'_{n}$, we have 
\begin{eqnarray*}
\delta^{n}((m\otimes \alpha\otimes \beta)^{\ast}\otimes \overline{E}_{21})(m'\otimes\alpha\otimes \beta) & = & \delta^{n-2}(m^{\ast}\otimes \overline{E}_{21})(m') \\ 
\delta^{n}((m\otimes \alpha\otimes \beta)^{\ast}\otimes \overline{E}_{21})(l\otimes\gamma) & = & 0 \\ 
\delta^{n}((m'\otimes \gamma)^{\ast}\otimes \overline{E}_{21})(m''\otimes\alpha\otimes \beta) & = & 
\left\{ 
\begin{array}{cc}
(-1)^{n}\overline{E}_{21} & (m'=m'') \\
0 & (m'\neq m'') \\ 
\end{array} 
\right. \\ 
\delta^{n}((m'\otimes \gamma)^{\ast}\otimes \overline{E}_{21})(l\otimes\gamma) & = & \delta^{n-1}(m'^{\ast}\otimes \overline{E}_{21})(l). 
\end{eqnarray*} 
Thus, we obtain $A_n = \delta^{n-2}$, $B_n= (-1)^n I_{F'_{n-1}}$, $C_n = 0$, and $D_n = \delta^{n-1}$.  
Hence 
\[
\delta^{n} = 
\left(
\begin{array}{cc}
\delta^{n-2} & (-1)^n I_{F'_{n-1}}  \\
0 & \delta^{n-1}  \\
\end{array} 
\right)  
\] 
for $n \ge 3$. 
Multiplying $\delta^{n}$ by invertible matrices, we have 
\[
\left(
\begin{array}{cc}
I_{F'_{n-1}} & 0 \\
(-1)^{n+1} \delta^{n-1} & I_{F'_{n}} \\ 
\end{array}
\right) 
\delta^{n} 
\left(
\begin{array}{cc}
0 & (-1)^{n+1} I_{F'_{n-2}}  \\
(-1)^{n} I_{F'_{n-1}} & \delta^{n-2}  \\ 
\end{array}
\right) = 
\left(
\begin{array}{cc}
I_{F'_{n-1}} & 0  \\
0 & 0  \\ 
\end{array}
\right). 
\]
Thereby, ${\rm Ker} \: \delta^{n}$ and ${\rm Im} \: \delta^{n}$ are $R$-free modules of rank $F'_{n-2}$ and 
$F'_{n-1}$, respectively. 
We also see that 
the induced surjection ${\rm Hom}_{E^{e}}(r^{\otimes n}, M_{21})/{\rm Im} \: \delta^{n-1} \to {\rm Hom}_{E^{e}}(r^{\otimes n}, M_{21})/{\rm Ker} \: \delta^{n}$ is an $R$-homomorphism of free $R$-modules of the same rank $F'_{n-1}$.  
Hence, the surjection is an isomorphism and ${\rm Ker} \: \delta^n = {\rm Im} \: \delta^{n-1}$. 
Hence $H^{n}({\rm S}_{11}(R), M_{21})=0$ for $n\ge 3$. This completes the proof. 
\qed 

\bigskip 

In the same way as Lemma \ref{lemma:S11M21}, we can prove the following lemma. 

\begin{lemma}\label{lemma:S11M32} 
Let $M_{32} = R\overline{E}_{32}$ be as above. Then 
\[
H^{n}({\rm S}_{11}(R), M_{32}) \cong \left\{ 
\begin{array}{cc} 
R & ( n = 1 )  \\
0 & ( n \neq 1 ). \\
\end{array} 
\right. 
\] 
\end{lemma} 

By Lemmas \ref{lemma:S11M21} and \ref{lemma:S11M32}, we have the following corollary. 

\begin{corollary}\label{cor:S11M''} 
Let $M''={\rm M}_3(R)/{\rm B}_3(R) = R\overline{E}_{21} \oplus R\overline{E}_{31} \oplus R\overline{E}_{32}$ be as above. Then  
\[
H^{n}({\rm S}_{11}(R), M'') \cong \left\{ 
\begin{array}{cc} 
R & ( n = 1 )  \\
0 & ( n \neq 1 ). \\
\end{array} 
\right. 
\] 
\end{corollary} 

\proof 
Let $r = R\alpha \oplus R\beta \oplus R\gamma$ be as above. 
Since $e_1M''e_1 = R\overline{E}_{31}$, $e_1M''e_2 = R\overline{E}_{32}$, and $e_2M''e_1 = R\overline{E}_{21}$, 
we have $M''^{E} = e_1M''e_1 = R\overline{E}_{31}$. 
On the other hand, there exists an isomorphism 
\[
\begin{array}{ccc}
{\rm Hom}_{E^{e}}(r, M'') & \stackrel{\cong}{\to} & R\overline{E}_{32} \oplus R\overline{E}_{21} \oplus R\overline{E}_{31} \\ 
f & \mapsto & (f(\alpha), f(\beta), f(\gamma)). 
\end{array} 
\]
Let us calculate $\delta^{0} : M''^{E} = R\overline{E}_{31} \to {\rm Hom}_{E^{e}}(r, M'') \cong 
R\overline{E}_{32} \oplus R\overline{E}_{21} \oplus R\overline{E}_{31}$. Since 
\begin{eqnarray*}
\delta^{0}(\overline{E}_{31})(\alpha) & = & \alpha \overline{E}_{31} - \overline{E}_{31}\alpha = - \overline{E}_{32} \\  
\delta^{0}(\overline{E}_{31})(\beta) & = & \beta \overline{E}_{31} - \overline{E}_{31}\beta = \overline{E}_{21} \\
\delta^{0}(\overline{E}_{31})(\gamma) & = & \gamma \overline{E}_{31} - \overline{E}_{31}\gamma = \overline{E}_{11} -\overline{E}_{33} \equiv 0,  
\end{eqnarray*} 
${\rm Ker} \: \delta^{0} = 0$. 
Hence $H^{0}(S_{11}(R), M'') = 0$. 

Similarly, we have an isomorphism 
\[
\begin{array}{ccc}
{\rm Hom}_{E^{e}}(r\otimes r, M'') & \stackrel{\cong}{\to} & R\overline{E}_{31}\oplus R\overline{E}_{21} \oplus R\overline{E}_{32} \oplus R\overline{E}_{31} \\ 
f & \mapsto & (f(\alpha\otimes \beta), f(\beta\otimes \gamma), f(\gamma\otimes \alpha), f(\gamma\otimes \gamma)).  
\end{array} 
\]
By calculating $\delta^{1} : {\rm Hom}_{E^{e}}(r, M'') =R(\alpha^{\ast}\otimes \overline{E}_{32}) \oplus R(\beta^{\ast}\otimes \overline{E}_{21}) \oplus R(\gamma^{\ast}\otimes \overline{E}_{31})  
\to {\rm Hom}_{E^{e}}(r\otimes r, M'')$, we have 
\begin{eqnarray*}
\delta^{1}(\alpha^{\ast}\otimes \overline{E}_{32})(\alpha\otimes \beta) & = & \overline{E}_{32}\beta = \overline{E}_{33} \equiv 0 \\ 
\delta^{1}(\alpha^{\ast}\otimes \overline{E}_{32})(\beta \otimes \gamma) & = & 0 \\  
\delta^{1}(\alpha^{\ast}\otimes \overline{E}_{32})(\gamma\otimes \alpha) & = & \gamma\overline{E}_{32} =  \overline{E}_{12} \equiv 0 \\  
\delta^{1}(\alpha^{\ast}\otimes \overline{E}_{32})(\gamma\otimes \gamma) & = & 0 \\  
\delta^{1}(\beta^{\ast}\otimes \overline{E}_{21})(\alpha\otimes \beta) & = & \alpha\overline{E}_{21} = \overline{E}_{11} \equiv 0 \\ 
\delta^{1}(\beta^{\ast}\otimes \overline{E}_{21})(\beta\otimes \gamma) & = & \overline{E}_{21}\gamma = \overline{E}_{23} \equiv 0 \\ 
\delta^{1}(\beta^{\ast}\otimes \overline{E}_{21})(\gamma\otimes \alpha) & = & 0 \\ 
\delta^{1}(\beta^{\ast}\otimes \overline{E}_{21})(\gamma\otimes \gamma) & = & 0 \\ 
\delta^{1}(\gamma^{\ast}\otimes \overline{E}_{31})(\alpha\otimes \beta) & = & - \overline{E}_{31}  \\ 
\delta^{1}(\gamma^{\ast}\otimes \overline{E}_{31})(\beta\otimes \gamma) & = & \beta\overline{E}_{31}=\overline{E}_{21}  \\ 
\delta^{1}(\gamma^{\ast}\otimes \overline{E}_{31})(\gamma\otimes \alpha) & = & \overline{E}_{31}\alpha = \overline{E}_{32}  \\ 
\delta^{1}(\gamma^{\ast}\otimes \overline{E}_{31})(\gamma\otimes \gamma) & = & \gamma \overline{E}_{31}+\overline{E}_{31}\gamma =\overline{E}_{11}+\overline{E}_{33} \equiv 0.  
\end{eqnarray*} 
Since ${\rm Ker}\; \delta^{1} = R(\alpha^{\ast}\otimes \overline{E}_{32}) \oplus R(\beta^{\ast}\otimes \overline{E}_{21})$ and ${\rm Im} \; \delta^{0} = R(-\alpha^{\ast}\otimes \overline{E}_{32}+\beta^{\ast}\otimes \overline{E}_{21})$, 
$H^{1}({\rm S}_{11}(R), M'') \cong R$. 

The short exact sequence $0 \to M_{21}\oplus M_{32} \to M'' \to M_{31} \to 0$ induces a long exact sequence 
\[
\cdots \to H^{n}({\rm S}_{11}(R), M_{21}\oplus M_{32}) \to H^{n}({\rm S}_{11}(R), M'') \to 
H^{n}({\rm S}_{11}(R), M_{31}) \to \cdots . 
\] 
For $n \ge 2$, we see that $H^{n}({\rm S}_{11}(R), M_{21}\oplus M_{32}) \cong H^{n}({\rm S}_{11}(R), M_{21})\oplus H^{n}({\rm S}_{11}(R), M_{32}) = 0$ and that $H^{n}({\rm S}_{11}(R), M_{31}) \cong H^{n}({\rm S}_{11}(R), M')=0$ by 
Lemmas \ref{lemma:S11M11}, \ref{lemma:S11M21}, and \ref{lemma:S11M32}. 
Hence $H^{n}({\rm S}_{11}(R), M'') = 0$ for $n\ge 2$. This completes the proof. 
\qed 

\bigskip 

By the discussions above, we have 

\begin{theorem}\label{th:S11} 
\[
H^{n}({\rm S}_{11}(R), {\rm M}_3(R)/{\rm S}_{11}(R)) \cong \left\{ 
\begin{array}{cc} 
R & ( n = 0, 1 )  \\
0 & ( n \ge 2 ). \\
\end{array} 
\right. 
\] 
\end{theorem} 

\proof 
The short exact sequence $0 \to M' \to M \to M'' \to 0$ induces a long exact sequence: 
\[
\begin{array}{rl}
0 & \to H^{0}({\rm S}_{11}(R), M') \to H^{0}({\rm S}_{11}(R), M) \to  H^{0}({\rm S}_{11}(R), M'')  \\
  & \to H^{1}({\rm S}_{11}(R), M') \to H^{1}({\rm S}_{11}(R), M) \to  H^{1}({\rm S}_{11}(R), M'')  \\ 
  & \to H^{2}({\rm S}_{11}(R), M') \to H^{2}({\rm S}_{11}(R), M) \to  H^{2}({\rm S}_{11}(R), M'')  \to \cdots . 
\end{array} 
\] 
Using Lemma \ref{lemma:S11M11} and Corollary \ref{cor:S11M''}, we have 
\[
\begin{array}{rl}
0 & \to R \to H^{0}({\rm S}_{11}(R), M) \to  0  \\
  & \to 0 \to H^{1}({\rm S}_{11}(R), M) \to R  \\ 
  & \to 0 \to H^{2}({\rm S}_{11}(R), M) \to 0  \to \cdots . 
\end{array} 
\] 
Thereby, $H^{0}({\rm S}_{11}(R), M) \cong H^{1}({\rm S}_{11}(R), M) \cong R$. 
By using Lemma \ref{lemma:S11M11} and Corollary \ref{cor:S11M''} again, 
$H^{n}({\rm S}_{11}(R), M') \cong H^{n}({\rm S}_{11}(R), M'') = 0$ for $n \ge 2$. 
Hence $H^{n}({\rm S}_{11}(R), M)=0$ for $n \ge 2$. 
\qed 

\subsection{The case $A = {\rm S}_{13}(k)$}\label{subsection:S13}
Let us consider the quiver 
\[
Q = 
\begin{xy}
{  (5, 5) \ar @{->} _{\alpha} (0, 0)   },  
{  (0, 0) \ar @{{*}-{*}} (5, 5)    },  
{  (5, -5) \ar @{{*}->} ^{\beta} (0, 0)    },  
{  \put(8, -3){$e_1$} }, 
{  \put(20, 13){$e_2$} }, 
{  \put(20, -18){$e_3$} } 
\end{xy} \hspace*{3ex} . 
\] 
Let $\Lambda$ be the incidence algebra associated to the ordered quiver $Q$ over a commutative ring $R$. 
Then we can regard $\Lambda \cong RQ = Re_1 \oplus Re_2 \oplus Re_3 \oplus R\alpha \oplus R\beta$ 
as $S_{13}(R) = \left\{  
\left(
\begin{array}{ccc}
\ast & \ast & \ast \\
0 & \ast & 0 \\
0 & 0 & \ast \\ 
\end{array}
\right)  
\right\}$ by $e_1 \mapsto E_{11}$, $e_2 \mapsto E_{22}$, $e_3 \mapsto E_{33}$, $\alpha \mapsto E_{12}$, and $\beta \mapsto E_{13}$.  
By Theorem \ref{th:incidence}, $H^{n}(S_{13}(R), {\rm M}_3(R)/S_{13}(R)) = 0$ for $n \ge 0$.

%\[
%Q = 
%\begin{xy}
%{  \ar @(lu,ld)@{{*}->}  _{\alpha} }, 
%{ \ar  @{-{*}} (10, 0) |{>} ^{\beta}  }, 
%{  \put(0, 5){$e_1$} }, 
%{  \put(30, 5){$e_2$} } 
%\end{xy} \hspace*{3ex} . 
%\] 
%
%\[
%Q = 
%\begin{xy}
%{  \ar @(lu,ld)@{{*}->}  _{\alpha} }, 
%{ \ar  !~:{@{-}@|{>}} (10, 0) ^{\beta}  }, 
%{  \put(0, 5){$e_1$} }, 
%{  \put(30, 5){$e_2$} } 
%\end{xy} \hspace*{3ex} . 
%\] 
%
%\[
%  \xymatrix{
%  \ar @(lu,ld)@{{*}->}  _{\alpha}  
%    A \ar[r] \ar[rrd] & B & \ar[d] \\
%    C \ar[u] & A \ar[l] \ar[lu] & D
%   }
%\]

\section{Appendix: Results on $H^{i}(A, {\rm M}_n(R)/A)$}

In this appendix, we show the tables on Hochschild cohomology $H^{\ast}(A, {\rm M}_2(R)/A)$ for $R$-subalgebras $A$ of ${\rm M}_n(R)$ over a commutative ring $R$ in the case $n=2, 3$.   
The ${}^t A$ column denotes the equivalence classes of ${}^t A$. The $N(A)$ column denotes the normalizer $N(A) = \{ b \in {\rm M}_n(R) \mid [b, a] = ba - ab \in A \mbox{ for any } a \in A \}$ of $A$.  
We also define ${\rm S}_{i}(R)$, ${\rm N}_3(R)$, ${\rm J}_3(R)$, etc. for a commutative ring $R$ in the same way as 
the case that $R$ is a field. 

\begin{table}[h]
\caption{Hochschild cohomology $H^{\ast}(A, {\rm M}_2(R)/A)$ for $R$-subalgebras $A$ of ${\rm M}_2(R)$}
\label{table:deg2} 
\tiny 
\begin{tabular}{|c|c|c|c|c|c|c|}\hline 
$A$ & $d = {\rm rank} A$ & $H^{\ast} = H^{\ast}(A, {\rm M}_2(R)/A)$ &  ${}^{t} A$ 
& $N(A)$ &  $\dim T_{{\rm Mold}_{2, d}/{\Bbb Z}, A}$  \\ 
\hline ${\rm M}_2(R)$ & $4$ & 
$H^{i} = 0$ for $i \ge 0$
& ${\rm M}_2(R)$ & ${\rm M}_2(R)$ & $0$   \\
\hline ${\rm B}_2(R) = \left\{ \left( 
\begin{array}{cc}
\ast & \ast   \\
0 & \ast 
\end{array} 
\right) \right\}$ & $3$ & $H^{i} = 0$ for $i \ge 0$ & ${\rm B}_2(R)$ & ${\rm B}_2(R)$ & $1$   \\
\hline ${\rm D}_2(R) = 
\left\{ \left( 
\begin{array}{cc}
\ast & 0   \\
0 & \ast 
\end{array} 
\right) \right\}$ & $2$ & $H^{i} = 0$ for $i \ge 0$ & ${\rm D}_2(R)$ 
& ${\rm D}_2(R)$ & $2$   \\
\hline ${\rm N}_{2}(R) = 
\left\{ \left( 
\begin{array}{cc}
a & b   \\
0 & a  
\end{array} 
\right) \right\}$ & $2$ & $H^{i} \cong \left\{ 
\begin{array}{cc}
R \oplus {\rm Ann}(2) & (i : \mbox{even}) \\ 
R \oplus (R/2R) & (i : \mbox{odd})
\end{array} 
\right.$ & ${\rm N}_{2}(R)$ & $\left\{ 
\begin{array}{c|c} 
\left( 
\begin{array}{cc}
\ast & \ast   \\
a & \ast  
\end{array} 
\right) & 2a = 0 
\end{array} \right\}$ & $2$   \\
\hline ${\rm C}_2(R) = \left\{ \left( 
\begin{array}{cc}
a & 0   \\
0 & a   
\end{array} 
\right) \right\}$ & $1$ & $H^{i} \cong \left\{ 
\begin{array}{cc}
R^3 & (i=0) \\
0 & (i \ge 1) 
\end{array} 
\right.$ & ${\rm C}_2(R)$ & ${\rm M}_2(R)$ & $0$  \\ \hline 
\end{tabular}
\end{table} 

%
%  Table deg3 
%

\begin{table}[hp]  
%\tiny 
\caption{Hochschild cohomology $H^{\ast}(A, {\rm M}_3(R)/A)$ for $R$-subalgebras $A$ of ${\rm M}_3(R)$}
\label{table:deg3}  
\scalebox{0.65}{
\begin{tabular}{|c|c|c|c|c|c|c|}\hline 
$A$ & $d={\rm rank} A$ & $H^{\ast} = H^{\ast}(A, {\rm M}_3(R)/A)$ &  ${}^{t} A$ 
& $N(A)$ &  $\dim T_{{\rm Mold}_{3, d}/{\Bbb Z}, A}$   \\ 
\hline ${\rm M}_3(R)$ & $9$ & 
$H^{i} = 0$ for $i \ge 0$
& ${\rm M}_3(R)$ & ${\rm M}_3(R)$ & $0$   \\
\hline ${\rm P}_{2, 1}(R) =  
\left\{ \left( 
\begin{array}{ccc}
\ast & \ast & \ast  \\
\ast & \ast & \ast  \\
0 & 0 & \ast 
\end{array} 
\right) \right\}$ & $7$ & $H^{i} = 0$ for $i \ge 0$ & ${\rm P}_{1, 2}(R)$ & 
${\rm P}_{2, 1}(R)$ & $2$   \\ 
\hline ${\rm P}_{1, 2}(R) = 
\left\{ \left( 
\begin{array}{ccc}
\ast & \ast & \ast  \\
0 & \ast & \ast  \\
0 & \ast & \ast 
\end{array} 
\right) \right\}$ & $7$ & $H^{i} = 0$ for $i \ge 0$ 
& ${\rm P}_{2, 1}(R)$ & ${\rm P}_{1, 2}(R)$ & $2$   \\ 
\hline ${\rm B}_3(R) = \left\{ \left( 
\begin{array}{ccc}
\ast & \ast & \ast  \\
0 & \ast & \ast  \\
0 & 0 & \ast 
\end{array} 
\right) \right\}$ & $6$ & $H^{i} = 0$ for $i \ge 0$ & ${\rm B}_3(R)$ & ${\rm B}_3(R)$ & $3$   \\
\hline $({\rm M}_2 \times {\rm D}_1)(R) = 
\left\{ \left( 
\begin{array}{ccc}
\ast & \ast & 0  \\
\ast & \ast & 0  \\
0 & 0 & \ast 
\end{array} 
\right) \right\}$ & $5$ & $H^{i} = 0$ for $i \ge 0$ & $({\rm M}_2 \times {\rm D}_1)(R)$ 
& $({\rm M}_2 \times {\rm D}_1)(R)$ & $4$   \\
\hline ${\rm S}_{10}(R) = 
\left\{ \left( 
\begin{array}{ccc}
a & b & c  \\
0 & a & d  \\
0 & 0 & e 
\end{array} 
\right) \right\}$ & $5$ & $H^{i} \cong \left\{ 
\begin{array}{cc}
R \oplus {\rm Ann}(2) & (i : \mbox{even}) \\ 
R \oplus (R/2R) & (i : \mbox{odd})
\end{array} 
\right.$ & ${\rm S}_{12}(R)$ & $\left\{ 
\begin{array}{c|c} 
\left( 
\begin{array}{ccc}
\ast & \ast & \ast  \\
a & \ast & \ast  \\
0 & 0 & \ast 
\end{array} 
\right) & 2a = 0 
\end{array} \right\}$ & $4$   \\
\hline ${\rm S}_{11}(R) = 
\left\{ \left( 
\begin{array}{ccc}
a & b & c  \\
0 & e & d  \\
0 & 0 & a 
\end{array} 
\right) \right\}$ & $5$ & $H^{i} \cong \left\{ 
\begin{array}{cc}
R & (i=0, 1) \\
0 & (i\ge 2) 
\end{array} 
\right.$ & ${\rm S}_{11}(R)$ & ${\rm B}_3(R)$ & $4$   \\
\hline ${\rm S}_{12}(R) = 
\left\{ \left( 
\begin{array}{ccc}
a & b & c  \\
0 & e & d  \\
0 & 0 & e 
\end{array} 
\right) \right\}$ & $5$ & $H^{i} \cong \left\{ 
\begin{array}{cc}
R \oplus {\rm Ann}(2) & (i : \mbox{even}) \\ 
R \oplus (R/2R) & (i : \mbox{odd})
\end{array} 
\right.$ & ${\rm S}_{10}(R)$ & $\left\{ 
\begin{array}{c|c} 
\left( 
\begin{array}{ccc}
\ast & \ast & \ast  \\
0 & \ast & \ast  \\
0 & a & \ast 
\end{array} 
\right) & 2a = 0 
\end{array} \right\}$ & $4$   \\
\hline ${\rm S}_{13}(R) = \left\{ \left( 
\begin{array}{ccc}
\ast & \ast & \ast  \\
0 & \ast & 0  \\
0 & 0 & \ast 
\end{array} 
\right) \right\}$ & $5$ & $H^{i} = 0$ for $i \ge 0$ & ${\rm S}_{14}(R)$ & ${\rm S}_{13}(R)$ & $4$   \\
\hline ${\rm S}_{14}(R) = \left\{ \left( 
\begin{array}{ccc}
\ast & 0 & \ast  \\
0 & \ast & \ast  \\
0 & 0 & \ast 
\end{array} 
\right) \right\}$ & $5$ & $H^{i} = 0$ for $i \ge 0$ & ${\rm S}_{13}(R)$ & ${\rm S}_{14}(R)$ & $4$   \\ 
\hline $({\rm B}_2 \times {\rm D}_1)(R) = \left\{ \left( 
\begin{array}{ccc}
\ast & \ast & 0  \\
0 & \ast & 0  \\
0 & 0 & \ast 
\end{array} 
\right) \right\}$ & $4$ &  $H^{i} = 0$ for $i \ge 0$ & $({\rm B}_2 \times {\rm D}_1)(R)$ 
& $({\rm B}_2 \times {\rm D}_1)(R)$ & $5$ \\
\hline ${\rm N}_3(R) = 
\left\{ \left( 
\begin{array}{ccc}
a & b & c  \\
0 & a & d  \\
0 & 0 & a 
\end{array} 
\right) \right\}$ & $4$ &  $H^{i} \cong \left\{ 
\begin{array}{cc}
R^2 & (i=0) \\
R^{i+1} & (i \ge 1) 
\end{array} 
\right.$ & ${\rm N}_3(R)$ & ${\rm B}_3(R)$ & $5$   \\ 
\hline ${\rm S}_6(R) = 
\left\{ \left( 
\begin{array}{ccc}
a & c & d  \\
0 & a & 0  \\
0 & 0 & b 
\end{array} 
\right) \right\}$ & $4$ & $H^{i} \cong R$ for $i \ge 0$ & ${\rm S}_{9}(R)$ & ${\rm S}_{13}(R)$ & $5$   \\
\hline ${\rm S}_7(R) = 
\left\{ \left( 
\begin{array}{ccc}
a & 0 & c  \\
0 & a & d  \\
0 & 0 & b 
\end{array} 
\right) \right\}$ & $4$ & $H^{i} \cong \left\{ 
\begin{array}{cc}
R^3 & (i=0) \\
0 & (i \ge 1) 
\end{array} 
\right.$ & ${\rm S}_{8}(R)$ & ${\rm P}_{2, 1}(R)$ & $2$   \\
\hline ${\rm S}_8(R) = \left\{ \left( 
\begin{array}{ccc}
a & c & d  \\
0 & b & 0  \\
0 & 0 & b 
\end{array} 
\right) \right\}$ & $4$ & $H^{i} \cong \left\{ 
\begin{array}{cc}
R^3 & (i=0) \\
0 & (i \ge 1) 
\end{array} 
\right.$ & ${\rm S}_{7}(R)$ & ${\rm P}_{1, 2}(R)$ & $2$   \\
\hline ${\rm S}_9(R) = \left\{ \left( 
\begin{array}{ccc}
a & 0 & c  \\
0 & b & d  \\
0 & 0 & b 
\end{array} 
\right) \right\}$ & $4$ & $H^{i} \cong R$ for $i \ge 0$ & ${\rm S}_{6}(R)$ & ${\rm S}_{14}(R)$ & $5$   \\
\hline ${\rm D}_3(R) = \left\{ \left( 
\begin{array}{ccc}
\ast & 0 & 0  \\
0 & \ast & 0  \\
0 & 0 & \ast 
\end{array} 
\right) \right\}$ & $3$ & $H^{i} = 0$ for $i \ge 0$ & ${\rm D}_{3}(R)$ & ${\rm D}_{3}(R)$ & $6$   \\
\hline $({\rm N}_2 \times {\rm D}_1)(R) = \left\{ \left( 
\begin{array}{ccc}
a & c & 0  \\
0 & a & 0  \\
0 & 0 & b 
\end{array} 
\right) \right\}$ & $3$ & $H^{i} \cong \left\{ 
\begin{array}{cc}
R \oplus {\rm Ann}(2) & (i : \mbox{even}) \\ 
R \oplus (R/2R) & (i : \mbox{odd})
\end{array} 
\right.$ & $({\rm N}_2 \times {\rm D}_1)(R)$ 
& $\left\{ 
\begin{array}{c|c} 
\left( 
\begin{array}{ccc}
\ast & \ast & 0  \\
a & \ast & 0  \\
0 & 0 & \ast 
\end{array} 
\right) & 
\begin{array}{c} 
2a = 0 
\end{array} 
\end{array} \right\}$ & $6$  \\ 
\hline ${\rm J}_3(R) = \left\{ \left( 
\begin{array}{ccc}
a & b & c  \\
0 & a & b  \\
0 & 0 & a 
\end{array} 
\right) \right\}$ & $3$ & $H^{i} \cong \left\{ 
\begin{array}{cc}
R^2 \oplus {\rm Ann}(3) & (i : \mbox{even}) \\ 
R^2 \oplus (R/3R) & (i : \mbox{odd})
\end{array} 
\right.$ & ${\rm J}_3(R)$ &  $\left\{ 
\begin{array}{c|c} 
\left( 
\begin{array}{ccc}
a & \ast & \ast  \\
c & a+b & \ast  \\
0 & -c & a+2b 
\end{array} 
\right) & 
\begin{array}{c} 
a, b, c \in R \\ 
3c = 0 
\end{array} 
\end{array} \right\}$ & $6$   \\
\hline ${\rm S}_2(R)  = \left\{ \left( 
\begin{array}{ccc}
a & 0 & 0  \\
0 & a & c  \\
0 & 0 & b 
\end{array} 
\right) \right\}$ & $3$ & $H^{i} \cong \left\{ 
\begin{array}{cc}
R^2 & (i=0) \\
0 & (i \ge 1) 
\end{array} 
\right.$ & ${\rm S}_{3}(R)$ &  
$\left\{ 
\left( 
\begin{array}{ccc}
\ast & 0 & 0  \\
\ast &\ast & \ast  \\
0 & 0 & \ast  
\end{array} 
\right) 
\right\} \sim {\rm S}_{13}(R)$ & $4$   \\
\hline ${\rm S}_3(R)  = \left\{ \left( 
\begin{array}{ccc}
a & 0 & c  \\
0 & b & 0  \\
0 & 0 & b 
\end{array} 
\right) \right\}$ & $3$ &  $H^{i} \cong \left\{ 
\begin{array}{cc}
R^2 & (i=0) \\
0 & (i \ge 1) 
\end{array} 
\right.$ & ${\rm S}_{2}(R)$ & ${\rm S}_{14}(R)$ & $4$   \\
\hline ${\rm S}_4(R) = \left\{ \left( 
\begin{array}{ccc}
a & b & c  \\
0 & a & 0  \\
0 & 0 & a 
\end{array} 
\right) \right\}$ & $3$ & $H^{i} \cong \left\{ 
\begin{array}{cc}
R^4 & (i=0) \\
R^{3\cdot 2^{i}} & (i \ge 1) 
\end{array} 
\right.$ & ${\rm S}_{5}(R)$ & ${\rm P}_{1, 2}(R)$ & $8$   \\
\hline ${\rm S}_5(R) = \left\{ \left( 
\begin{array}{ccc}
a & 0 & b  \\
0 & a & c  \\
0 & 0 & a 
\end{array} 
\right) \right\}$ & $3$ & $H^{i} \cong \left\{ 
\begin{array}{cc}
R^4 & (i=0) \\
R^{3\cdot 2^{i}} & (i \ge 1) 
\end{array} 
\right.$ & ${\rm S}_{4}(R)$ & ${\rm P}_{2, 1}(R)$ & $8$  \\
\hline $({\rm C}_2 \times {\rm D}_1)(R) = \left\{ \left( 
\begin{array}{ccc}
a & 0 & 0  \\
0 & a & 0  \\
0 & 0 & b 
\end{array} 
\right) \right\}$ & $2$ &  $H^{i} \cong \left\{ 
\begin{array}{cc}
R^3 & (i=0) \\
0 & (i \ge 1) 
\end{array} 
\right.$ & $({\rm C}_2 \times {\rm D}_1)(R)$  
& $({\rm M}_2\times {\rm D}_1)(R)$ & $4$  \\
\hline ${\rm S}_1(R) = \left\{ \left( 
\begin{array}{ccc}
a & b & 0  \\
0 & a & 0  \\
0 & 0 & a 
\end{array} 
\right) \right\}$ & $2$ & $H^{i} \cong \left\{ 
\begin{array}{cc}
R^4 & (i=0) \\
R & (i \ge 1) 
\end{array} 
\right.$ & ${\rm S}_{1}(R)$ & $\left\{ 
\left( 
\begin{array}{ccc}
\ast & \ast & \ast \\
0 &\ast & 0  \\
0 & \ast & \ast  
\end{array} 
\right) 
\right\} \sim {\rm B}_3(R)$ & $4$   \\
\hline ${\rm C}_3(R) = \left\{ \left( 
\begin{array}{ccc}
a & 0 & 0  \\
0 & a & 0  \\
0 & 0 & a 
\end{array} 
\right) \right\}$ & $1$ & $H^{i} \cong \left\{ 
\begin{array}{cc}
R^8 & (i=0) \\
0 & (i \ge 1) 
\end{array} 
\right.$ & ${\rm C}_3(R)$ & ${\rm M}_3(R)$ & $0$  \\ \hline 
\end{tabular}
} 
\end{table} 

%
%\newpage 
%

\end{document}